\documentclass[microtype]{gtpart}
\usepackage{amssymb, amsmath, amscd, latexsym, mathrsfs, appendix, amsxtra, color}
\usepackage{graphics}
\usepackage{amsfonts}
\usepackage[all]{xy}
\title[Rational Pontryagin classes]{Dalian notes on rational Pontryagin classes}

\author[Michael Weiss]{Michael S. Weiss}
\givenname{Michael}
\surname{Weiss}
\address{Mathematics Institute, WWU M\"unster, Einsteinstrasse 62, 48149 M\"unster, Germany}
\email{m.weiss@uni-muenster.de}


\theoremstyle{plain}
\newtheorem{thm}{Theorem}[subsection]
\newtheorem{prop}[thm]{Proposition}
\newtheorem{cor}[thm]{Corollary}
\newtheorem{lem}[thm]{Lemma}

\theoremstyle{definition}
\newtheorem{defn}[thm]{Definition}

\newtheorem{rem}[thm]{Remark}
\newtheorem{expl}[thm]{Example}

\newcommand{\lra}{\longrightarrow}

\newcommand{\smin}{\smallsetminus}

\newcommand{\cone}{\textup{cone}}
\newcommand{\id}{\textup{id}}
\newcommand{\im}{\textup{im}}

\newcommand{\holim}{\textup{holim}}
\newcommand{\hocolim}{\textup{hocolim}}
\newcommand{\colim}{\textup{colim}}
\newcommand{\hofiber}{\textup{hofiber}}
\newcommand{\mor}{\textup{mor}}

\newcommand{\map}{\textup{map}}
\newcommand{\rmap}{\textup{\textsf{R}}\textup{map}}

\newcommand{\str}{\textup{str}}

\newcommand{\fimm}{\textup{fimm}}

\newcommand{\TOP}{\textup{TOP}}
\newcommand{\STOP}{\textup{STOP}}

\newcommand{\GL}{\textup{GL}}
\newcommand{\PL}{\textup{PL}}

\newcommand{\SOr}{\textup{SO}}
\newcommand{\Or}{\textup{O}}

\newcommand{\diff}{\textup{diff}}
\newcommand{\emb}{\textup{emb}}
\newcommand{\loc}{\textup{loc}}
\newcommand{\imm}{\textup{imm}}

\newcommand{\parpun}{{g,1+\varepsilon}}
\newcommand{\inc}{\textup{in}}

\newcommand{\config}{\mathsf{con}}

\newcommand{\man}{\mathsf{Man}} 
\newcommand{\fin}{\mathsf{Fin}} 
\newcommand{\finplus}{{\fin_*}}
\newcommand{\rel}{\textup{ver}} 
\newcommand{\rhh}{\textrm{rhh}}

\newcommand{\intr}{\textup{int}}

\newcommand{\sA}{\mathsf{A}}
\newcommand{\sC}{\mathsf{C}}
\newcommand{\sD}{\mathsf{D}}
\newcommand{\sJ}{\mathcal J} 
\newcommand{\sE}{\mathcal E} 
\newcommand{\sH}{\mathcal H} 
\newcommand{\sK}{\mathcal K} 
\newcommand{\sL}{\mathcal L} 
\newcommand{\sP}{\mathcal P} 
\newcommand{\sQ}{\mathcal Q}

\newcommand{\op}{\textup{op}}
\newcommand{\haut}{\textup{haut}}
\newcommand{\homeo}{\textup{homeo}}
\newcommand{\tree}{\mathsf{Tree}}

\newcommand{\sSet}{\mathsf{sSet}}

\newcommand{\NN}{\mathbb N}
\newcommand{\RR}{\mathbb R}
\newcommand{\ZZ}{\mathbb Z}
\newcommand{\QQ}{\mathbb Q}
\newcommand{\bound}{\textup{bound}}
\newcommand{\cobound}{\textup{cobound}}

\newcommand{\twosub}[2]{\begin{array}{cc}
\scriptstyle{#1} \\  [-1mm] \scriptstyle{#2}  \end{array}}

\newcommand{\holimsub}[1]{\begin{array}[t]{cc} \textup{holim} \\ [-1mm]
\scriptstyle{#1} \end{array}}

\newcommand{\hocolimsub}[1]{\begin{array}[t]{cc} \textup{hocolim} \\
[-1mm] \scriptstyle{#1} \end{array}}

\newcommand{\uli}{\underline}

\begin{document}
\begin{abstract} The Pontryagin classes of an $m$-dimensional real vector bundle satisfy well-known vanishing relations
due to their close relationship with the Chern classes of the complexified vector bundle. It is also known that the
\emph{rational} Pontryagin classes qualify as characteristic classes for Euclidean fiber bundles, i.e.,
bundles whose fibers are homeomorphic to $m$-dimensional Euclidean space. In this more general setting of fiber bundles,
the vanishing relations which one has for vector bundles fail to hold.
\end{abstract}

\maketitle

\section{Introduction} \label{sec-intro}
\subsection{Main result} \label{subsec-mainresult}
The ring $H^*(B\Or;\QQ)$ is a polynomial ring $\QQ[p_1,p_2,p_3,\dots]$
where $p_i\in H^{4i}(B\Or;\QQ)$ is the Pontryagin class. The restrictions of the
classes $p_i$ to the classifying spaces for finite-dimensional (oriented) vector bundles satisfy well-known relations:
\begin{equation} \label{eqn-relations} p_n=e^2 \in H^{4n}(B\SOr(2n);\QQ), \quad \forall k>0:~p_{n+k}=0\in H^{4n+4k}(B\Or(2n);\QQ)
\end{equation}
where $e$ denotes the Euler class.

Let $\TOP(m)$ be the topological or simplicial group of homeomorphisms from $\RR^m$ to $\RR^m$.
The classifiying space $B\TOP(m)$ is also a classifying space for \emph{fiber bundles} with fiber $\RR^m$ (and
transition group $\TOP(m)$). Let
\[ B\TOP=\bigcup_{m\ge 0} B\TOP(m) \]
be the colimit of the spaces $B\TOP(m)$. In the late 1950s and early 1960s it emerged that the rational Pontryagin
classes of an $m$-dimensional vector bundle depend only on the underlying fiber bundle (with transition group
$\TOP(m)$); see \cite{RohlinSvarc}, \cite{Thom}, \cite{Novikov}. More to the point, thanks largely to the work of Sullivan \cite[Thm.18]{SullivanThesis} and Kirby-Siebenmann \cite{KirbySiebenmann} in the late 1960s it was understood that the inclusion $B\Or\to B\TOP$ induces an isomorphism in rational cohomology,
\[
\xymatrix{
H^*(B\TOP;\QQ) \ar[r]^-\cong & H^*(B\Or;\QQ).
}
\]
Therefore we can write unambiguously $p_i\in H^{4i}(B\TOP;\QQ)$. We can also write
unambiguously $p_i\in H^{4i}(B\TOP(m);\QQ)$ using the restriction map from
$H^*(B\TOP;\QQ)$ to $H^*(B\TOP(m);\QQ)$.
In other words, the rational Pontryagin classes can be
viewed as characteristic classes for fiber bundles with fiber $\RR^m$ for some $m$,
and as such they are stable by construction; they do not change under fiberwise product of such fiber bundles
with trivial line bundles.

\emph{From now on homology and cohomology are
taken with rational coeffficients unless otherwise stated.} ---
The main result of this article is that the analogues of relations \eqref{eqn-relations}
fail to hold in the cohomology or $B\TOP(2n)$ and $B\STOP(2n)$, respectively. (Here
$\STOP(m)$ stands for the topological group of orientation preserving homeomorphisms
from $\RR^m$ to $\RR^m$; a mild abuse of notation.) For large enough
$n$ we have $p_n\ne e^2$ in $H^{4n}(B\STOP(2n))$
and even more surprisingly, $p_{n+k}\ne 0$ in $H^{4n+4k}(B\TOP(2n))$
where $4n+4k$ can be nearly as big as $9n$. Here is a more precise statement.

\begin{thm} \label{thm-pmain} There exist positive constants $c_1$ and $c_2$ such that,
for all positive integers $n$ and $k$ where $n\ge c_1$ and $k<5n/4-c_2$, the class
$p_{n+k}$ is nonzero in $H^{4n+4k}(B\TOP(2n))$.
\end{thm}

\begin{cor} \label{cor-eulpon} Let $e\in H^{2n+2k}(B\STOP(2n+2k))$ be the Euler class.
For $n$ and $k$ satisfying the conditions of theorem~\ref{thm-pmain}, we have
\[ p_{n+k}\ne e^2 \in H^{4n+4k}(B\STOP(2n+2k)). \]
\end{cor}  

\proof By a well known transfer argument,
$p_{n+k}\ne 0$ in $H^{4n+4k}(B\TOP(2n))$ implies $p_{n+k}\ne 0$ in
$H^{4n+4k}(B\STOP(2n))$. It follows that $p_{n+k}\ne e^2$
in the cohomo\-logy ring $H^*(B\STOP(2n+2k))$. \qed

\begin{expl} It is allowed to take $c_1=83$ and $c_2=11$ in theorem~\ref{thm-pmain}.
Consequently for $k=1,2,\dots,92$ the classes $p_{83+k}\in H^{4(83+k)}(B\TOP(2\cdot83))$ are nonzero
and $p_{83+k}\ne e^2$ in $H^{4(83+k)}(B\STOP(2\cdot(83+k)))$.
But we will also discover that $p_{n+k}$ is nonzero in $H^{4n+4k}(B\TOP(2n))$ if $n\ge 53$ and $k\le \lfloor 3n/4\rfloor$.
Consequently for $k=1,2,\dots,39$ the classes $p_{53+k}\in H^{4(53+k)}(B\TOP(2\cdot 53))$ are nonzero
and $p_{53+k}\ne e^2$ in $H^{4\cdot(53+k)}(B\STOP(2\cdot(53+k)))$.

These bounds, e.g. $c_1=83$ and $c_2=11$, emerge in section~\ref{subsec-stabi}, more precisely in
lemmas~\ref{lem-appligoeppl} and~\ref{lem-estimate}. The least $n$ for which \emph{it is shown here} that
$p_{n+k}$ is nonzero in $H^{4n+4k}(B\TOP(2n))$ for some $k>0$ is $n=53$. But this is probably
far from the best possible; see remark~\ref{rem-KupKran}.
\end{expl}

\smallskip
There is a geometric formulation of theorem~\ref{thm-pmain} which is in fact slightly stronger.

\begin{thm} \label{thm-gmain} For $n$ and $k$ satisfying the conditions in theorem~\ref{thm-pmain}, there exists a
fiber bundle $E\to M$ where
\begin{itemize}
\item $M$ is a closed smooth
stably parallelized manifold of dimension $2n+4k$,
\item the fiber is a closed oriented topological manifold
of dimension $2n$,  
\item the signature of the total space $E$ is nonzero, but
\item all decomposable $($tangential$)$ Pontryagin numbers of $E$ are zero.
\end{itemize}
\end{thm}
The \emph{decomposable} Pontryagin numbers of $E$ are those corresponding to monomials in the Pontryagin
classes, of cohomological degree $4n+4k$ but distinct from $p_{n+k}$. There is no assumption in theorem~\ref{thm-gmain}
that $M$ be connected, but it would be alright to add the condition that all fibers of $E\to M$ are oriented homeomorphic to the same
oriented $2n$-dimensional manifold. The structure group of the bundle
is then the group of orientation-preserving homeomorphisms of that manifold.

In our construction of the fiber bundle $E\to M$ in theorem~\ref{thm-gmain}, the fiber is $W_g$~, the connected
sum of $g$ copies of $S^n\times S^n$, for some large $g$.   

\smallskip
\proof[Deducing theorem~\ref{thm-pmain} from theorem~\ref{thm-gmain}] Adopt the notation of theorem~\ref{thm-gmain}.
Tangent bundles of topological manifolds such as $E$ can be interpreted as tangent microbundles \cite{Milnormicro}, but see also \cite{Kister}.
Specifically, let $TE$ be the tangent (micro-)bundle of $E$.
The Hirzebruch signature theorem \cite{Hirz},\cite{MilStash} is usually stated for closed orientable \emph{smooth} manifolds, but it is also
valid for closed orientable topological manifolds. (The inclusion $B\SOr \to B\STOP$ is a rational homotopy equivalence and consequently
the inclusion of Thom spectra $M\SOr\to M\STOP$ is also a rational homotopy equivalence. By transversality,
these Thom spectra can be interpreted as bordism spectra, although it must be said that transversality
in the setting of topological manifolds is a difficult story. See for example \cite{QuinnTrans}.)
Therefore the scalar product of $p_{n+k}(TE)\in H^{4n+4k}(E;\QQ)$ with the fundamental
class of $E$ must be nonzero (since the other contributors in Hirzebruch's formula for the signature, the
decomposable tangential Pontryagin numbers of $E$, are all zero by assumption). In particular, $p_{n+k}(TE)$ is nonzero.
But $TE$ is the Whitney sum of $T^{\textup{vert}}E$ (the vertical tangent microbundle) and a stably trivial
(vector) bundle, since the base $M$ of $p\co E\to M$ is stably framed. Therefore $p_{n+k}(TE)=p_{n+k}(T^{\textup{vert}}E)$
is nonzero. Since $T^{\textup{vert}}E$ has fiber dimension $2n$, it follows
that $p_{n+k}$ is nonzero in $H^*(B\STOP(2n);\QQ)$. \qed

\begin{prop} \label{prop-concess} In theorem~\ref{thm-gmain} the following
condition on $E\to M$ can be added to the list: there exists a topological embedding $u\co \RR^{2n}\times M\to E$, over $M$, such that the fiber bundle $E\smin u(0\times M)\to M$ obtained by restriction of $E\to M$ admits a fiberwise
smooth structure.
\end{prop}
This will be fully established and used only in section~\ref{sec-decomp} to prove that the classes
$p_{n+k}\in H^{4n+4k}(B\TOP(2n))$ of theorem~\ref{thm-pmain} detect elements in the rationalized homotopy groups of $B\TOP(2n)$.

\begin{rem} \label{rem-Morlet}
It is a theorem of smoothing theory, Morlet style
\cite{MorletCR},\cite[Essay V]{KirbySiebenmann},\cite[Thm.4.4]{BurghLash} that $\diff_\partial(D^m)\simeq \Omega^{m+1}(\TOP(m)/\Or(m))$
for $m\ne 4$, where $\TOP(m)/\Or(m)$ is short for the homotopy fiber of
$B\Or(m)\to B\TOP(m)$. This fact makes the spaces $B\TOP(m)$ important in differential
topology. Proposition~\ref{prop-diskdiff} below spells out some consequences of theorem~\ref{thm-gmain}
(and the related proposition~\ref{prop-Crowleyquest}) for the rational homotopy type of $\diff_\partial(D^{2n})$. See also remark~\ref{rem-otherdisks}.
\end{rem}

\begin{rem} \label{rem-otherdisks} Farrell-Hsiang \cite{FaHsi} and Watanabe \cite{Watanabe1}, \cite{Watanabe2}
show in very different ways that for many $m$ the inclusion $B\Or(m) \to B\TOP(m)$ fails to be a rational
homotopy equivalence. The Farrell-Hsiang argument uses the smoothing theory connection between $\TOP(m+1)/\TOP(m)$
and the space of smooth pseudo-isotopies of $D^m$ (similar in spirit to remark~\ref{rem-Morlet}) and then,
following Waldhausen, the deep relationship between smooth pseudo-isotopy spaces and algebraic K-theory.
The force of this approach is limited by the known stability range (currently \cite{Igu}, Igusa's stability theorem)
for smooth pseudo-isotopy spaces. That is to say, Farrell and Hsiang get a complete description of $H^*(B\TOP(m))$ for
$\star\le 4m/3-c$ for a constant $c$,
and nothing beyond that range. (The range was improved recently in \cite{RW17} and \cite{Kra20}.
Igusa \cite{Igu02} independently used higher Franz-Reidemeister torsion classes to detect non-trivial rational cohomology
in $B\diff_\partial(D^{2n})$, in agreement with the Farrell-Hsiang result.) ---
Watanabe writes about $\diff_\partial(D^m)$, related to $\TOP(m)$ by remark~\ref{rem-Morlet},
and uses configuration space integrals to
find nontrivial rational cohomology classes in rather high degrees. It is not absurd to assume a close connection between
configuration space integrals and
functor calculus. See \cite{Volic1}, \cite{Volic2} for example. Therefore the articles \cite{Watanabe1}, \cite{Watanabe2} by
Watanabe could be closely related to this article here, which relies on functor calculus (manifold calculus).
\end{rem}

\subsection{Overview} Section~\ref{sec-homhirz} is a discussion of
Pontryagin numbers and the related kappa classes (of certain bundles of manifolds)
from a computational point of view.
Section~\ref{sec-calculus} is a summary of results from manifold calculus, with emphasis on applications to
self-embeddings of a specific smooth manifold $W_\parpun$ which is also prominent in section~\ref{sec-homhirz}.
Section~\ref{sec-conhoto} makes a link between
sections~\ref{sec-calculus} and~\ref{sec-homhirz}. Although this is rather weak, it turns out to be useful later where we
need to understand geometric properties of manifold bundles whose construction relies to some extent on non-geometric methods.
Section~\ref{sec-diss} draws on all the previous sections and provides a few computations
which, together, make up the proof of theorem~\ref{thm-gmain}. (Section~\ref{sec-diss} also has a preview of its own.)
Section~\ref{sec-decomp} is a review of the proof which leads us, essentially without further computations,
to a stronger conclusion: the ``surreal'' Pontryagin classes of theorem~\ref{thm-pmain}
detect elements in the rational homotopy groups $\pi_{4n+4k}(B\TOP(2n))\otimes\QQ$. Appendix~\ref{sec-compodesc}
is an abstract essay on the theme \emph{semi-direct product formulas for monoids of homotopy
automorphisms}. Semi-direct product expressions are heavily used in this article from section~\ref{sec-homhirz} onwards.

\subsection{Acknowledgments and apologies}
This article grew out of notes intended to clarify an obscure story about rational Pontryagin classes, manifold calculus
and parametrized surgery, told by me to S{\o}ren Galatius and Oscar Randal-Williams at the ICM satellite conference in Dalian, China, in August 2014.
For many years, from around 1994 to 2012, I believed that $p_n=e^2$ in
$H^{4n}(B\STOP(2n))$ for all $n$, and consequently that $p_{n+k}=0$
in $H^{4n+4k}(B\TOP(2n))$ for all $n,k>0$. An elaborate strategy for proving this was developed, based on remark~\ref{rem-Morlet}
and using ideas from singularity theory. I offer my apologies to those who sat through talks or series of talks on this unfortunate project,
and to Rui Reis who, more misguided than guided by me, participated in the project.

I am grateful to both S{\o}ren Galatius and Oscar Randal-Williams for their encouragement and interest, and for their
work on parametrized surgery and homological stability \cite{GalRW2014}, \cite{GalRW2018} which is used here.
I am also indebted to S{\o}ren Galatius for reading an older version of this article and making
many suggestions for improvement. Two of these are incorporated in the formulation
of theorem~\ref{thm-pmain}. In the older formulation, $n+k$ was subject to a curious divisibility condition, satisfied in
most cases but not all. S{\o}ren Galatius supplied a lemma which made this restriction superfluous. (His lemma is now
subsumed in a much more general statement on Hirzebruch $L$-polynomials \cite{BeBe} which was proved recently by
Berglund and Bergstr\"om.) In the older formulation of theorem~\ref{thm-pmain}, there was a bound of the form
$k\le n/2-$const.; Galatius had some ideas on how that could be improved and they were implemented. He also pointed
out that my proof of theorem~\ref{thm-pmain} proved the stronger statement theorem~\ref{thm-gmain}.

More recently, at a workshop in Bonn, September 2019, Alexander Kupers and Manuel Krannich gave some new arguments
suggesting that theorems~\ref{thm-pmain} and~\ref{thm-gmain} are still far from optimal. See remark~\ref{rem-KupKran}.

An earlier version of this article contained arguments based on the multiplicativity
of the sequence of the Hirzebruch $\sL$-polynomials which were incomplete or careless.
This was kindly pointed out to me by Martin Olbermann.

Finally, I am hugely indebted to the exacting and wise anonymous reviewer who had many suggestions
for making the paper more readable. Most of these I hope to have implemented.

\emph{Financial support.} This work was supported by the A.v.Humboldt foundation through a Humboldt Professorship award, 2012-2017; at
an earlier stage, 2008-2011, by the Engineering and Physical Sciences Research Council (UK),
Grant EP/E057128/1; and at the editing stage, 2019-2020, by the Deutsche Forschungsgemeinschaft
(DFG, German Research Foundation) under Germany's Excellence Strategy EXC 2044-390685587, Mathematics M\"unster: Dynamics-Geometry-Structure.

\section{Homotopical description of some kappa classes} \label{sec-homhirz}
\subsection{Two useful manifolds and their automorphisms} \label{subsec-deleting} For a positive
integer $g$, let $W_g$ be the connected sum of $g$ copies of $S^n\times S^n$ (an oriented closed manifold).
Let $W_{g,1}$ be the manifold obtained from $W_g$ by deleting the interior
of a smoothly embedded codimension zero (compact) disk $D$ in $W_g$:
\[   W_{g,1}= W_g\smin \intr(D). \]
Then $\partial W_{g,1}\cong S^{2n-1}$. The manifolds $W_g$ and $W_{g,1}$ are ubiquitous in the work of
Galatius and Randal-Williams on parametrized surgery, for example \cite{GalRW2014}, \cite{GalRW2018}.
They are also important here, for large but otherwise unspecified $g$, because the fiber of the
bundle $E\to M$ in theorem~\ref{thm-gmain} is $W_g$ or equivalently $W_{g,1}/\partial W_{g,1}$.
Indeed the construction of $E\to M$ makes some use of parametrized surgery.

Select a point
$z$ in $\partial W_{g,1}$ once and for all. Put
\[ W_\parpun =W_{g,1}\smin\{z\}. \]
There is a forgetful or restriction-induced map
$B\diff_\partial(W_{g,1}) \to B\diff_{\partial}(W_\parpun)$.
Here $\diff_{\partial}(-)$ generally refers to (topological or simplicial groups of) diffeomorphisms from a manifold to
itself which extend the identity on the boundary. (For the topological group interpretation, use the $C^\infty$
compact-open topology \cite[\S2.1]{Hirsch}. For the simplicial group interpretation, a $k$-simplex in $\diff_\partial(M)$
is a diffeomorphism $\Delta^k\times M \to \Delta^k\times M$ which is over $\Delta^k$ and which extends the
identity on $\partial M\times\Delta^k$.) Therefore $\diff_{\partial}(W_\parpun)$ is already defined, no
less than $\diff_\partial(W_{g,1})$. But it is often useful to think of $W_{g,1}$ as the one-point compactification
of $W_\parpun$~, and so to think of elements of $\diff_{\partial}(W_\parpun)$ as \emph{homeomorphisms} $W_{g,1}\to W_{g,1}$ which
extend the identity on $\partial W_{g,1}$ and restrict to diffeomorphisms $W_\parpun\to W_\parpun$\,.

For comparison, we also need spaces of smooth embeddings. Suppose that $M$ and $N$ are smooth manifolds with boundary,
and that a preferred smooth embedding $\partial M\to \partial N$ has been selected. Let $\emb_\partial(M,N)$ be the
space of neat smooth embeddings $f\co M\to N$ which extend the
preferred embedding $\partial M\to N$. (Neatness means that $\partial M$ is the transverse preimage of $\partial N$
under $f$.) Again we use the $C^\infty$ Whitney topology. Alternative: $\emb_\partial(M,N)$ is a fibrant simplicial set
where a $k$-simplex is a smooth embedding $\Delta^k\times M \to \Delta^k\times N$ which is over $\Delta^k$, extends the preferred embedding
$\Delta^k\times \partial M \to \Delta^k\times \partial N$ and is fiberwise neat. 

[Showing that the simplicial set $X:=\emb_\partial(M,N)$ is fibrant: this is left as an exercise with the following suggestions.
Without loss of generality, $N$ is a neat smooth submanifold of $Q=\RR^p\times[0,\infty)$ for some $p\gg 0$. Introduce a simplicial set $Y$ where a $k$-simplex is
a smooth \emph{map} $\Delta^k\times M \to \Delta^k\times N$ over $\Delta^k$ which takes
$\Delta^k\times \partial M$ to $\Delta^k\times \partial N$, extends the preferred embedding
$\Delta^k\times \partial M \to \Delta^k\times \partial N$ and is transverse to $\Delta^k\times\partial N$.
Introduce another simplicial set $Z$ where a $k$-simplex is
a smooth map $\Delta^k\times M \to \Delta^k\times \RR^p\times\RR$ over $\Delta^k$ which extends
the composition $\Delta^k\times \partial M\to \Delta^k\times\partial N\hookrightarrow \Delta^k\times\RR^p\times\{0\}$.
There are inclusions of simplicial sets $X\to Y\to Z$. Show that if $Y$ is
fibrant, then $X$ is fibrant. Using a smooth neighborhood retraction for $N$ in $Q$, show that
if $Z$ is fibrant, then $Y$ is fibrant. But $Z$ is fibrant because it is isomorphic to a simplicial group.]

\begin{lem} \label{lem-obsercalc} The inclusion
$\diff_\partial(W_\parpun)\to \emb_\partial(W_\parpun,W_\parpun)$ is a homotopy equivalence, provided $n\ge 3$.
\end{lem}

\proof Let $U$ be a standard open neighborhood (a half-disk) of $z$ in $W_{g,1}$, so that $V:=W_{g,1}\smin U$ is a compact smooth manifold
with corners. The boundary $\partial V$ is a union $\partial_0 V\cup \partial_1 V$ where $\partial_0V=V\cap \partial W_{g,1}$
and $\partial_1V$ is the closure of the complement of $\partial_0V$ in $\partial V$. It is easy to see that in the diagram
\[  \diff_\partial(W_\parpun) \hookrightarrow \emb_\partial(W_\parpun,W_\parpun) \lra \emb_{\partial_0}(V,W_\parpun), \]
both the right-hand arrow and the composite arrow are homotopy equivalences. (Indeed, a smooth embedding
$f\co V\to W_\parpun$ which is the identity on $\partial_0V$ must preserve the intersection form on the $n$-dimensional
integral homology. Since the intersection form is nondegenerate, $f$ must be a homotopy equivalence.
Since $n\ge 3$, this implies that the closure of the complement of $\im(f)$ is a collar, diffeomorphic to
$\partial_1V\times[0,1)\cong f(\partial_1V)\times[0,1)$. The space of diffeomorphisms from
$\partial_1V\times[0,1)$ to itself
extending the identity on $\partial_1V\times\{0\}$ is contractible.) \qed

\medskip
Manifold calculus (not to be confused with calculus on
manifolds) is a type of functor calculus concerned with contravariant functors
from certain categories of manifolds (of a fixed dimension) to spaces.
Among the most successful applications of manifold calculus are homotopical descriptions of
spaces of smooth embeddings $\emb(M,N)$.
In such a case the codimension should be at least three for the machine to function. But it is important to have the correct
interpretation of codimension. The geometric dimension of the target $N$ matters, and the
handle dimension of the source $M$; see definition~\ref{defn-hadim} just below.
In our situation, $M=W_\parpun=N$, the geometric dimension of the target is certainly $2n$. Since we are
after a space of embeddings relative to a fixed diffeomorphism of the boundaries, $\emb_\partial(W_\parpun,W_\parpun)$,
the handle dimension of the source must be computed relative to the boundary and it is $n$. Therefore the codimension count is $2n-n=n$.

\smallskip
In this way it turns out that both $B\diff_\partial(W_{g,1})$ and
$B\diff_{\partial}(W_\parpun)$ are homologically rather accessible, but for very different reasons: the access is through
parametrized surgery in the first case, through manifold calculus in the second.

\begin{defn} \label{defn-hadim} The absolute handle dimension
of a compact smooth $m$-dimensional $M$ with boundary $\partial M$ is the maximum index of handles in a finite handle decomposition of $M$,
as expected. But we have to define the handle dimension of a smooth $m$-dimensional $M$ \emph{relative} to $\partial M$, and $M$ as well as
$\partial M$ may be noncompact.
In such a case we will say that $M$ has handle dimension $\le k$ \emph{relative to $\partial M$}
if $M$ can be obtained from the trivial cobordism $(\partial M\times[0,1];\partial M\times\{0\},\partial M\times\{1\})$
in finitely many steps of the following kind.
\begin{itemize}
\item[(a)] Simultaneously and properly attaching a family of $m$-dimensional handles, all of index $\le k$, to the outgoing boundary $L$
of a cobordism with incoming boundary $\partial M$, and smoothing corners, thereby creating a new $m$-dimensional cobordism
with incoming boundary $\partial M$ and a new outgoing boundary. (The word \emph{properly} means that every compact region of $L$ meets
only finitely many of the new handles.)
\item[(b)] Last step (compulsory): deleting the outgoing boundary of a cobordism with incoming boundary $\partial M$.
\end{itemize}
For example, the $2n$-dimensional $W_\parpun$ with boundary $\cong \RR^{2n-1}$ can be constructed from
the trivial cobordism $(\RR^{2n-1}\times[0,1];\RR^{2n-1}\times\{0\},\RR^{2n-1}\times\{1\})$ in one step of type (a),
where $2g$ handles of dimension $2n$ and index $n$ are attached simultaneously to the outgoing boundary $\RR^{2n-1}\times\{1\}$,
then one step of type (b). Another example: if $M$ is compact, then its handle dimension relative to $\partial M$ is
always $m=\dim(M)$. In this case the last step (b), though compulsory, is superfluous because the outgoing boundary at that stage is empty.
\end{defn}

\begin{rem} There is a homotopy fiber sequence
\begin{equation} \label{eqn-weissfib}
 F\lra B\diff_\partial(W_{g,1}) \lra B\diff_{\partial}(W_\parpun) \end{equation}
where $F$ is a delooping of $\diff_\partial(D^{2n})$. Indeed the restriction map
\[ \diff_\partial(W_{g,1}) \to \emb_{\partial_0}(V,W_{g,1}) \]
(notation as in the proof of lemma~\ref{lem-obsercalc}) is a fibration and its
fiber is homotopy equivalent to $\diff_\partial(D^{2n})$. Here $\emb_{\partial_0}(V,W_{g,1})$ can be replaced by the homotopy equivalent
$\diff_{\partial}(W_\parpun)$.

For $n\ne 3$, the set $\pi_0F$ can be identified with the
\emph{inertia group} of $W_g$, the subgroup of $\Theta_{2n}$ (group of oriented 2n-dimensional smooth
homotopy spheres modulo oriented diffeomorphism) consisting of the elements represented by homotopy $2n$-spheres
$M$ such that $W_g\# M$ is diffeomorphic to $W_g$. (Exercise for the reader.) The reviewer points out that
the inertia group of $W_g$ was determined by Wall \cite{Wal62} and Kosinski \cite{Kos67}. It is trivial, so that
$F$ is path connected.

At an earlier stage of this investigation, the plan was to use the
homotopy fiber sequence \eqref{eqn-weissfib}, parametrized surgery
(for an analysis of $B\diff_\partial(W_{g,1})$), manifold calculus
(for an analysis of $B\diff_{\partial}(W_\parpun)$) and Morlet's result as in remark~\ref{rem-Morlet}
to obtain information on $B\TOP(2n)$. This plan was abandoned later for a more direct strategy. Meanwhile Kupers \cite{Kupers}
found some very good uses for \eqref{eqn-weissfib} as it is.
\end{rem}

\medskip
It seems appropriate to end this (sub)section with a few words on topological monoids and $A_\infty$-spaces,
a.k.a. $\mathcal E_1$-algebras in the category of spaces.
A good example for our purposes is $\emb_\partial(W_\parpun,W_\parpun)$. In section~\ref{sec-calculus} we will develop
homotopical descriptions of this based on manifold calculus, and in doing so we will also employ
structural concepts like semi-direct products of topological monoids. It is not a problem
to set up a homotopical framework for topological monoids or for $A_\infty$-spaces
(e.g. a Quillen model category structure, as done in appendix~\ref{sec-compodesc} for \emph{simplicial} monoids), but
in the body of the paper no such framework is specified. We only make one step towards that in saying that
a (continuous) homomorphism of topological monoids or an $A_\infty$-map between $A_\infty$ spaces is a
\emph{weak equivalence} if the underlying map of based spaces is a weak equivalence. This is also meant as
an invitation to readers to make sense of names (for topological monoids or $A_\infty$-spaces) in their own way, within reason,
and as an excuse for an author reluctant to clutter the text with details of little consequence.

\subsection{Decomposable vs.~indecomposable}
Let $\sP$ be a polynomial with rational coefficients in the Pontryagin classes $p_1,p_2,p_3,\dots$, homogeneous of
degree $4n+4k$ in the cohomological sense. (For example, if $n+k=10$, then the polynomial
$\frac{5}{8}p_{10}+\frac{3}{7}p_2p_3p_5-21p_2^5$ qualifies.)
By
\[ \kappa_t(\sP)\in H^{2n+4k}(B\homeo_\partial(W_{g,1})) \]
we mean the class obtained as follows. Let
$(E,\partial E)\to B\homeo_\partial(W_{g,1})$ be the tautological fiber bundle pair with fiber pair
$(W_{g,1},\partial W_{g,1})$. Evaluate $\sP$ on the vertical tangent (micro-)bundle of $E$ (which is stably trivalized
over $\partial E$) to obtain a class in $H^{4n+4k}(E,\partial E)$. Apply integration along the
fibers of $(E,\partial E)\to B\homeo_\partial(W_{g,1})$
to get a class in $H^{2n+4k}(B\homeo_\partial(W_{g,1}))$. This is $\kappa_t(\sP)$.

A practical description is as follows. Think of $H_{2n+4k}(B\homeo_\partial(W_{g,1}))$ as
rationalized stable homotopy $\pi^s_{2n+4k}(B\homeo_\partial(W_{g,1}))\otimes\QQ$
and represent an element $x$ of the stable homotopy group
$\pi^s_{2n+4k}(B\homeo_\partial(W_{g,1}))$ by a map from $M$ to $B\homeo_\partial(W_{g,1})$,
where $M$ is a stably framed closed smooth manifold of dimension $2n+4k$. (It may be necessary to replace $x$ by $cx$
where $c$ is a nonzero integer.) This determines a fiber bundle
$E_M\to M$ with fiber $W_g\cong W_{g,1}/\partial W_{g,1}$ (by pulling back $E\to B\homeo_\partial(W_{g,1})$ and collapsing boundary
spheres to points). Therefore $E_M$ is a closed oriented manifold of dimension $4n+4k$. Then we have
\[  \langle \kappa_t(\sP),x\rangle = \langle \sP(TE_M),[E_M]\,\rangle,  \]
i.e., $\langle \kappa_t(\sP),x\rangle$ is the tangential Pontryagin number determined by $E_M$ and $\sP$. Furthermore,
we can also write $\sP(T^{\textup{vert}}E_M)$ instead of $\sP(TE_M)$, since the tangent bundle of $M$ is stably trivialized.
This point of view will be important later.

\medskip
In the following we also write $\kappa_t(\sP)\in H^{2n+4k}(B\diff_\partial(W_\parpun))$ for the
image of $\kappa_t(\sP)$, as just defined, under the homomorphism in cohomology induced by the inclusion of
$B\diff_\partial(W_\parpun)$ in $B\homeo_\partial(W_{g,1})$.
Here we obtain a description of the classes $\kappa_t(\sP)$ in some special cases (after some preliminaries).

\begin{defn} \label{defn-semidirect}
Let $X$ and $Y$ be topological monoids; let $\alpha$ be a \emph{right} action of $X$ on $Y$ by
monoid endomorphisms. The \emph{semi-direct product} $X\ltimes Y$ is the topological monoid with underlying
space $X\times Y$, unit element $(1_X,1_Y)$ and associative multiplication given by
$(x_0,y_0)\cdot(x_1,y_1):= (x_0x_1,\alpha(y_0,x_1)y_1)$. 
\end{defn}

Here is an instance of a semi-direct product that we will need.
Let $X=\haut_\partial(W_{g,1})$ consist of the homotopy automorphisms of $W_{g,1}$ relative to the boundary.
View $W_{g,1}$ as a based space by taking the antipode of $z$ in $\partial W_{g,1}=S^{2n-1}$ as the base point.
View $Y:=\map_*(W_{g,1},\SOr(2n))$ as a topological group with pointwise multiplication.
The group-like monoid $X$ acts on the right of $Y$ by group endomorphisms, i.e., maps from $W_{g,1}$ to
$\SOr(2n)$ can be composed with homotopy automorphisms of $W_{g,1}$. We want to relate the group
$\diff_\partial(W_\parpun)$ to
\[  X\ltimes Y= \haut_\partial(W_{g,1})\ltimes \map_*(W_{g,1},\SOr(2n)). \]
To that end we introduce some modifications of $X$ and $Y$. Let $X'\subset X$ consist of the
homotopy automorphisms of $W_{g,1}$ relative to the boundary which map interior to interior. Let
$\GL^{\textup{pos}}(2n)$ be the group of orientation-preserving linear automorphisms of $\RR^{2n}$
and let $A$ be the contractible subgroup of those which restrict to the identity on the
standard $\RR^{2n-1}\subset\RR^{2n}$. Let $Y'$ be topological group of maps from $W_\parpun$ to
$\GL^{\textup{pos}}(2n)$ which take the base point to $A$. We have an inclusion $X'\to X$
and a forgetful map $Y\to Y'$ (both monoid homomorphisms and weak equivalences). The action
of $X$ on $Y$ restricts to an action of $X'$ on $Y$ which descends or extends to an action of $X'$ on $Y'$.
As a result we obtain monoid
homomorphisms  $X\ltimes Y \leftarrow X'\ltimes Y \to X'\ltimes Y'$
which are also weak equivalences.

Now choose a trivialization of the tangent vector bundle $TW_\parpun$ once and for all, i.e.,
an isomorphism with a trivial vector bundle $W_\parpun\times \RR^{2n}\to W_\parpun$. (Suppose or arrange that
this takes the tangent space of $\partial W_\parpun$ at the base point to the standard $\RR^{2n-1}\subset \RR^{2n}$.)
This allows us to think of the differential $df$ of any $f\in \diff_\partial(W_\parpun)$ as a \emph{map}
from $W_\parpun$ to $\GL^{\textup{pos}}(2n)$,
and more precisely as an element of $Y'$. In this way, $\diff_\partial(W_\parpun)$ becomes
a submonoid of $X'\ltimes Y'$, by dint of $f\mapsto (f,df)$. There is an
induced map of classifying spaces
\[  v\co B\diff_\partial(W_\parpun) \lra B(X'\ltimes Y'). \]
For ease of notation and by abuse of notation, we often write this in the form
\begin{equation}  \label{eqn-oldsemidir1}
v\co B\diff_\partial(W_\parpun) \lra B(\haut_\partial(W_{g,1})\ltimes \map_*(W_{g,1},\SOr(2n))) \end{equation}
neglecting the small distinction between $X'\ltimes Y'$ and $X\ltimes Y$.
In addition we have the forgetful projection
\begin{equation} \label{eqn-oldsemidir2}
 u\co B(\haut_\partial(W_{g,1})\ltimes \map_*(W_{g,1},\SOr(2n)))\lra B(\haut_\partial(W_{g,1})). \end{equation}

\begin{rem} \label{rem-easyuv} The maps $u$ and $v$ can be understood in a more geometric fashion by viewing
$B(\haut_\partial(W_{g,1})\ltimes \map_*(W_{g,1},\SOr(2n)))$ as a classifying space
for data of the following type: on a space $X$,
\begin{itemize}
\item[(i)] a fibration pair $(E,\partial E)\to X$ with fiber pair
$\simeq (W_{g,1},\partial W_{g,1})$,
\item[(ii)] a trivialization of the boundary fibration $\partial E\to X$ so that $\partial E\cong X\times \partial W_{g,1}$
\item[(iii)] and a map $\gamma$ from $E$ to $B\SOr(2n)$ taking the base point of each fiber $E_b$
to the base point. (The base point of $E_b$ is the base point of $\partial E_b\cong \partial W_{g,1}$.)
\end{itemize}
The map $\gamma$ in (iii) is subject to a mild condition: the restriction of $\gamma$ to any
fiber of $E\to X$ is based nullhomotopic. --- For even more precision, write $H$ for the semidirect product
$\haut_\partial(W_{g,1})\ltimes \map_*(W_{g,1},\SOr(2n))$.
We may write $h=(h_1,h_2)$ for $h\in H$,
where $h_1\in \haut_\partial(W_{g,1})$ and $h_2\in\map_*(W_{g,1},\SOr(2n))$.
In any case $H$ is a group-like topological monoid acting on $W_{g,1}$.
The action defines a \emph{transport category} $\texttt{T}$
(a.k.a. Grothendieck construction) with object space $W_{g,1}$ and morphism space $H\times W_{g,1}$,
so that the source of a morphism $(h,x)$ is $x$ and the target is $h_1(x)$, and composition
is defined by $(g,h_1(x))\circ (h,x)= (gh,x)$.
Let $E_u$ be the classifying space of the transport category; this is the universal example of $E$ in (i).
Viewing $\SOr(2n)$ as a category with one object, let
$\gamma_u\co E_u\to B\SOr(2n)$ be the map induced by the continuous functor $\varphi\co\texttt{T}\to \SOr(2n)$
taking a morphism $(h,x)$ in $\texttt{T}$ to $h_2(x)\in \SOr(2n)$.
Functor property:
\[  \begin{array}{c}
\varphi((g,h_1(x))\circ (h,x))=\varphi(gh,x)=(gh)_2(x)= ((g_2\circ h_1)\cdot h_2)(x)\\
=g_2(h_1(x))\cdot h_2(x) = \varphi(g,h_1(x))\cdot \varphi(h,x).
\end{array}
\]
Of course $\gamma_u$ is the universal example of $\gamma$ in (iii).
\end{rem}

\smallskip
Let $\sL_{n+k}$ be the Hirzebruch polynomial of cohomological degree $4n+4k$ in the Pontryagin classes, so that
the tangential Pontryagin number associated with $\sL_{n+k}$ and a closed oriented manifold of dimension $4n+4k$ is the
signature of that manifold. (This works for topological manifolds just as it does for smooth manifolds,
as pointed out in the introduction.)

Let $\sJ_{\TOP}\subset H^*(B\TOP)\cong\QQ[p_1,p_2,\dots]$ be the ideal of
decomposable elements and let $\sJ_{\SOr(2n)}\subset H^*(B\SOr(2n))\cong \QQ[e,p_1,p_2,\dots,p_{n-1}]$ be the
ideal of decompo\-sable elements.

\begin{prop} \label{prop-uv} \emph{(i)} In the notation of displays~\eqref{eqn-oldsemidir1} and ~\eqref{eqn-oldsemidir2},
the class
\[ \kappa_t(\sL_{n+k})\in H^{2n+4k}(B\diff_\partial(W_\parpun)) \]
is in the image of $(uv)^*$,
i.e., it comes from $H^{2n+4k}(B\haut_\partial(W_{g,1}))$.

\emph{(ii)} \emph{(See also \cite[\S4.3]{RandWillNote}.)} There is a commutative diagram of homomorphisms of the following shape:
\[
\xymatrix@C=35pt{
\sJ^{4n+4k}_{\TOP} \ar[d]_-{\textup{restriction}} \ar[r]^-{\kappa_t}  &
H^{2n+4k}(B\diff_\partial(W_\parpun)) \\
\sJ^{4n+4k}_{\SOr(2n)} \ar@{..>}[r]^-{\kappa_\sJ} & \ar[u]_-{v^*}
H^{2n+4k}(B(\haut_\partial(W_{g,1})\ltimes \map_*(W_{g,1},\SOr(2n))))
}
\]
\end{prop}

\proof Statement (i) is an easy consequence of the Hirzebruch signature theorem.
For $x\in H_{2n+4k}(B\diff_\partial(W_\parpun))$ represented by a map
$M \to B\diff_\partial(W_\parpun)$, where $M$ is a closed smooth stably framed manifold of dimension $2n+4k$, we have
\[  \langle \kappa_t(\sL_{n+k}), x\rangle= \textup{signature of }E_M  \]
as explained earlier. Here $E_M\to M$ is the fiber bundle with fiber $W_g\cong W_{g,1}/\partial W_{g,1}$
determined by $M \to B\diff_\partial(W_\parpun)\subset B\homeo_\partial(W_{g,1})$.
This description of the class $\kappa_t(\sL_{n+k})$ makes it clear that $\kappa_t(\sL_{n+k})$
comes from $B\haut_\partial(W_{g,1})$. Indeed a bordism class of maps $M\to B\haut_\partial(W_{g,1})$ with $M$ as above still determines
a fibration $E_M\to M$ where the fibers are oriented Poincar\'e duality spaces of formal dimension $2n$ (homotopy
equivalent to $W_{g,1}/\partial W_{g,1}$ alias $W_g$),
and $E_M$ is therefore an oriented Poincar\'e duality space of formal dimension $4n+4k$. Oriented
Poincar\'e duality spaces, too, have signatures. (See also remark~\ref{rem-sigconf}.)

For the proof of statement (ii), let $B=B(\haut_\partial(W_{g,1})\ltimes \map_*(W_{g,1},\SOr(2n)))$.
There is a canonical fibration pair (beware new notation) $\pi\co (E,\partial E) \to B$
with fiber pair $(W_{g,1},\partial W_{g,1})$ and trivialized boundary fibration $\partial E\to B$.
This is pulled back from $B\haut_\partial(W_{g,1})$. There is a preferred choice of map
$\tau\co E\to B\SOr(2n)$. (See remark~\ref{rem-easyuv}.)
Let $\sP$ be a monomial in the classes $e,p_1,\dots p_{n-1}\in H^*(B\SOr(2n))$, of
cohomological degree $4n+4k$ and of ordinary degree $r>0$ (i.e., the sum of the exponents of $e$, $p_1$,
$p_2$, \dots, $p_{n-1}$ in $\sP$ is $r$). More precisely write
\[  \sP= e^{r_0}p_1^{r_1}\cdots p_{n-1}^{r_{n-1}}  \]
where $r_0+\cdots+r_{n-1}=r$.
Then $\tau^*\sP\in H^*(E)$.
We try to make sense of the fiber integral
\[
\int_\pi \tau^*\sP
\]
as an element of $H^*(B)$, although $\tau$ is not claimed or required to be nullhomotopic on $\partial E$.

Now we follow \cite{RandWillNote}. Randal-Williams makes the following observations (related to similar but much more complicated
statements in earlier drafts of this article here):
\begin{itemize}
\item[(a)] the elements $\tau^*e,\tau^*p_1,\dots,\tau^*p_{n-1}\in H^*(E)$ can all be lifted
across the inclusion-induced map
$H^*(E,\partial E)\to H^*(E)$;
\item[(b)] whenever we choose such lifts denoted $\bar\tau^*e,\bar\tau^*p_1,\dots, \bar\tau^*p_{n-1}$, the product
\[  (\bar\tau^*e)^{r_0}(\bar\tau^* p_1)^{r_1}\cdots (\bar\tau^*p_{n-1})^{r_{n-1}} \in H^*(E,\partial E) \]
is the same, independent of the choice of lifts, provided $r>1$.
\end{itemize}
Proof of (a): this uses specific properties of $\partial E\to B$. For us it is a trivalized
fibration and the fibers are spheres $S^{2n-1}$. So there is a preferred isomorphism
\[  H^*(\partial E) \xrightarrow{(s,t)} H^*(B)\oplus H^{*-(2n-1)}(B) \]
where $s$ is induced by the preferred section $B\to \partial E$
and $t$ is integration along the fibers.
By Stokes' theorem, the image of the restriction homomorphism from $H^*(E)$ to $H^*(\partial E)$
is contained in $\ker(t)$. Our assumptions
on $\tau$ imply that the same restriction homomorphism takes the elements
$\tau^*e,\tau^*p_1,\dots,\tau^*p_{n-1}$ to $\ker(s)$. Therefore they all map to
zero in $H^*(\partial E)$.  ---
The proof of (b) is easier. The difference of two homogeneous elements
in $H^*(E,\partial E)$ with the same image in $H^*(E)$ is in the image of the boundary homomorphism
$H^{*-1}(\partial E) \to H^*(E,\partial E)$.
Products in $H^*(E,\partial E)$ involving at least one factor coming from $H^{*-1}(\partial E)$,
and at least one more factor, are always zero. Therefore: change one lift at a time to verify that
it does not influence the product in (b). ---
Now we can define
\[
\int_\pi \tau^*\sP :=  \int_\pi (\bar\tau^*e)^{r_0}(\bar\tau^* p_1)^{r_1}\cdots (\bar\tau^*p_{n-1})^{r_{n-1}}
\]
provided $r=r_0+\cdots+r_{n-1}$ is $>1$. This amounts to the definition
of the lower horizontal arrow in the diagram in (ii). Commutativity of the diagram in (ii) can be
established as follows. We set up a larger diagram
\[
\xymatrix@C=35pt@R=15pt{
\sJ^{4n+4k}_{\TOP} \ar[d]_-{=} \ar[r]^-{\kappa_t}  &
H^{2n+4k}(B\diff_\partial(W_\parpun)) \\
\sJ^{4n+4k}_{\TOP} \ar[d]_-{\textup{restriction}} \ar[r]^-{\kappa_\sJ} &
H^{2n+4k}(B(\haut_\partial(W_{g,1})\ltimes \map_*(W_{g,1},\TOP))) \ar[d]^-{j^*} \ar[u]_-{v_1^*}   \\
\sJ^{4n+4k}_{\SOr(2n)} \ar[r]^-{\kappa_\sJ} &
H^{2n+4k}(B(\haut_\partial(W_{g,1})\ltimes \map_*(W_{g,1},\SOr(2n))))
}
\]
where $j^*$ is inclusion-induced, $v_1^*$ is analogous to $v^*$ in the diagram of (ii)
and the middle horizontal arrow is constructed like the lower horizontal arrow (using a specific selection of polynomial generators for
$H^*(B\TOP)$, namely, the Pontryagin classes). This larger diagram is evidently commutative.
But now $v_1^*=v^* j^*$, and the commutativity of the diagram in (ii) follows. \qed

\subsection{A naturality property}
Let $B\SOr(2n)_\QQ$ be the rationalization of the based space $B\SOr(2n)$. Let $q$ be a based map
from $B\SOr(2n)_\QQ$ to $B\SOr(2n)_\QQ$. It induces an $A_\infty$-map from
$\SOr(2n)_\QQ$ to $\SOr(2n)_\QQ$ and consequently a map
\[
\xymatrix@R=16pt{
B(\haut_\partial(W_{g,1})\ltimes \map_*(W_{g,1},\SOr(2n)_\QQ)) \ar[d]^-{q_\sharp} \\
B(\haut_\partial(W_{g,1})\ltimes \map_*(W_{g,1},\SOr(2n)_\QQ))
}
\]
(Compare remark~\ref{rem-easyuv}; for the present purposes $\gamma$ in \ref{rem-easyuv}~(iii) is a map
from $E$ to $B\SOr(2n)_\QQ$.) This leads us to a diagram
\begin{equation} \label{eqn-natu}
\begin{aligned}
\xymatrix@C=35pt{
\sJ^{4n+4k}_{\SOr(2n)} \ar[r]^-{\kappa_\sJ} \ar[d]^-{q^*} &
H^{2n+4k}(B(\haut_\partial(W_{g,1})\ltimes \map_*(W_{g,1},\SOr(2n)_\QQ))) \ar[d]_-{H^*(q_\sharp)}  \\
\sJ^{4n+4k}_{\SOr(2n)} \ar[r]^-{\kappa_\sJ} &
H^{2n+4k}(B(\haut_\partial(W_{g,1})\ltimes \map_*(W_{g,1},\SOr(2n)_\QQ)))
}
\end{aligned}
\end{equation}
where $\kappa_\sJ$ is the homomorphism of proposition~\ref{prop-uv} part (ii).

\begin{prop} \label{prop-muchridiculed} Diagram~\eqref{eqn-natu} commutes. \end{prop}

\proof For this purpose we can view $B\SOr(2n)_\QQ$ as a product of rational Eilenberg-MacLane spaces.
The map $f$ is determined, up to homotopy, by the graded ring homomorphism
$f^*\co H^*(B\SOr(2n)) \to H^*(B\SOr(2n))$.
That in turn is determined by its values on the preferred polynomial generators,
\[
\begin{array}{rcl}
e & \mapsto & \sQ_e(e,p_1,\dots,p_{n-1}) \\
p_1 & \mapsto & \sQ_1(e,p_1,\dots,p_{n-1}) \\
p_2 & \mapsto & \sQ_2(e,p_1,\dots,p_{n-1})  \\
& \vdots & \\
p_{n-1} & \mapsto & \sQ_n(e,p_1,\dots,p_{n-1}).
\end{array}
\]
where $\sQ_e,Q_1\dots,\sQ_{n-1}$ are rational polynomials in $n$ variables. ---
In the notation of the proof of proposition~\ref{prop-uv} part (ii), let us choose lifts
$\bar\tau^* e,\bar\tau^* p_1,\dots,\bar\tau^* p_{n-1}$ in $H^*(E,\partial E)$ of $\tau^*e,\tau^*p_1,\dots,\tau^*p_{n-1}\in H^*(E)$
respectively, once and for all.
Then the lower itinerary in diagram~\eqref{eqn-natu} is the ($\QQ$-linear) map given by
\[  \sJ^{4n+4k}_{\SOr(2n)}\ni  \sP \quad\mapsto
\quad \int_\pi \big(\sP\circ(\sQ_e,\sQ_1,\dots,\sQ_{n-1})\big)(\bar\tau^*e,\bar\tau^*p_1,\dots,\bar\tau^*p_n). \]
To describe the upper itinerary in diagram~\eqref{eqn-natu}, we use the abbreviation $\sigma:=f\tau$. It is the ($\QQ$-linear) map given by
\[  \sJ^{4n+4k}_{\SOr(2n)}\ni  \sP \quad\mapsto\quad \int_\pi \sP(\bar\sigma^*e,\bar\sigma^* p_1,\dots,\bar\sigma^* p_n)\,. \]
However, $\bar\sigma^*e$ and $\sQ_e(\bar\tau^*e,\bar\tau^*p_1,\dots,\bar\tau^*p_n)$ are elements of $H^*(E,\partial E)$ mapping to
the same element $(f\tau)^*e$ of $H^*(E)$; similarly $\bar\sigma^* p_j$ and $\sQ_j(\bar\tau^* e,\bar\tau^* p_1,\dots,\bar\tau^* p_n)$ are elements of $H^*(E,\partial E)$ mapping to
the same element $(f\tau)^*p_j$ of $H^*(E)$, for every $j=1,2,\dots,n-1$.
Therefore the argument which established observation (b) in the proof of proposition~\ref{prop-uv} part (ii) also shows that
\[  \big(\sP\circ(\sQ_e,\sQ_1,\dots,\sQ_{n-1})\big)(\bar\tau^*e,\bar\tau^* p_1,\dots,\bar\tau^* p_n)
=  \sP(\bar\sigma^*e,\bar\sigma^* p_1,\dots,\bar\sigma^* p_n). \]
This completes the proof. \qed

\section{Manifold calculus and configuration categories} \label{sec-calculus}
\subsection{Smooth embeddings vs. formal immersions} \label{subsec-embvsimm}
We begin with some vocabulary. Let $M$ and $N$ be smooth manifolds with boundary. Suppose that a preferred smooth
embedding $\partial M\to \partial N$ has been selected.

\begin{defn} \label{defn-fimm} A \emph{formal immersion} from $M$ to $N$, rel $\partial$, consists of a (continuous) map
$f\co M\to N$ and a map $\varphi\co TM\to TN$ covering $f$, together subject to the following conditions:
\begin{itemize}
\item[] $\varphi$ is linear and injective on the fibers;
\item[] $f$ extends the prescribed embedding $\partial  M \to \partial N$, while $\varphi$ agrees with the derivative
of $f$ on $T(\partial M)\subset TM$, and maps the outward half-spaces of $TM|_{\partial M}$ to the
outward half-spaces of $TN|_{\partial N}$.
\end{itemize}
Notation for the space of these formal immersions is $\fimm_\partial(M,N)$, or $\fimm(M,N)$
if $\partial M$ and $\partial N$ are empty.
\end{defn}
A genuine \emph{neat} smooth immersion $f\co M\to N$ rel $\partial$ determines a formal immersion $(f,\varphi)$
where $\varphi$ is the derivative $df$. (The neatness
condition on $f$ means that $\partial M$ is the \emph{transverse} preimage of $\partial N$ under $f$.)
Consequently there is an inclusion map $\imm_\partial(M,N)\to \fimm_\partial(M,N)$,
on the understanding that $\imm_\partial(M,N)$ is for the space of neat smooth immersions rel $\partial$.

\medskip
With lemma~\ref{lem-obsercalc} in mind, we will look for smooth embeddings or families of smooth embeddings
from $W_\parpun$ to itself, extending the identity on the boundary, with unusual derivatives. (This will in the end allow us
to establish theorem~\ref{thm-gmain}.) Equivalently, with lemma~\ref{lem-obsercalc} and definition~\ref{defn-fimm} in mind,
we will look for unusual families of formal immersions from $W_\parpun$ to itself,
rel $\partial$, which can be deformed into families of smooth embeddings, rel $\partial$. \newline
In \cite{BoavidaWeiss}, the standard theorems of manifold calculus applied to spaces of smooth
embeddings were reformulated in such a way that the obstructions to deforming families of formal immersions
to families of embeddings are particularly visible.
Here is a sample which is close to being the most useful variant for us. It is for manifolds without boundary
$M$ and $N$; later a version for manifolds with boundary will be stated. It uses certain categories $\config(M;r)$
and $\config(N;r)$ of ordered configurations in $M$ and $N$, of cardinality bounded above by a positive integer $r$.
Taken by itself it is neither useful nor difficult. To make it useful, combine it with the difficult theorem,
\cite{GoWe} relying on \cite{GoKl1}, \cite{GoKl2} and \cite{GoDisj}, that the standard
comparison map $\emb(M,N)\to T_r\emb(M,N)$
is highly-connected. Here $T_r\emb(M,N)$ denotes the $r$-th Taylor approximation of $\emb(M,N)$, both viewed
as contravariant functors in $M$. (In fact this comparison map is $((r+1)(c-2)+3-\dim(N))$-connected where $c$ is the
codimension, wisely determined as the difference between geometric dimension of $N$ and
handle dimension of $M$.) Recall or accept that $T_1\emb(M,N)$ is another name for $\fimm(M,N)$.

\begin{thm} \label{thm-main} \cite[Thm.5.1]{BoavidaWeiss} For any integer $r\ge 1$, there is a square
\[
\xymatrix{
T_r\emb(M,N) \ar[r] \ar[d] &  \ar[d]^-{\textup{specialization}} {\rmap_\fin(\config(M;r),\config(N;r))} \\
T_1\emb(M,N) \ar[r] & \rmap_\fin(\config^\loc(M;r),\config^\loc(N;r))
}
\]
which commutes up to a specified homotopy and is homotopy cartesian as such.
\end{thm}
It is time to explain what these configuration categories are.
Let $k\in\NN$. Write $\uli k=\{1,2,\dots,k\}$. The space of maps from $\uli k$ to $M$
comes with an obvious stratification. There is one stratum for each equivalence relation $\eta$ on $\uli k$~. The
points of that stratum are precisely the maps $\uli k\to M$ which can be factorized as projection from $\uli k$ to
$\uli k\,/\eta$ followed by an injection of $\uli k\,/\eta$ into $M$. \newline
We construct a topological category $\config(M)$
(category object in the category of topological spaces) whose object space is
\[  \coprod_{k\ge 0} \emb(\uli k\,,M), \]
that is, the topological disjoint union of the ordered configuration spaces of $M$ for each cardinality $k\ge 0$.
By a morphism from $f\in \emb(\uli k\,,M)$ to $g\in\emb(\uli\ell\,,M)$ we mean a pair consisting of a map $v\co \uli k\to \uli\ell$
and a Moore homotopy $\gamma=(\gamma_t)_{t\in [0,a]}$ from $f$ to $gv$ such that $\gamma$ in reverse is an
\emph{exit path} in the stratified space of \emph{all} maps from $\uli k$ to $M$. Equivalently, if $\gamma_s(x)=\gamma_s(y)$
for some $s\in [0,a]$ and $x,y\in \uli k$\,, then $\gamma_t(x)=\gamma_t(y)$ for all $t\in[s,a]$. (Andrade \cite{Andrade}
uses the expression \emph{sticky homotopy}.)
The space of all morphisms is therefore a coproduct
\[ \coprod_{\twosub{k,\ell\ge 0}{v\co\uli k\to\uli\ell}} P(v) \]
where $P(v)$ consists of triples $(f,g,\gamma)$ as above: $f\in\emb(\uli k\,,M)$,
$g\in\emb(\uli\ell\,,M)$ and $\gamma$ is the reverse of an exit path in $\map(\uli k\,,M)$ from $gv$ to $f$. Composition of morphisms
is obvious. It is obvious that in $\config(M)$, the maps \emph{source} and \emph{target} from morphism space
to object space are fibrations. (This is very useful. In general terms, if for a topological category at least
one of these two maps is a fibration, then the nerve is a Segal space.
See definition~\ref{defn-Segalsp}.)

\begin{prop} \label{prop-Miller1}
For fixed objects $f\in\emb(\uli k\,,M)$ and $g\in\emb(\uli\ell\,,M)$, the space of morphisms from $f$ to $g$
is homotopy equivalent $($per inclusion$)$ to the space of pairs $(v,\gamma)$
where $v\co \uli k\to \uli\ell$ as before and $\gamma=(\gamma_t)_{t\in[0,1]}$ is a path in $\map(\uli k,M)$ from $f$ to $gv$ such that
$\gamma_t$ is \emph{injective} for all $t$ strictly less than $1$.
\end{prop}
\proof This is a special case of \cite{Millershort}. \qed

\begin{cor} \label{cor-Miller1} The space of all morphisms in $\config(M)$ with fixed target $g\in \emb(\uli\ell\,,M)$
is homotopy equivalent to $\coprod_{j\ge 0}\emb(\uli j\,,U)$,
where $U$ is a standard $($open tubular$)$ neighborhood of the finite set $\im(g)$ in $M$.
\end{cor}
In short, we have a good understanding of morphism
spaces and object spaces in $\config(M)$, and they all somehow boil down to ordered configuration spaces of $M$ or
something very closely related.

\smallskip
Next some variants: $\config(M;r)$ is the full subcategory of $\config(M)$ with object space
$\coprod_{k=0}^r  \emb(\uli k\,,M)$.

The localized form $\config^\loc(M)$ is a type of comma category. Its objects are the morphisms in $\config(M)$
whose target has cardinality 1, in other words has the form $x\co \uli 1\to M$. The morphisms are commutative triangles
\[
\xymatrix@C=-10pt@R=30pt{ (f\co \uli k\to M) \ar[dr]_{\tau\circ\gamma} \ar[rr]^\gamma && (g\co \uli\ell\to M) \ar[dl]^\tau  \\
& (x\co \uli 1\to M)
}
\]
in $\config(M)$. (In making the step from $\config(M)$ to $\config^\loc(M)$ we lose one good property:
the map \emph{source} from morphism space to object space is no longer a fibration. But we still have a good degreewise
understanding of the nerve as a simplicial set.) By proposition~\ref{prop-Miller1},
the object space of $\config^\loc(M)$ is homotopy equivalent to the total space of a fiber bundle on $M$ whose fiber
at $x\in M$ is
\[  \coprod_{k\ge 0} \emb(\uli k\,,T_xM)~. \]
And $\config^\loc(M;r)$ is like $\config^\loc(M)$ but with cardinalities of configurations bounded above by $r$.
We write $\fin$ for the category whose objects are the finite sets $\uli\ell$, where $\ell\ge 0$, and whose morphisms
are all maps between these (not required to be order preserving). There are obvious forgetful functors
\[  \config^\loc(M) \to \config(M) \to \fin\,. \]
Next, the meaning of $\rmap_\fin$ in the right-hand column of the square in theorem~\ref{thm-main} ought to be explained.
It wants to say \emph{space of simplicial maps over the nerve of $\fin$ in the right derived sense}. More details
are given in section~\ref{subsec-Rezk}. \newline
Last not least, a few words on the horizontal maps in the square of the theorem are in order. Looking at
the top row for example, the point is that we have a natural transformation
\[ \emb(M,N) \lra \rmap_\fin(\config(M;r),\config(N;r)) \]
of contravariant functors in the variable $M$. It is easy to verify that the target, as a functor of $M$, satisfies
the conditions for a polynomial functor of degree $r$, in the sense of manifold calculus. Therefore, by
a universal property (in the derived sense) of $T_r$~, that natural map factors canonically through $T_r\emb(M,N)$.
The reasoning for the lower row is similar.
Let us note that $T_1\emb(M,N)$ is another way
to write $\imm(M,N)$, space of smooth immersions from $M$ to $N$,
if we assume that $\dim(M)\ne \dim(N)$ or, in case $\dim(M)=\dim(N)$, that $M$ has no closed component.
For the construction of the diagram it is wiser to stick to the $T_1\emb(M,N)$ notation because this emphasizes the
(aforementioned) universal property. How do we think of the lower horizontal arrow
in the diagram? The idea is as follows: inspired by the Smale-Hirsch $h$-principle,
we think of a smooth immersion $M\to N$ mainly as a \emph{formal} smooth immersion, i.e., a continuous map $f\co M\to N$ together with
linear injections $\psi_x\co T_xM\to T_{f(x)}N$, one for each $x\in M$ and depending continuously on $x\in M$.
The lower right-hand term in the square of the theorem has a very similar description: an element of it can be seen
as a continuous map $f\co M\to N$ together with functors (in the right derived sense)
\[ \psi_x\co \config(T_xM)\to \config(T_{f(x)}N) \]
depending continuously on $x\in M$. (See \cite[4.1]{BoavidaWeiss} and proposition~\ref{prop-unravelBoaWe} below, and further
explanations early on in section~\ref{subsec-locandlittle}.)
So the lower horizontal map in the square is obtained by noting that a linear injection
$T_xM\to T_{f(x)}N$ determines a functor between configuration categories,
$\config(T_xM)\to \config(T_{f(x)}N)$.

\medskip
The theorem has a variant for manifolds with boundary. For that we assume that $M$ and $N$ are manifolds with boundary
and a smooth embedding $\partial M\to \partial N$ has been fixed in advance. Write $M_-:=M\smin\partial M$.
We are interested in the space $\emb_\partial(M,N)$
of smooth neat embeddings which extend the selected embedding $\partial M\to \partial N$. (One reason for avoiding the
situation where we look at all neat embeddings $M\to N$, taking boundary to boundary, is that it may lead to a kind
of two-variable manifold calculus.) Here we need a new definition of the configuration category $\config(M)$.
\newline
Let $\finplus$ be the
category with objects $[k]=\{0,1,\dots,k\}$ for $k\ge 0$. These objects are viewed as based sets with base point $0$,
so the morphisms are the based maps. \newline
We make a topological category with object space
\[  \coprod_{k\ge 0} \emb_*([k],M/\partial M) \]
where $\emb_*(...)$ is for spaces of \emph{based} embeddings. (The base point of $M/\partial M$ is the class of
$\partial M$. Note that $M/\partial M$ is the pushout of $M\leftarrow \partial M\to *$.) By definition, a morphism from $f\co [k]\to M/\partial M$ to $g\co [\ell]\to M/\partial M$
consists of a based map $v\co [k] \to [\ell]$ and a Moore homotopy
$\gamma=(\gamma_t)_{t\in [0,a]}$ from $f$ to $gv$ such that each $\gamma_t$ is a based map,
and if $\gamma_s(x)=\gamma_s(y)$ for some $s\in [0,a]$ and $x,y\in \uli k$\,, then
$\gamma_t(x)=\gamma_t(y)$ for all $t\in[s,a]$.

The nerve of this topological category is $\config(M)$. It
comes with a forgetful map to $N\finplus$. (If $\partial M$ happens to be empty, this factors
through the inclusion $N\fin\to N\finplus$.) We define $\config^\loc(M)$ to be the same as $\config^\loc(M_-)$, so this comes
with a forgetful functor to $\fin$ as before.

\begin{thm} \label{thm-mainbdry} \cite[Thm.6.4]{BoavidaWeiss} For any integer $r\ge 1$, there is a square
\[
\xymatrix{
T_r\emb_\partial(M,N) \ar[r] \ar[d] &  \ar[d]^-{\textup{specialization}} \rmap^\partial_\finplus(\config(M;r),\config(N;r)) \\
T_1\emb_\partial(M,N) \ar[r] & \rmap^\partial_\fin(\config^\loc(M;r),\config^\loc(N;r))
}
\]
which commutes up to a specified homotopy and is homotopy cartesian as such. \qed
\end{thm}
Perhaps the right-hand column is self-explanatory. If not, let us agree that $M$ is equipped with a closed
collar $\bar U\cong \partial M\times[0,1]$, closure of an open collar $U$. We can focus attention on the neat
smooth embeddings $M\to N$ which are prescribed on $\bar U$. (There was a prescription on $\partial M$ to begin with,
in the shape of a smooth embedding $\partial M \to \partial N$. Choose an extension of that to $\bar U$ once and for all,
transverse to $\partial N$ and taking $\bar U\smin \partial M$ to $N\smin\partial N$.)
Then we have also prescribed
functors from $\config(U;r)$ to $\config(N;r)$ and from $\config^\loc(U;r)$ to $\config^\loc(N;r)$. In the right-hand column,
the notation $\rmap^\partial$ means that
we look for functors in the derived sense which extend these specified functors from $\config(U;r)$ to
$\config(N;r)$ and from $\config^\loc(U;r)$ to $\config^\loc(N;r)$. The
upper right hand term for example is the homotopy fiber of the restriction map
\[ \rmap_{\finplus}(\config(M;r),\config(N;r))\lra \rmap_{\finplus}(\config(U;r),\config(N;r)) \]
over the point determined by that selected functor $\config(U;r)\to \config(N;r)$.

\begin{rem} \label{rem-Miller2}
The description of $\config(M)$ just given is isomorphic to (but less wordy than)
what is called \emph{the particle model} in \cite[\S6]{BoavidaWeiss}. But there it is assumed that $\partial M$ is compact.
Now this compactness condition appears to be superfluous, and it is also unwanted here. It is true that $M/\partial M$
is metrizable if $\partial M$ is compact and non-metrizable if $\partial M$ is not compact.
This implies that $M/\partial M$ and powers $(M/\partial M)^k$ with various
sensible stratifications do not qualify as \emph{Quinn} stratified spaces if $\partial M$ is
noncompact, so that the work of Miller \cite{Miller} is not directly applicable to them. But in that situation
we can choose an increasing sequence of compact subsets $K_j$ of $M$ such that $K_j$ is contained in the
interior of $K_{j+1}$ and $\bigcup_j K_j=M$. Every morphism in $\config(M)$ can be realized in $K_j/(\partial M\cap K_j)$
for some $j\gg 0$. Reasoning along these lines one can show, using \cite{Miller} after all,
that the space of morphisms in $\config(M)$ with fixed target $g\in \emb_*([\ell],M/\partial M)$
is weakly homotopy equivalent to
\[ \coprod_{j\ge 0}\emb_*([j]\,,V) \]
where $V\subset M/\partial M$ is a standard neighborhood of $\im(g)$ in $M/\partial M$. (In other words,
the preimage of $V$ in $M$ is the disjoint union of an open collar on $\partial M$ and a tubular
neighborhood of $\{g(x)~|~x\in \uli j\,\}$\,.) This generalizes proposition~\ref{prop-Miller1}.
\end{rem}

\begin{rem} In theorem~\ref{thm-mainbdry},
the notation $\rmap^\partial_\finplus(\config(M;r),\config(N;r))$ with $\partial$ in the
superscript position is consistent with \cite{BoavidaWeiss}. It is inconsistent with the notation
$\emb_\partial(M,N)$, which is widely accepted. In defence of the superscript position it can be said
that the superscript $\partial$ refers to categories or Segal spaces \emph{over} $\config(M;r)$, respectively
$\config(N;r)$, whereas the subscript $\finplus$ refers to a category (or Segal space $N\finplus$) \emph{under}
$\config(M;r)$, respectively
$\config(N;r)$.
\end{rem}

Let us unravel what theorem~\ref{thm-mainbdry} means for $\emb_\partial(W_\parpun,W_\parpun)$ and for
the delooping $B\emb_\partial(W_\parpun,W_\parpun)$. We take $M=N=W_\parpun$ in theorem~\ref{thm-mainbdry}. Beware that $W_\parpun$ has
dimension $2n$, but its handle dimension relative to the boundary is $n$. It follows that
$\emb_\partial(W_\parpun,W_\parpun)$ maps to the homotopy pullback (holim) of
\begin{equation} \label{eqn-mainbdryspec1}
\begin{aligned}
\xymatrix@C=-90pt@R=40pt{
\fimm_\partial(W_\parpun,W_\parpun) \ar@/^1pc/[dr] &&
\ar@/_1pc/[dl]^-{\textup{ specialization}} \rmap^\partial_\finplus(\config(W_\parpun;r),\config(W_\parpun;r)) \\
 & \rmap^\partial_\fin(\config^\loc(W_\parpun;r),\config^\loc(W_\parpun;r))
}
\end{aligned}
\end{equation}
by a $((r+1)(c-2)+3-2n)$-connected map. Here $c=2n-n=n$, so that $(r+1)(c-2)+3-2n$ simplifies to $(r-1)n-2r+1$.
By lemma~\ref{lem-obsercalc},
the topological monoid $\emb_\partial(W_\parpun,W_\parpun)$ is group-like. The three spaces in \eqref{eqn-mainbdryspec1} are also
$A_\infty$-spaces (if not exactly topological monoids) and all maps in sight respect the $A_\infty$-structures. Therefore, assuming $(r-1)n-2r+1>0$,
we may discard the non-invertible components
in these three spaces without changing the meaning of the homotopy pullback.
We learn that the group-like monoid or group-like $A_\infty$-space $\emb_\partial(W_\parpun,W_\parpun)$ maps to the homotopy pullback of
\begin{equation} \label{eqn-mainbdryspec2}
\begin{aligned}
\xymatrix@C=-50pt@R=40pt{
\fimm_\partial^\times(W_\parpun,W_\parpun) \ar@/^1pc/[dr] &&
\ar@/_1pc/[dl]^-{\textup{specialization}} \haut^\partial_\finplus(\config(W_\parpun;r)) \\
 & \haut^\partial_\fin(\config^\loc(W_\parpun;r))
}
\end{aligned}
\end{equation}
by a $((r-1)n-2r+1)$-connected map which respects the $A_\infty$ structures.

We can make~\eqref{eqn-mainbdryspec2} more explicit by choosing a trivialization of
the tangent bundle $TW_\parpun$ as in the discussion leading up to~\eqref{eqn-oldsemidir1},
and by using \cite[Lem.4.1]{BoavidaWeiss} which is analogous to the $h$-principle of immersion theory.
In these circumstances we have compatible semidirect product expressions for two of the three terms in
~\eqref{eqn-mainbdryspec2}. (One of these expressions is fairly clear. The other one, although very
similar, is elaborately justified at an abstract level in appendix~\ref{sec-compodesc}. This takes care of
the compatibility, too. See especially corollary~\ref{cor-endoandframed} and~\ref{expl-endoconfigloc}.)
Then the homotopy pullback of \eqref{eqn-mainbdryspec2}, which we view as a homotopy pullback of group-like $A_\infty$-spaces,
turns into the homotopy pullback of
\begin{equation} \label{eqn-mainbdryspec2coord}
\begin{aligned}
\xymatrix@C=-70pt@R=40pt{
{\begin{array}{l} \haut_*(W_\parpun) \,\,\, \ltimes \\ \map_*(W_\parpun,\SOr(2n))
\end{array}} \ar@/^1pc/[dr] &&  \ar@/_1pc/[dl]^-{\textup{specialization}} \haut^\partial_\finplus(\config(W_\parpun;r)) \\
 & {\begin{array}{l} \haut_*(W_\parpun) \,\,\, \ltimes \\
\map_*(W_\parpun,\haut_\fin(\config(\RR^{2n};r)))
\end{array}}
}
\end{aligned}
\end{equation}
which we still view as a homotopy pullback of group-like $A_\infty$-spaces. This shows:
\begin{prop} \label{prop-unravelBoaWe}
The space $B\emb_\partial(W_\parpun,W_\parpun)$ maps to the \emph{base point path-component} of the
homotopy pullback of
\begin{equation} \label{eqn-mainbdryspec3}
\begin{aligned}
\xymatrix@C=-85pt@R=40pt{
{B\left(\begin{aligned} \begin{array}{l}  \haut_*(W_\parpun)\,\,\, \ltimes \\
\map_*(W_\parpun,\SOr(2n))
\end{array}\end{aligned}\right)} \ar@/^1pc/[dr] &&
\ar@/_1pc/[dl]^-{\textup{specialization}} B\haut^\partial_\finplus(\config(W_\parpun;r)) \\
 & {B\left(\begin{aligned}
\begin{array}{l}  \haut_*(W_\parpun) \,\,\, \ltimes \\
\map_*(W_\parpun,\haut_\fin(\config(\RR^{2n};r)))
\end{array} \end{aligned} \right)}
}
\end{aligned}
\end{equation}
by an $((r-1)n-2r+2)$-connected map. \qed
\end{prop}

\begin{rem} Configuration categories were not invented in \cite{BoavidaWeiss}. We took the
concept as we found it in Andrade's thesis \cite{Andrade}. Many alternative but weakly equivalent descriptions are given
in~\cite{BoavidaWeiss}. Configuration categories and closely related notions have also been very important in the recent development of
factorization homology. See for example \cite{AyFrRo}. Perhaps the underlying idea
that the ordered configuration spaces $\emb(\uli k\,,M)$ of a manifold $M$,
taken together for all $k$, should be organized into some bigger structure goes back to Fulton-MacPherson
\cite{FulMac} in a complex algebraic geometry setting. Axelrod and Singer \cite{AxSing} exported this idea to
differential topology. Sinha used the Axelrod-Singer formulation and emphasized its usefulness in
manifold calculus, concentrating on the case of a 1-dimensional source manifold. See for example \cite{Sinha}.
\end{rem}

\subsection{Functors in the right derived sense} \label{subsec-Rezk}
Generally we like to replace topological categories
$\mathcal A, \mathcal B, ...$ by their nerves $N\mathcal A, N\mathcal B, ...$
which are simplicial spaces. If these nerves are well behaved, we interpret \emph{functor from $\mathcal A$ to $\mathcal B$
in the derived sense} to mean \emph{simplicial map from $\mathcal A$ to $\mathcal B$ in the derived sense}. The space
of these is denoted
\[  \rmap(\mathcal A,\mathcal B). \]
It remains to be said what is meant by \emph{well-behaved} and what is meant by
the notation $\rmap(X,Y)$ for two simplicial spaces $X$ and $Y$.
\newline
Rezk has invented the concept of \emph{complete Segal space}. This is a simplicial space
which has roughly the properties that we expect from a nerve, but formulated in a homotopical way.
The following is copied from \cite{BoavidaWeiss}.

\begin{defn} \label{defn-Segalsp} A \emph{Segal space} is a simplicial space $X$ satisfying condition ($\sigma$) below.
If condition ($\kappa$) below is also satisfied, then $X$ is a \emph{complete Segal space}.
\end{defn}

\medskip
\begin{itemize}
\item[($\sigma$)] For each $n\ge 2$ the map $(u_1^*,u_2^*,\dots,u_n^*)$ from $X_n$ to the homotopy inverse limit of
the diagram
\[
\xymatrix{
X_1\ar[r]^-{d_0} & X_0 & \ar[l]_{d_1} X_1 \ar[r]^-{d_0} & \cdots & \cdots\ar[r]^-{d_0}& X_0 & \ar[l]_-{d_1}  X_1
}
\]
 is a weak homotopy equivalence. (The $u_i^*$ are iterated face operators corresponding to
the weakly order-preserving maps $u_i\co\{0,1\}\to\{0,1,2,\dots,n\}$ defined by $u_i(0)=i-1$ and $u_i(1)=i$.)
\end{itemize}
In order to formulate condition ($\kappa$) we introduce some vocabulary based on ($\sigma$). We call an element
$z$ of $\pi_0X_1$ \emph{homotopy left invertible} if there is an element $x$ of $\pi_0X_2$ such that $d_0x=z$
and $d_1x$ is in the image of $s_0\co \pi_0X_0\to \pi_0 X_1$. (In such a case $d_2x$ can loosely be thought of as
a left inverse for $z=d_0x$. Indeed $d_1x$ can loosely be thought of as the composition $d_2x\circ d_0x$, and by
assuming that this is in the image of $s_0$ we are saying that it is in the path component of an identity morphism.
We have written $d_0$, $d_1$, $s_0$ etc.~for maps induced on $\pi_0$ by the face and degeneracy operators.)
We call $z$ \emph{homotopy right invertible} if there is
an element $y$ of $\pi_0X_2$ such that $d_2y=z$
and $d_1x$ is in the image of $s_0\co \pi_0X_0\to \pi_0 X_1$.
Finally $z\in \pi_0X_1$ is \emph{homotopy invertible} if it is both homotopy left invertible and homotopy right invertible.
Let $X_1^w$ be the union of the homotopy invertible path components of $X_1$\,. It is a subspace of $X_1$\,.
\begin{itemize}
\item[($\kappa$)] The map $d_0$ restricts to a weak homotopy equivalence from $X_1^w$ to $X_0$\,.
\end{itemize}

\begin{defn} \label{defn-rezkfunctor}
A \emph{functor} from a complete Segal space $X$ to another complete Segal space $Y$ is just
a simplicial map $f\co X\to Y$~.
\begin{itemize}
\item[(i)] Such a functor is a \emph{weak equivalence} if and only if $f_n\co X_n\to Y_n$ is a weak homotopy
equivalence for all $n\ge 0$.
\item[(ii)] Suppose that $X$ and $Y$ are complete Segal spaces. By $\map(X,Y)$,
the space of all functors from $X$ to $Y$, we mean a simplicial set whose set of
$n$-simplices is the set of simplicial maps from $\Delta^n\times X$ to $Y$ (where $\Delta^n\times X$
has $k$-th term equal to $\Delta^n\times X_k$).
For a homotopy invariant (= right derived) notion
of \emph{space of all functors from $X$ to $Y$}\,, we require a cofibrant replacement of
$X^c$ of $X$ and a fibrant replacement $Y^f\to Y$ of $Y$.
Then we define $\rmap(X,Y)$ as $\map(X^c,Y^f)$, the space of simplicial maps
from $X^c$ to $Y^f$. Here we rely on a standard model category structure (to make sense of
cofibrant and fibrant replacements) on the category of
simplicial spaces in which a morphism $X\to Y$ is a weak equivalence, respectively fibration,
if $X_n\to Y_n$ is a weak equivalence, respectively fibration, for every $n\ge 0$.
The cofibrant and
fibrant replacements can be constructed in a functorial way.
\end{itemize}
\end{defn}
(\emph{End of quotation from \cite{BoavidaWeiss}}.) \newline
These ideas are not directly applicable to the nerves of $\config(M)$ and $\config(N)$,
which are ``Segal'' but not complete in the sense of Rezk. In fact the nerve of $\fin$
is ``Segal'' but not complete. But the reference functor $\config(M)\to \fin$
induces a map of the nerves (which are Segal spaces) which is
\emph{fiberwise complete} as in the following definition. (Again this is essentially a quotation from \cite{BoavidaWeiss}, but the
wording is not strictly the same.)

\medskip
Let $Y$ and $Z$ be simplicial spaces which satisfy ($\sigma$).  Let $f\co Y\to Z$
be a simplicial map. The following condition is an obvious variation on condition ($\kappa$) above.
\begin{itemize}
\item[($\kappa_{\rel}$)] The square
\[
\xymatrix{
Y_1^w \ar[d]^{d_0}  \ar[r] & Z_1^w \ar[d]^{d_0} \\
Y_0  \ar[r] & Y_0
}
\]
is a homotopy pullback.
\end{itemize}

\begin{defn} \label{defn-rezkfunctorover} We say that $f\co Y\to Z$ constitutes a
a fiberwise complete Segal space over $Z$ if $Y$ and $Z$ satisfy ($\sigma$) and $f$ satisfies ($\kappa_\rel$).
A \emph{functor} from a simplicial space $X$ over $Z$ to a fiberwise
complete Segal space $Y$ over $Z$ is just
a simplicial map $g\co X\to Y$ over $Z$.
\end{defn}

To make sense of $\rmap_Z(X,Y)$, space of derived simplicial maps from $X$ to $Y$ over $Z$,
 we use a model category structure on the category of
simplicial spaces over $Z$ (and at this point we should try forget the fact that $Y$ is Segal and fiberwise
complete over $Z$). The notion of weak equivalence is in any case degreewise. Other details are omitted.

\subsection{Local configuration categories and little disk operads} \label{subsec-locandlittle}
From the definitions and a result \cite[Lemma 4.1]{BoavidaWeiss} which expresses the locality of local
configuration categories, the space
\[ \rmap^\partial_\fin(\config^\loc(M;r),\config^\loc(N;r)) \]
in the diagram of theorem~\ref{thm-mainbdry} can be described as the space of sections of a
certain fibration $E\to M$. Namely, $E$ is the space of triples $(x,y,g)$ such that $x\in M$, $y\in M$ and
\[ g\in \rmap_\fin(\config(T_xM;r),\config(T_yN;r)); \]
the projection $E\to M$ takes $(x,y,g)$ to $x$.
More precisely, these are sections defined on all of $M$ and \emph{prescribed} over $\partial M$. The prescription
on $\partial M$ comes from the fact that a preferred embedding $\partial M\to \partial N$ has been selected.

\begin{thm} \cite[\S7]{BoavidaWeiss} \cite{Weiss_embsuppl} \label{thm-operpest}
Let $U$ and $V$ be finite dimensional real vector
spaces. Let $\sE_U$ and $\sE_V$ be the operads of little disks in $U$ and $V$, respectively.
There is a weak equivalence of derived mapping spaces
\[ \rmap(\sE_U,\sE_V) \lra  \rmap_\fin(\config(U),\config(V)). \]
Similarly, for an integer $r\ge 1$ there is a weak equivalence
\[ \rmap(\sE_{U,\le r},\sE_{V,\le r}) \lra \rmap_\fin(\config(U;r),\config(V;r)) \]
where $\sE_{U,\le r}$ is the variant of $\sE_U$ truncated at cardinality $r$; see remark~\ref{rem-ops} below.
Analogous statements hold for the rationalizations of $\sE_U,\sE_V$ and of $\config(U),\config(V)$
$($and their truncations$)$, provided $\dim(U)\ge 3$ and $\dim(V)\ge 3$.
\end{thm}
(Note that the configuration spaces $\emb(\uli k,U)$, which are the essential constituents of $\sE_U$ and $\config(U)$,
are simply connected if $\dim(U)\ge 3$, connected but not simply connected if $\dim(U)=2$ and homotopy discrete if $\dim(U)=1$.)

As a consequence, the space in the lower right-hand corner
of the diagram of theorem~\ref{thm-mainbdry} can be described as the space of sections, prescribed over $\partial M$, of a
fibration on $M$ whose total space is the space of triples $(x,y,f)$ where $x\in M$,
$y\in N$ and $f\in \rmap(\sE_{T_xM,\le r},\sE_{T_yN,\le r})$.

\begin{rem} \label{rem-ops} The operads in this theorem are plain operads in the category of spaces.
In order to make sense of spaces of derived maps between such operads, we need at least a notion of weak
equivalence for maps between such operads \cite{DwyKa}. This is \emph{levelwise} equivalence.
Cisinski and Moerdijk \cite{CisinskiMoerdijk1}, \cite{CisinskiMoerdijk2}, \cite{CisinskiMoerdijk3}
have developed a more tractable setting for operads in spaces by associating to
an operad $P$ its dendroidal nerve $N_dP$, a contravariant functor from a certain category
of trees to spaces. Briefly, the functor $N_d$ gives (in the case of plain operads)
a faithful translation so that $\rmap(P,Q)$ can be
identified with $\rmap(N_dP,N_dQ)$. And here the correct interpretation of $\rmap(N_dP,N_dQ)$ can be
obtained by using any of the familiar model structures with \emph{levelwise} weak equivalences on
the category of contravariant functors from that tree category to spaces. Also, the correct interpretation
of the \emph{truncation at level $r$} of an operad $P$ is obtained by restricting $N_dP$ to a certain
full subcategory of the tree category (consisting of those trees where every edge has at most $r$ branches,
a.k.a.~incoming edges).
\end{rem}

\begin{rem} \label{rem-ops2} The abstract relationship between the operad $\sE_V$ in theorem~\ref{thm-operpest}
and the topological category (or Segal space) $\config(V)$ is as follows. A monochromatic (= plain) topological operad $P$ determines
a category $P^\sharp$ enriched over spaces with object set $\NN=\{0,1,2,\dots\}$. This is the PROP associated with the
operad; see for example \cite[\S2.3]{Adams}. The space of morphisms from $a$ to $b$ in $P^\sharp$ is the disjoint union, taken over all maps
$f$ from $\uli a$ to $\uli b$, of the products
\[   \prod_{j=1}^b P(\,|f^{-1}(j)|\,)\,. \]
The category $P^\sharp$ comes with a monoidal product such that the induced monoid structure on the set of objects
is the ordinary addition in $\NN$. It is easy to recover the operad $P$ from $P^\sharp$ with the monoidal product.
But the process which takes us from $P$ to $P^\sharp$ \emph{without the monoidal product} and then further to the
comma category $P^\sharp\!\downarrow\! 1$ is a forgetful process. That forgetful process describes the relationship between $\sE_V$ and $\config(V)$:
if we take $P=\sE_V$ then the comma category $P^\sharp\!\downarrow\! 1$ is weakly equivalent (in the setting of Segal spaces) to
$\config(V)$. Theorem~\ref{thm-operpest} is a special case of a statement \cite[Thm. 7.5]{BoavidaWeiss} saying roughly that a forgetful map
of derived mapping spaces
\[   \rmap(P,Q) \lra  \rmap_\fin(P^\sharp\!\downarrow\! 1,Q^\sharp\!\downarrow\! 1) \]
is a weak equivalence if the spaces $P(0)$, $P(1)$, $Q(0)$ and $Q(1)$ are contractible.
\end{rem}

\subsection{Spaces of derived maps between certain topological operads} \label{subsec-derivmaps}
Let $P$ and $Q$ be (plain) operads
in the category of spaces. We ask for a practical description of the homotopy fibers $\Phi(r)$ of the forgetful map
\[
\xymatrix@M=8pt{ \rmap(P_{\le r},Q_{\le r}) \ar[r] &
\rmap(P_{\le  r-1},Q_{\le r-1})
}
\]
where $P_{\le r}$ for example means the truncation of $P$ at level $ r$ (where operations of arity greater than
$ r$ are suppressed). Note that there is a plural here: the map has many homotopy fibers, and we ought to
specify a point in the target
to specify one of them, and we do but the notation does not show it. ---
To make it easier, we assume that $P$ and $Q$ are both contractible in degrees 0 and 1. (Beware: for an
operad, being contractible in degree 0 is quite different from being empty in
degree 0. The operad $\sE_s$ of little disks in $\RR^s$ is the standard example of an operad which is contractible
both in degree 0 and 1.)

\smallskip The answer which follows, condensed in theorem~\ref{thm-Goeppl} below, is a quotation from \cite{Goeppl}.
[The reviewer observed: ``The notation in this section is not consistent with \cite{Goeppl}.'' It is impossible to
make it consistent with \cite{Goeppl} without causing new inconsistencies, but here
is a tiny dictionary. Where G\"oppl writes $\textup{Rhom}(-,-)$ it is $\rmap(-,-)$ here; he writes $E_n$ for the operad, or rather
the family of weakly equivalent operads, which is denoted $\sE_n$ here, and $\bar\Omega$ for a category which
is denoted $\tree^c$ here and in \cite{BoavidaWeiss}. Where G\"oppl and others write $\Omega(T)$ for a certain operad associated with a tree $T$,
its is $\omega(T)$ here.]

Following \cite{MoerdijkWeiss}, \cite{CisinskiMoerdijk1}, \cite{CisinskiMoerdijk2}, \cite{CisinskiMoerdijk3}
we associate to a (plain) operad $P$ in spaces its dendroidal nerve $N_dP$. As already mentioned
in remark~\ref{rem-ops}, this is a contravariant functor from a certain category of trees to spaces.
Because we are concerned exclusively with operads $P$ wich satisfy $P(0)=*$, we can get away
with a category of trees, notation $\tree^c$, which is smaller and easier to understand than the one used mainly in
\cite{MoerdijkWeiss}, \cite{CisinskiMoerdijk1}, \cite{CisinskiMoerdijk2}, \cite{CisinskiMoerdijk3},
whose name there is $\Omega$. (For the reduction see \cite[Lem.7.12]{BoavidaWeiss}; although this
makes a slightly bigger reduction, it is easy to adapt the proof.)

\begin{defn} An object of $\tree^c$ is a finite poset $S$ such that
\begin{itemize}
\item[-] for every $x\in S$ the set $\{y\in S~|~y\le x\}$ is totally ordered
(with the partial order induced from $S$);
\item[-] there exists $z\in S$ such that $z\le x$ for all $x\in S$. (This is the \emph{root}.)
\end{itemize}
(The elements of $S$ can be thought of as \emph{edges}.)
A morphism in $\tree^c$ from $S$ to $T$ is an order-preserving map $f\co S\to T$ which satisfies the
following condition: if $x$ and $y$ are unrelated elements in $S$ (so that neither $x\le y$ nor $y\le x$ holds),
then $f(x)$ and $f(y)$ are also unrelated. (\emph{Warning}: suppose that $f$ satisfies the
condition and $x,y$ are elements of $S$ such that $f(x)\le f(y)$. It does not follow that
$x\le y$, because it can happen that $f(x)=f(y)$ and $y<x$.)
\end{defn}

(Moerdijk uses the expression \emph{closed tree} for objects of $\tree^c$ and writes $\bar\Omega$
instead of $\tree^c$.)

\begin{defn} For an object $S$ in $\tree^c$ and $x\in S$, let $\inc(x)$ be the set of edges just above $x$,
that is, the set of minimal elements in $\{y\in S~|~y>x\}$. Let $\tree^c_{\le r}$ be the full
subcategory of $\tree^c$ spanned by the objects $S$ such that $|\inc(x)|\le r$ for every $x\in S$.
\end{defn}

\begin{defn} For every $r\ge 0$ there is an object $t_r$ of $\tree^c$ characterized up to isomorphism
by the following property: $|\inc(z)|=r$ if $z$ is the root, and $|\inc(z)|=0$ in all
other cases. This is called the \emph{$r$-corolla}. It belongs to $\tree^c_{\le r}$ but not to $\tree^c_{\le r-1}$\,.
\end{defn}

\begin{defn} \label{defn-grafting} Let $S$ be an object of $\tree^c$ and let $f\co t_r\to S$ be any morphism in $\tree^c$.
Think of this as a selection of one element $z\in S$ (corresponding to the root of $t_r$) and
further elements $y_1,y_2,\dots, y_r$ in $S$, all $\ge z$ and
\emph{pairwise unrelated} in the poset $S$. For every $j\in \{1,2,\dots,r\}$, let $g_j\co t_{s_j}\to S$ be a morphism
taking the root to $y_j$. Let $q=s_1+s_2+\cdots+s_r$. Then there is a morphism $h\co t_q\to S$ obtained from $f$
and the $g_j$ by \emph{grafting}. To make this explicit we identify $\{1,2,\dots,q\}$ with the coproduct
of the sets $\{1,2,\dots,s_j\}$ for $j=1,2,\dots,r$, and specify $h(x)=g_j(x)$ if $x$ belongs to the $j$-th summand.
\end{defn}

For a (plain) operad $P$ in spaces,
and a finite set $U$ with $r$ elements, we may write $P(U)$ to mean
$P(r)\times_{\Sigma_r} (\textup{set of bijections $U\to \uli r$})$.
In particular $P(\uli r)$ is canonically identified with $P(r)$.

\smallskip
For an operad $P$ in spaces which has $P(0)=*$, the value $(N_dP)_S$ of the nerve $N_dP$ at an object $S$ in $\tree^c$ is the product
$\prod_{x\in S} P(\inc(x))$. Elements of $(N_dP)_S$ are sometimes called $P$-\emph{decorations} of the tree $S$.
This description of $N_dP$ does not explain very well how $N_dP$ is a contravariant functor
on $\tree^c$. The following alternative description does explain this, but it is a little more advanced.
We observe that every object $S$ of $\tree^c$ determines a \emph{multi-colored} operad $\omega(S)$
whose objects (alias colors) are the elements of $S$. Operations in $\omega(S)$, with an ordered list
$y_1,\dots,y_r$ of $r$ sources and one target $z$,
correspond to morphisms $t_r\to S$ in $\tree^c$ taking the root to $z$ and the other edges to $y_1,\dots,y_r$, respectively.
Composition of operations is by grafting as in definition~\ref{defn-grafting}.
(By construction, $\omega(P)$ does not have any operations with repeated sources.)
Now $(N_dP)_S$ can be defined as the space of morphisms from $\omega(S)$ to $P$. (Such a morphism
from $\omega(S)$ to $P$ must take all objects of $\omega(S)$ to the same object of $P$, since $P$ has only one object.)
A morphisms $S\to T$ in $\tree^c$ induces a map of operads $\omega(S)\to \omega(T)$ and then a map of spaces
$(N_dP)_T\to (N_dP)_S$. ---
Relationship between the two descriptions of $(N_dP)_S$: every $x\in S$ determines a ``subtree'' $S_x$ of $S$ consisting
of $x$ and the elements of $\inc(x)$. By restricting morphisms $\omega(S)\to P$ along the inclusion $\omega(S_x)\to \omega(S)$,
we obtain elements in $P(\inc(x))$. In this way, the space of morphisms $\omega(S)\to P$ turns out to be a product of
factors $P(\inc(x))$, for $x\in S$.

\begin{defn} \label{defn-boundcobound} For an operad $P$ in spaces which has $P(0)=*$, let
\[  \bound_r(P): = \hocolimsub{\twosub{t_r\to S}{S\textup{ in }\tree^c_{\le r-1}}} (N_dP)_S \]
where $t_r$ is the $r$-corolla. (The homotopy colimit is indexed by the ``under'' category $t_r\downarrow F_{r-1}$
where $F_{r-1}\co \tree^c_{\le r-1}\to \tree^c$ is the inclusion functor.)
Similarly let
\[ \cobound_r(P):= \holimsub{\twosub{S\to t_r}{S\textup{ in } \tree^c_{\le r-1}}} (N_dP)_S\,. \]
(The homotopy limit is indexed by the ``over'' category $F_{r-1}\downarrow t_r$.)
The morphisms $t_r\to S$ and $S\to t_r$ for $S$ in the fine print of the definitions of $\bound_r(P)$
and $\cobound_r(P)$ determine canonical maps
\[  \bound_r(P) \lra P(r) \lra \cobound_r(P). \]
These should be viewed as maps between spaces with an action of the symmetric group $\Sigma_r$
(the automorphism group of $t_r$).
\end{defn}

\begin{prop} \label{prop-bound} \cite[Ex.1.1.6, Ex.2.1.13]{Goeppl}
If $P=\sE_m$ then $P(r)$ is weakly equivalent to $\emb(\uli r,\RR^m)$ and the map
$\bound_r(P)\to P(r)$
is the Axelrod-Singer $($-Fulton-MacPherson$)$ boundary inclusion, up to weak equivalence of maps. \qed
\end{prop}
\proof[Sketch proof] (The sketch proof of this in \cite{Goeppl} gives the right idea, but it is too negligent
about the indexing categories.) Let $\texttt{C}_0=t_r\downarrow F_{r-1}$ and let $\Phi_0$ be the contravariant
functor taking $t_r\to S$ in $\texttt{C}_0$ to $(N_dP)_S$\,, so that $\bound_r(P)= \hocolim~\Phi_0$.
The category $\texttt{C}_0$ has a full subcategory $\texttt{C}_1$ spanned by the objects $u\co t_r\to S$
of $\texttt{C}_0$ such that all maximal elements of the poset $S$ are in the image of $u$. Let $\Phi_1$ be the
restriction of $\Phi_0$ to $\texttt{C}_1$. The inclusion
$\texttt{C}_1\to \texttt{C}_0$ admits a right adjoint. This implies that $(\texttt{C}_1)^\op$ is \emph{homotopy terminal}
in $(\texttt{C}_0)^\op$. It follows that the inclusion $\hocolim~\Phi_1\to \hocolim~\Phi_0$ is a weak equivalence. Next,
the category $\texttt{C}_1$ has a subcategory $\texttt{C}_2$ consisting of the objects $t_r\to S$ where $|\inc(x)|\ne 1$
for all $x\in S$. Let $\Phi_2$ be the restriction of $\Phi_1$ to $\texttt{C}_2$.
The inclusion $\texttt{C}_2\to \texttt{C}_1$ has a left adjoint
$\alpha$ and the unit maps $(t_r\to S)\to \alpha(t_r\to S)$ for objects $t_r\to S$ in $\texttt{C}_1$
induce weak equivalences $\Phi_1\alpha(t_r\to S)\to \Phi_1(t_r\to S)$. It follows that the inclusion
$\hocolim~\Phi_2\to \hocolim~\Phi_1$ is a weak equivalence. (The conversion from $\hocolim~\Phi_0$
to $\hocolim~\Phi_2$ did not use the assumption $P=\sE_m$; it used only $P(0)=*$ and $P(1)\simeq *$.) ---
Now we can follow the reasoning
of \cite[Ex.2.1.13]{Goeppl}. It is convenient to use a specific model of $P=\sE_m$
as in \cite[Ex.1.1.6]{Goeppl}. In this setting the functor $\Phi_2$ is a \emph{cofibrant} object
in a suitable model structure on the category of contravariant functors from $\texttt{C}_2$ to spaces.
It follows that the projection $\hocolim~\Phi_2\to \colim~\Phi_2$ is a weak equivalence. This completes the sketch proof
because $\colim~\Phi_2$ is exactly the Axelrod-Singer-Fulton-MacPherson boundary. More precisely,
the canonical map from $\colim~\Phi_2$ to $(N_dP)_{t_r}=P(r)$ is exactly the A.S.F.M.~boundary inclusion.  \qed

\smallskip
In the definition of $\cobound_r(P)$, the category $F_{r-1}\downarrow t_r$
has a small subcategory $\sA_r$ which is isomorphic to the poset of proper subsets of $\uli r$.
Namely, we may identify the set of non-root edges of $t_r$ with $\uli r$. Every proper subset $U$ of
$\uli r$ determines a ``sub''-tree of $t_r$ (and an inclusion morphism) consisting of the root
and the edges in $U$. --- The inclusion $\sA_r\to  F_{r-1}\!\downarrow\!t_r$ has a left adjoint. This leads to the
following.

\begin{lem} \label{lem-cobound} The restriction map
\[ \cobound_r(P)~~=\!\!\!\!\!\holimsub{\twosub{S\to t_r}{S\textup{ in } \tree^c_{\le r-1}}} \!\!\!\!(N_dP)_S
\quad\lra\quad \holimsub{(S\to t_r)\textup{ in }\sA_r} (N_dP)_S \]
is a weak equivalence. \qed
\end{lem}

\begin{thm}
\cite{Goeppl} \label{thm-Goeppl} The following self-explanatory commutative square is homotopy Cartesian:
\[
\xymatrix{
\rmap(P_{\le r},Q_{\le r}) \ar[d] \ar[r] & {\rmap^{\Sigma_r}(P(r),Q(r))} \ar[d]  \\
{\rmap(P_{\le r-1},Q_{\le r-1})} \ar[r] & {\rmap^{\Sigma_r}(\bound_r(P)\to P(r),~Q(r)\to \cobound_r(Q))}
}
\]
\end{thm}
\emph{Interpretation}: each homotopy fiber in the left-hand column can be described as a space of fillers
(with $\Sigma_r$ invariance):
\[
\xymatrix{
\bound_r(P) \ar[d] \ar[r] \ar[dr] & P(r) \ar@{..>}[d] \ar[dr] \ar[r] & \cobound_r(P) \ar[d] \\
\bound_r(Q) \ar[r] &   Q(r) \ar[r] & \cobound_r(Q )
}
\]
\emph{Remark}: G\"oppl writes $\rmap_{\Sigma_r}$ where we have $\rmap^{\Sigma_r}$. \newline
\emph{Remark}: In the upper right-hand term of the diagram of theorem~\ref{thm-Goeppl}, the spaces $P(r)$ and $Q(r)$
are understood to be functors from $\Sigma_r$ to spaces. (A group is a category.) In the
category of such functors, the notion of weak equivalence that we want here is \emph{levelwise}
weak equivalence; a $\Sigma_r$-map $X\to Y$ of $\Sigma_r$-spaces is a weak equivalence, by definition, if
the underlying map of spaces is a weak equivalence. This determines the meaning of the
derived mapping space $\rmap^{\Sigma_r}(P(r),Q(r))$, as in \cite{DwyKa}. For a ``computation''
in a model category setting, we look for a cofibrant replacement of $P(r)$,
for example (if the projective model structure is used)
a $\Sigma_r$-CW-space $X$ with free $\Sigma_r$-cells and a $\Sigma_r$-map $X\to P(r)$ which is a weak equivalence as a map
of spaces. (We also need a fibrant replacement $Q(r)\to Y$, but ``fibrant'' is a weak condition in the projective
model structure.) Similarly, in the lower right-hand term of the diagram, we are looking at two functors from
$[1]\times\Sigma_r$ to spaces, where $[1]$ is the totally ordered set $\{0,1\}$. (An ordered set is a category.)
In the category of such functors, the notion of weak equivalence that we want here is again levelwise weak equivalence.

\medskip

\section{Configuration categories and homotopy automorphisms} \label{sec-conhoto}
\subsection{The case of unrestricted cardinalities} \label{subsec-unrescar}
Let $K$ be a smooth compact manifold with boundary and let $g\co K\to K$ be a homeomorphism. The
homeomorphism induces a functor of configuration categories,
\[  \config(K\smin\partial K)\to \config(K\smin\partial K). \]
Can we recover the \emph{homotopy} automorphism of the pair $(K,\partial K)$ determined by $g$
from the above functor of configuration categories?

The question can be sharpened in more than one way. Let $K_-:= K\smin \partial K$. The boundary of $K$
can be recovered from $K_-$ in a homotopical sense as the space $\partial^h K_-$ of proper maps from a half-open interval
$[0,1)$ to $K_-$\,. Indeed the pair $(K,\partial K)$ can be recovered from $K_-$ in a homotopical sense as the map
\begin{equation}  \label{eqn-hbdry} \partial^h K_- \lra K_- \end{equation}
given by evaluation at $0\in [0,1)$. --- Therefore we reconsider; let $g\co K_-\to K_-$ be a homeomorphism. The
homeomorphism induces a functor of configuration categories,
\[  \config(K\smin\partial K)\to \config(K\smin\partial K). \]
Can we recover the \emph{homotopy} automorphism of the map~\eqref{eqn-hbdry} determined by $g$
from the above functor of configuration categories? (A map of spaces can be viewed as a functor from $\{0,1\}$ to spaces, where $\{0,1\}$
has the usual ordering.) To formulate an even more precise question, we
introduce $\haut(\partial^h K_-\to K_-)$, the space of homotopy
automorphisms of~\eqref{eqn-hbdry}.
Let $\haut_{\fin}(\config(K\smin\partial K))$ be the union of the weakly invertible path components
of $\rmap_{\fin}(\config(K\smin\partial K))$.
We look for an arrow making the following diagram
commutative (up to specified homotopies):
\begin{equation} \label{eqn-sentries1}
\begin{aligned}
\xymatrix@R=18pt@C=20pt@M=5pt{
B\homeo(K_-) \ar[d] \ar@/^1pc/[dr]   \\
B\haut(\partial^h K_- \to K_-) \ar[d]^-{\textup{res.}} & \ar@{..>}[l] B\haut_{\fin}(\config(K_-)) \ar[d] \\
B\haut(K_-) & \ar[l]_-\simeq  B\haut_{N\fin}(\config(K_-;1))
}
\end{aligned}
\end{equation}
See remark~\ref{rem-commu} below.

This problem has a trivial solution when $\partial K$ is empty. In the general case,
the term $B\haut(\partial^h K_- \to K_-)$ can also be described (up to weak equivalence) as the classifying
space of the topological monoid or $A_\infty$-space of
\emph{proper} homotopy automorphisms of $K_-$ (homotopy automorphisms in the world of locally compact
spaces and proper maps). But this does not make it easier; the word \emph{proper} does not translate
well into the language of configuration categories. Therefore the cases where $\partial K$ is nonempty are
harder. Here is one rather severe condition on $K$ which is sufficient (to guarantee the existence of the broken arrow):
\begin{itemize}
\item[] \emph{$K$ is the total space of a smooth fiber bundle on a closed smooth manifold
$L$ where the fibers are smooth manifolds homeomorphic to a disk $D^c$ for some $c\ge 3$.}
\end{itemize}
See \cite{TillWeissSen1} and especially \cite[Thm 4.2.1]{WeissSen2}. 
A weaker but still sufficient condition is given in \cite{Tillmann}. Roughly stated, the weaker condition of \cite{Tillmann} is that
$K$ be a thickening of a simplicial complex embedded smoothly in $K\smin\partial K$, of codimension
$c\ge 3$.

There is a more complicated version for a smooth compact $K$ with a codimension zero smooth compact
submanifold $\partial_0K$ of $\partial K$. Let $\partial_1K$ be the closure of $\partial K\smin\partial_0K$
in $\partial K$. \emph{For the purposes of this section} and especially diagram~\eqref{eqn-sentries2}
just below, we redefine $K_-$ as $K\smin \partial_1K$\,. This is a manifold with
boundary, the boundary being $\partial K_-=\partial_0 K\smin \partial_1K$. Let $\partial^h_1K_-$ be the space of proper maps from $[0,1)$ to $K_-$\,.
This contains a copy of $\partial^h(\partial K_-)$, the space of proper maps from $[0,1)$ to $\partial K_-$.
We look for an arrow making the following diagram
commutative (up to specified homotopies):
\begin{equation} \label{eqn-sentries2}
\begin{aligned}
\xymatrix@R=25pt@C=16pt@M=5pt{
B\homeo_\partial (K_-) \ar[d] \ar@/^1.5pc/[dr]   \\
{B\haut_{\textup{ row 1}}\left( \begin{aligned} \begin{matrix} \partial^h(\partial K_-) & \!\!\to\!\! & \partial K_- \\
\downarrow && \downarrow \\
\partial^h_1K_- & \!\!\to\!\! & K_- \end{matrix} \end{aligned} \right)} \ar[d] & \ar@{..>}[l] B\haut^\partial_\finplus(\config(K_-)) \ar[d] \\
B\haut_\partial(K_-) & \ar[l]_-\simeq  B\haut^\partial_\finplus(\config(K_-;1))
}
\end{aligned}
\end{equation}
(\emph{Notation}: the configuration category $\config(K_-)$ has an over category, a.k.a. comma category,
associated with the unique object taken to $[0]$ in $\finplus$ by the reference functor. The notation
$\haut$ with superscript ${}^\partial$ and subscript ${}_\finplus$ is for homotopy automorphisms which cover the identity of
$\finplus$ and induce the identity automorphism of that over category, in a derived sense. In the left-hand column,
$\haut_{\textup{ row }1}(-)$ means: the space of homotopy automorphisms of the square diagram which restrict to
the identity on the first row of the diagram. More precisely, the square diagram is a functor $F$ from $[1]\times[1]$ to spaces
and we are interested in the homotopy fiber of $\haut(F)\to \haut(F|_{\{0\}\times[1]})$. And of course, $\haut(F)$ is short
for the union of invertible path components of $\rmap(F,F)$. The preferred notion of weak equivalence in the
category of functors from $[1]\times[1]$ to spaces is once again levelwise weak equivalence.)

Despite appearances, the lower horizontal arrow in~\eqref{eqn-sentries2} is not completely obvious.
Let $Y=K\smin\partial K$ be the space of objects of $\config(K_-;1)$ taken to the object $[1]$ of $\finplus$,
and let $X$ be the space of morphisms in $\config(K_-;1)$ taken to the unique morphism $[1]\to [0]$ in $\finplus$.
There is a map $s\co X\to Y$ (taking a morphism to its source). Now there is a forgetful homomorphism of $A_\infty$-spaces
\[  \haut^\partial_\finplus(\config(K_-;1)) \lra \hofiber[\,\haut(s\co X\to Y) \lra \haut(X)\,] \]
which is also a weak equivalence.
Therefore we should (at least) specify a weak equivalence or a zigzag of weak equi\-valences relating
$s:X\to Y$ to the inclusion $\partial(K_-)\to K_-$\,.
To do so we introduce $X_1$~, the space of all maps $[0,1]\to K_-$
taking $0$ to some point in $\partial(K_-)$~, and $X_2$~, the space of all maps $[0,1]\to K_-$ taking $0$
and only $0$ to some point in $\partial(K_-)$. Then we have a commutative diagram
\[
\xymatrix@R=12pt@C=10pt{
X \ar[d] &\ar[l] X_2 \ar[d] \ar[r] & X_1 \ar[d] & \ar[l] \partial(K_-) \ar[d] \\
Y  & \ar@{=}[l] Y \ar[r] & K_- \ar@{=}[r] & K_-
}
\]
in which all horizontal arrows are weak equivalences.

\smallskip
Diagram~\eqref{eqn-sentries2} reduces to diagram~\eqref{eqn-sentries1} if $\partial_0K$ is empty.

\smallskip
Again, the problem (finding the broken arrow in diagram~\eqref{eqn-sentries2})
has a trivial solution if $\partial_1 K$ is empty. In the case of nonempty $\partial_1K$,
a sufficient but severe condition
is roughly as follows: there exists a smooth compact
submanifold $L$ of $K\smin\partial_1 K$, neatly embedded
with boundary $\partial L\subset \partial_0 K\smin\partial_1 K$, and of codimension $c\ge 3$ in $K$, such that
$K$ is a thickening (a.k.a.~regular neighborhood) of $\partial_0K\cup L$. See \cite[Thm.3.2.1]{TillWeissSen1} and \cite[Thm.5.3.1]{WeissSen2} for the details.
This will take care of the case that we are after, although \cite{Tillmann} gives a weaker condition
which is still sufficient.

\smallskip
The case that we are after
is $K=W_{g,1}$ and $\partial_1K$ equal to a closed hemisphere of $\partial W_{g,1}\cong S^{2n-1}$
containing the selected point $z$\,.
(Then $K_-=K\smin\partial_1K$ is diffeomorphic to $W_\parpun$. More precisely the inclusion of
$K\smin\partial_1K$ in $W_\parpun$ is smoothly isotopic to a diffeomorphism, as a smooth embedding taking
boundary to boundary.) The severe condition is satisfied in this case, provided $n\ge 3$. Namely, for $L$ we can take a
disjoint union of $2g$ neatly embedded $n$-disks in $K\smin \partial_1K$. (Therefore $c$ is $2n-n=n$.)
The restriction map
\[ \homeo_\partial (K)\to \homeo_\partial(K_-) \]
is a weak equivalence (exercise, relying on the Alexander trick), and
\[
B\haut_{\textup{ row 1}}\left( \begin{aligned} \begin{matrix} \partial^h(\partial K_-) & \!\!\to\!\! & \partial K_- \\
\downarrow && \downarrow \\
\partial^h_1K_- & \!\!\to\!\! & K_- \end{matrix} \end{aligned} \right)~~\simeq~~ B\haut_\partial(K)
\]
by an easier form of the Alexander trick for homotopy automorphisms of disks.
[Here is some notation which we will need in a moment:
for a manifold $M$ with boundary, let $\haut^\str_\partial(M)\subset \haut_\partial(M)$ be the
submonoid whose elements are the $f\in \haut_\partial(M)$ which satisfy $f(M\smin\partial M)\subset M\smin\partial M$.]
Therefore diagram~\eqref{eqn-sentries2} in this instance simplifies to
\begin{equation} \label{eqn-sentries2inst}
\begin{aligned}
\xymatrix@R=18pt{
 B\homeo_{\partial}(W_{g,1})  \ar@/^1pc/[dr] \ar[d]  \\
 B\haut^\str_{\partial}(W_{g,1}) \ar[d]^-{\textup{res.}} &
 \ar@{..>}[l] B\haut^\partial_\finplus(\config(W_\parpun)) \ar[d]^-{\textup{res.}}  \\
 B\haut^\str_\partial(W_\parpun) &  \ar[l]_-\simeq   B\haut^\partial_\finplus(\config(W_\parpun;1))
}
\end{aligned}
\end{equation}
Moreover we have $\haut^\str_\partial(W_\parpun)\simeq \haut_\partial(W_\parpun)\simeq \haut_*(W_\parpun)$.

\begin{rem} \label{rem-commu} In diagram \eqref{eqn-sentries1}, the outer
diagram (obtained by deleting the broken arrow) is already commutative up to a preferred homotopy.
We look for a solution (consisting of
the broken arrow and additional homotopies or higher homotopies) respecting this information.
A precise formulation is as follows. Diagram~\eqref{eqn-sentries1} can be shortened to a diagram
\[
\xymatrix@R=25pt@C=-5pt@M=5pt{
& B\homeo(K_-) \ar[dl]_-e \ar[dr]^-{\varphi_1} &   \\
B\haut(\partial^h K_-\to K_-) \ar[dr]_-{\varphi_2} && \ar@{..>}[ll] B\haut_{\fin}(\config(K_-)) \ar[dl]^-f \\
& B\haut(K_-) &
}
\]
That diagram (without the edges $e$ and $f$) defines a map
\begin{equation} \label{eqn-thanksTN}
\xymatrix@M=5pt{
\map\big(B\haut_{\fin}(\config(K_-)),\,B\haut(\partial^h K_-\to K_-)\big) \ar[r] &
\rmap(\varphi_1,\varphi_2).
}
\end{equation}
(The target space is a space of derived natural transformations between two functors from $\{0,1\}$ to spaces,
where $\{0,1\}$ has the usual poset structure. A poset can be viewed as a category. The functors are $\varphi_1$ and
$\varphi_2$.) What we should really look
for is this: \emph{a point in the homotopy fiber of~\eqref{eqn-thanksTN} over the point of $\rmap(\varphi_1,\varphi_2)$
determined by $e$, $f$ and that preferred homotopy.} Similar observations can be made in regard to diagrams~\eqref{eqn-sentries2}
and the special case in \eqref{eqn-sentries2inst}. They also apply to
the truncated variant in section~\ref{subsec-rescar} below and the special case~\eqref{eqn-Wsentriestrunc}. \newline
The main results of \cite{WeissSen2} are not stated with all the details that we need here. For example
\cite[Thm 4.2.1]{WeissSen2} provides the broken arrow in~\eqref{eqn-sentries1} and a homotopy making the upper triangle
commutative (up to that homotopy). It speaks of an $A_\infty$-action (by homotopy automorphisms)
of $\haut_{\fin}(\config(K_-))$ on the pair
$(K,\partial K)$, or equivalently on the map $(\partial^h K_-\to K_-)$, which extends a better known
$A_\infty$ action of $\homeo(K_-)$ on the same
pair or map. It fails to point out that the underlying $A_\infty$-action
of $\haut_{\fin}(\config(K_-))$ on the space $K_-$ is the tautological action of
$\haut_{\fin}(\config(K_-;1))\simeq \haut(K_-)$
on $K_-$, pulled back by means
of the forgetful $A_\infty$-map from $\haut_{\fin}(\config(K_-))$ to $\haut_{\fin}(\config(K_-;1))$.
\end{rem}

\begin{rem} The reviewer has asked: why is it relevant that $K$ is smooth, e.g. in diagram~\eqref{eqn-sentries1}~?
It is probably not relevant at all, but the statements from \cite{TillWeissSen1} and \cite[Thm 4.2.1]{WeissSen2} that we have
used here are for smooth manifolds. One excuse is that the proofs there rely on manifold calculus, which is well
established in a setting of smooth manifolds but not in a setting of topological manifolds.
\end{rem}

\subsection{The case of restricted cardinalities.} \label{subsec-rescar} These statements have truncated versions.
We truncate the configuration categories by allowing only configurations of cardinality $\le r$, say.
Under the usual severe condition on $K$ and $\partial_0K\subset \partial K$,
there is an arrow making the following diagram commutative up to specified homotopies:
\[       
\xymatrix@R=22pt@C=20pt@M=5pt{
B\homeo_\partial (K_-) \ar[d] \ar@/^1.5pc/[dr]   \\
{\tau_{\le b}\left(B\haut_{\textup{ row 1}}\left( \begin{aligned} \begin{matrix} \partial^h(\partial K_-) & \!\!\to\!\! & \partial K_- \\
\downarrow && \downarrow \\
\partial^h_1K_- & \!\!\to\!\! & K_- \end{matrix} \end{aligned}
\right)\!\!\right)} \ar[d] & \ar@{..>}[l] B\haut^\partial_\finplus(\config(K_-;r)) \ar[d] \\
\tau_{\le b}\big(B\haut_\partial(K_-)\big)  & \ar[l]_-\simeq  \tau_{\le b}\big(B\haut^\partial_\finplus(\config(K_-;1))\big)
}
\]
The prefix $\tau_{\le b}$ means that we kill homotopy groups (by Postnikov truncation) in
dimensions $>b$. The integer $b$ depends on $r$. It is advisable to use a functorial Postnikov truncation as in \cite{DwyKaIndag84}.
The estimates in \cite[Thm 5.4.1]{WeissSen2} show that taking
$b$ equal to $(r+1)(c-2)-(\dim(K)-1)$ is safe,  
where $c$ is the codimension of $L$ in $K$, from the condition which we impose on $K$. (Apologies as in
remark~\ref{rem-commu} apply; the reference to \cite[Thm 5.4.1]{WeissSen2} only justifies the upper triangle
in the diagram, not the lower square which is much easier. Note that
\cite[Thm 5.4.1]{WeissSen2} uses a relative Postnikov truncation, notation $\wp_{(j+1)(c-2)}$ there, applied to
the map
\[ \partial^h_1K_-\to K_-~, \]
there $\partial^h_1M_-\to M_-$. This is a less forgetful
construction than the truncation $\tau_{\le b}$ with $b=(r+1)(c-2)-(\dim(K)-1)$ which we apply here in the first
instance to the middle term in the left-hand column \emph{after} making a $B\haut$ construction.)

In the situation of diagram~\eqref{eqn-sentries2inst} 
we obtain a diagram
\begin{equation} \label{eqn-Wsentriestrunc}
\begin{aligned}
\xymatrix@R=18pt@M=6pt{
 B\homeo_{\partial}(W_{g,1})  \ar@/^1pc/[dr] \ar[d]  \\
{\tau_{\le b}\left(B\haut^\str_{\partial}(W_{g,1})\right)} \ar[d]^-{\textup{res.}} &
 \ar@{..>}[l] B\haut^\partial_\finplus(\config(W_\parpun)) \ar[d]^-{\textup{res.}}  \\
{\tau_{\le b}\left(B\haut^\str_\partial(W_\parpun)\right)} & \ar[l]_-\simeq
\tau_{\le b}\left(B\haut^\partial_\finplus(\config(W_\parpun;1))\right)
}
\end{aligned}
\end{equation}
(commutative up to preferred
homotopies) where
\[ b:=(r+1)(n-2)-(2n-1)=(r-1)n-2r-1\,. \]
Indeed the dimension of $W_{g,1}$ is $2n$ and for $c$ we can take $n$. That was already
pointed out just before diagram \eqref{eqn-sentries2inst}.

\subsection{An application}
Fix a smooth compact manifold $K$ as in diagram~\eqref{eqn-sentries2}.
Suppose that the inclusion-induced map $\pi_0(\partial_1K)\to \pi_0K$ is surjective; equivalently, every connected component of
$K\smin\partial_1 K$ is noncompact.

Take $M=N=K\smin \partial_1K=: K_-$ in theorem~\ref{thm-mainbdry}. Instead of
$T_1\emb_\partial(K_-,\, K_-)$ we may write
$\fimm_\partial(K_-,\, K_-)$ as in definition~\ref{defn-fimm}. Let
\[ \fimm^\times_\partial(K_-,\, K_-)\subset \fimm_\partial(K_-,\, K_-) \]
consist of the elements whose
underlying map (from $K_-$ to itself) is a homotopy equivalence. Composing with
\[  \diff_{\partial}(K_-) \hookrightarrow
\emb_\partial(K_-,\,K_-) \to T_r\emb_\partial(K_-,\,K_-), \]
we obtain (from the square in theorem~\ref{thm-mainbdry}) a commutative square
\begin{equation} \label{eqn-sentriesappl1}
\begin{aligned}
\xymatrix{
B\diff_\partial(K_-) \ar[r] \ar[d] &  \ar[d]^-{\textup{specialization}} B\haut_\finplus^\partial(\config(K_-;r)) \\
B\fimm^\times_\partial(K_-,\,K_-) \ar[r] & B\haut_\fin^\partial(\config^\loc(K_-;r))
}
\end{aligned}
\end{equation}   
We wish to compare this with the commutative square
\begin{equation} \label{eqn-sentriesappl2}
\begin{aligned}
\xymatrix@C=24pt@M=6pt{
B\diff_\partial(K_-) \ar[d]^-{\textup{incl.}} \ar[r] & \ar[d]
{\tau_{\le b}\left(B\haut_{\textup{ row 1}}\left( \begin{aligned} \begin{matrix} \partial^h(\partial K_-) & \!\!\to\!\! & \partial K_- \\
\downarrow && \downarrow \\
\partial^h_1K_- & \!\!\to\!\! & K_- \end{matrix} \end{aligned} \right)\right)} \\
B\fimm^\times_\partial(K_-,\,K_-)  \ar[r]^-{\textup{forget}} &  {\tau_{\le b}\left(B\haut_{\partial}(K_-)\right)}
}
\end{aligned}
\end{equation}
where $b$ depends on $r$ as specified in section~\ref{subsec-rescar}.
Both diagrams,~\eqref{eqn-sentriesappl1} and~\eqref{eqn-sentriesappl2}, can be viewed as covariant functors
from the poset $A:=\{1,2\}^2$ to the category of spaces, say $F_1$ and $F_2$. (Use matrix style indexing,
so that $(i,j)$ marks the position in row $i$ and column $j$.) Let $C$ be the full sub-poset
of $A$ with elements $(1,1)$ and $(2,1)$. It is a consequence of section~\ref{subsec-rescar} that
we can make a derived natural transformation
\begin{equation} \label{eqn-sentriesnu1} \nu\co F_1\to F_2 \end{equation}
which extends the identity
$F_1|_C\to F_2|_C$ (on first columns). More formally stated, this produces an element in
\begin{equation} \label{eqn-sentriesnu2}
\nu\in  \hofiber_\id\left[\,\rmap(F_1,F_2) \to \rmap(F_1|_C,F_2|_C)\,\right]. \end{equation}
The new arrow $\nu_{(1,2)}$ is the broken arrow
in the (first) diagram of section~\ref{subsec-rescar}. The new arrow $\nu_{(2,2)}$ is not really new;
it is the lower horizontal arrow in the same diagram, pre-composed with a forgetful map.

Now specialize as in diagram~\eqref{eqn-Wsentriestrunc} and make the simplifications proposed there.
A trivialization of the tangent bundle $TW_\parpun$ or $TW_{g,1}$ can be used as
in proposition~\ref{prop-unravelBoaWe}. Then diagram~\eqref{eqn-sentriesappl1} turns into
\begin{equation} \label{eqn-sentriesappl3}
\begin{aligned}
\xymatrix@C=14pt{
B\diff_\partial(W_\parpun) \ar[r] \ar[d] &  \ar[d]^-{\textup{specialization}}
B\haut^\partial_\finplus(\config(W_\parpun;r)) \\
{B\left(\!\!\!\begin{aligned} \begin{array}{l}  \haut_*(W_{g,1})\,\,\, \ltimes \\
\map_*(W_{g,1},\SOr(2n))
\end{array}\end{aligned}\!\!\!\right)} \ar[r] & {B\left(\!\!\!\begin{aligned}
\begin{array}{l}  \haut_*(W_{g,1}) \,\,\, \ltimes \\
\map_*(W_{g,1},\haut_\fin(\config(\RR^{2n};r)))
\end{array} \end{aligned} \!\!\!\right)}
}
\end{aligned}
\end{equation}
and diagram~\eqref{eqn-sentriesappl2} turns into
\begin{equation} \label{eqn-sentriesappl4}
\begin{aligned}
\xymatrix{
B\diff_\partial(W_\parpun) \ar[r] \ar[d] &  \ar[d]
{\tau_{\le b}\left(B\haut_{\partial}(W_{g,1})\right)} \\
{B\left(\begin{aligned} \begin{array}{l}  \haut_*(W_{g,1})\,\,\, \ltimes \\
\map_*(W_{g,1},\SOr(2n))
\end{array}\end{aligned}\right)} \ar[r] & {\tau_{\le b}\left(B\haut_*(W_{g,1})\right)}
}
\end{aligned}
\end{equation}
where $b=(r-1)n-2r-1$.

\smallskip
The notation here deviates slightly from the notation in proposition \ref{prop-unravelBoaWe} and thereabouts; no careful distinction
is made between $\haut_*(W_{g,1})$ and $\haut_*(W_\parpun)$, or between $\map_*(W_{g,1},-)$ and
$\map_*(W_\parpun,-)$.

The natural transformation $\nu$ of~\eqref{eqn-sentriesnu1} and~\eqref{eqn-sentriesnu2}, specialized to be a natural map from
diagram~\eqref{eqn-sentriesappl3} to diagram~\eqref{eqn-sentriesappl4}, is what we need from this section.
It looks like a small reward for a huge effort.
But notice how it is related to section~\ref{sec-homhirz}.

\section{Dissonance} \label{sec-diss}
In this section, we combine results from the previous sections to make some
constructions leading to a proof of theorem~\ref{thm-gmain}.

\subsection{Preview} An outline of the section is as follows. We begin with the diagram
\begin{equation}
B\diff_\partial(W_{g,1}) \lra B\diff_\partial(W_\parpun) \lra B\homeo_{\partial}(W_{g,1}).
\end{equation}
The map on the left is induced by restriction of diffeomorphisms,
and the other is induced by one-point compactification.
We pass to rational homology and obtain a diagram of rational vector spaces
\[
\xymatrix@C=15pt{
H_{2n+4k}(B\diff_\partial(W_{g,1}))
\ar[r]^-{\psi_1} & H_{2n+4k}(B\diff_\partial(W_\parpun)) \ar[r]^-{\psi_2} & H_{2n+4k}(B\homeo_\partial(W_{g,1})).
}
\]
It is understood that $n$ and $k$ satisfy the conditions of theorem~\ref{thm-gmain}, and we may assume $k>0$.

(i) Using \cite[Thm.1.1]{GalRW2014} we can find $x\in H_{2n+4k}(B\diff_\partial(W_{g,1}))$
such that $\langle\kappa_\sP,\psi_2\psi_1(x)\rangle$ is nonzero in $\QQ$ for selected $\sP\in H^{4n+4k}(B\TOP)$.

(The number $\langle\kappa_\sP,\psi_2\psi_1(x)\rangle$ depends only on the image of
$\sP$ in $H^{4n+4k}(B\SOr(2n))$. In particular it has to be zero if $\sP$
is a monomial $p_{a_1}p_{a_2}\cdots p_{a_j}$ in the Pontryagin classes where at least
one of $a_1,a_2,\dots,a_j$ is greater than $n$.)

(ii) Using mainly section~\ref{sec-calculus}, we are going to produce linear endomorphisms
$\varphi$ of $H_{2n+4k}(B\diff_\partial(W_\parpun))$. These endomorphisms are induced by self-maps
of a certain approximation to $B\diff_\partial(W_\parpun)$. The approximation comes
from proposition~\ref{prop-unravelBoaWe}; it is important that it is sharp in rational homology $H_{2n+4k}$.

(iii) Using sections~\ref{sec-homhirz} and~\ref{sec-conhoto}, we will be able to express $\langle\kappa_\sP,\psi_2\varphi\psi_1(x)\rangle$
in terms of $\langle\kappa_\sP,\psi_2\psi_1(x)\rangle$. With careful choices of $\varphi$ and $x$
we can achieve that
$\langle\kappa_\sP,\psi_2\varphi\psi_1(x)\rangle$ is nonzero for the indecomposable
monomial, $\sP=p_{n+k}$, and zero for all decomposable monomials in the Pontryagin classes, i.e.,
those with more than one
factor. For some nonzero integer $t$, we can represent $\psi_2\varphi\psi_1(tx)$
by a map $f$ from a stably framed (smooth) closed manifold $M$ of dimension $2n+4k$ to $B\homeo_\partial(W_{g,1})$.
Such an $f$ corresponds to a bundle with fiber $W_{g,1}$ on $M$ (so that the boundary bundle with fiber
$\partial W_{g,1}\cong S^{2n-1}$ is trivialized). Glue in disks $D^{2n}$ fiberwise to make a bundle $E\to M$ with
fiber $W_g$\,. This can be taken as $E\to M$ in theorem~\ref{thm-gmain}.

\subsection{Mapping spaces and rationalization} Mapping spaces do not always
behave in a very predictable way under rationalization of the arguments. An easy example: let $X=S^n$ and
let $Y$ be an Eilenberg-MacLane space $K(\ZZ,n)$ (for the same $n\ge 1$). Then the rationalizations $X_\QQ$
and $Y_\QQ$ are defined. We have $\map_*(X,Y)\simeq \ZZ$, but $\map_*(X_\QQ,Y_\QQ)\simeq \map_*(X,Y_\QQ)\simeq\QQ$.
It can not be claimed that $\map_*(X,Y_\QQ)$ is the rationalization of $\map_*(X,Y)$.
For another example, let $X=S^1\wedge B\ZZ/2$ and let $Y=\Omega^\infty\Sigma^\infty S^3$. Then
$\pi_2(\map_*(X,Y_\QQ))= 0$ but $\pi_2(\map_*(X,Y)))\otimes\QQ$ is uncountable (by the affirmed
Segal conjecture for the group $\ZZ/2$).

\smallskip
\emph{Terminology} (for this section). A map of spaces $X\to Y$ is an \emph{$\rhh$-equivalence}
(rhh for \emph{rational higher homotopy}) if, for every $x\in X$, the induced homomorphism
$\pi_j(X,x)\to \pi_j(Y,f(x))$ is a rational isomorphism for $j\ge 3$ and a rational injection for $j=2$.
Equivalently, $X\to Y$ is an $\rhh$-equivalence if, for every $y\in Y$, the homotopy fiber
$\hofiber_y[\,X\to Y\,]$
is either empty or has rationally trivial homotopy groups $\pi_j$ for every $j\ge 2$ and every
choice of base point (in that homotopy fiber).
A space $X$ is \emph{$\rhh$-local} if $\pi_j(X,x)$ is a rational vector space for every $j\ge 2$
and every $x\in X$. 
(This choice of words has drawbacks. If $f\co X\to Y$ is an $\rhh$-equivalence and $Z$ is $\rhh$-local, the map
$f^*\co [Y,Z]\to [X,Z]$ of homotopy sets can easily fail to be bijective.)

\begin{lem} \label{lem-makestart} Suppose that each of the vertical arrows in a commutative diagram of spaces
\[
\xymatrix@R=12pt{
X \ar[r] \ar[d] & Y \ar[d] & \ar[l] Z \ar[d] \\
X' \ar[r] & Y' & \ar[l] Z'
}
\]
is an $\rhh$-equivalence.
Then the induced map from the homotopy limit of the first row
to the homotopy limit of the second row is an $\rhh$-equivalence. If moreover
$X'$, $Y'$ and $Z'$ are $\rhh$-local, then $\holim(X'\to Y'\leftarrow Z')$ is also $\rhh$-local.  \qed
\end{lem}

\begin{lem} Let $X$ be a compact CW-space. Let $f\co Y_0\to Y_1$ be a map of spaces which is an $\rhh$-equivalence.
Then the induced map
\[   \map(X,Y_0)\lra \map(X,Y_1) \]
is an $\rhh$-equivalence. Moreover, if $Y_1$ is $\rhh$-local, then $\map(X,Y_1)$ is $\rhh$-local.
\end{lem}

\proof Induction over the skeletons of $X$, using lemma~\ref{lem-makestart}. \qed

\begin{lem} \label{lem-muchmalignedrat}
Let $\Gamma$ be a discrete group. Let $X$ be a $\Gamma$-CW-space with only free $\Gamma$-cells \cite[\S II.1]{tomDieckTrafo} Suppose
that the CW quotient $X/\Gamma$ is compact. Let $f\co Y_0\to Y_1$ be a $\Gamma$-map of spaces with action of $\Gamma$.
Suppose that $f$ is an $\rhh$-equivalence. Then the induced map
$\map^\Gamma(X,Y_0) \to \map^\Gamma(X,Y_1)$
is an $\rhh$-equivalence. Moreover, if $Y_1$ is $\rhh$-local, then $\map^\Gamma(X,Y_1)$ is also $\rhh$-local.
\end{lem}

\proof Induction over the skeletons of $X$, using lemma~\ref{lem-makestart}. \qed

\medskip
We come to similar statements for spaces of derived maps between topolo\-gical operads.
Let $P$, $Q$ and $Q'$ be operads in spaces satisfying the conditions of theorem~\ref{thm-Goeppl}. We impose
one more condition on $P$: for all $r\ge 2$, the homotopy orbit spaces
\[  P(r)_{h\Sigma_r} \qquad\textup{ and }\qquad \left(\bound_r(P)\right)_{h\Sigma_r} \]
are weakly equivalent to compact CW-spaces. (It follows that there exist two
$\Sigma_r$-CW-spaces with finitely many free $\Sigma_r$-cells and maps from these
to $P(r)$ and $\bound_r(P)$, respectively, which are ordinary weak homotopy equivalences.)
Let $f\co Q\to Q'$ be a map of operads.

\begin{prop} \label{prop-ratoper} Under these conditions, if $f_s\co Q(s)\to Q'(s)$ is an $\rhh$-equivalence
for all $s\ge 0$, then for finite $r\ge 0$ the induced map
\[   \rmap(P_{\le r},Q_{\le r}) \lra \rmap(P_{\le r},Q'_{\le r})  \]
is an $\rhh$-equivalence. And if $Q'(s)$ is $\rhh$-local for all $s\ge 0$, then
$\rmap(P_{\le r},Q'_{\le r})$ is also $\rhh$-local.
\end{prop}

\proof We concentrate on the first statement. (The proof of the second statement follows similar lines, but it is
easier.) The main idea is to use theorem~\ref{thm-Goeppl} and induction on $r$. In addition
we will use lemma~\ref{lem-makestart}, lemma~\ref{lem-muchmalignedrat} and two more observations.
\begin{itemize}
\item[(a)] A composition of two $\rhh$-equivalences is an $\rhh$-equivalence.
\item[(b)] Let $\sA$ be a finite poset, let $F_1$ and $F_2$ be functors from $\sA$
to spaces, and let $\nu\co F_1\to F_2$ be a natural transformation. If $\nu_a\co F_1(a)\to F_2(a)$
is an $\rhh$-equivalence for every $a\in \sA$, then the map from $\holim~F_1$ to $\holim~F_2$
induced by $\nu$ is an $\rhh$-equivalence. Moreover, if each value of $F_2$ is $\rhh$-local, then
$\holim~F_2$ is $\rhh$-local. (This is slightly stronger than lemma~\ref{lem-makestart}, but it can also
be deduced from lemma~\ref{lem-makestart} by a finite iteration.)
\end{itemize}
It follows from lemma~\ref{lem-cobound}
and observation (b) that the map from $\cobound_r(Q)$ to $\cobound_r(Q')$
induced by $f$ is an $\rhh$-equivalence.
Knowing that much, we can take a look at the commutative (up to preferred homotopy) square
\[
\xymatrix@R=14pt{
{\rmap^{\Sigma_r}(P(r),Q(r))} \ar[d] \ar[r] & {\rmap^{\Sigma_r}(P(r),Q'(r))} \ar[d] \\
{\rmap^{\Sigma_r}\left(\!\!\begin{array}{l} \bound_r(P)\to P(r), \\ \rule{1mm}{0mm} Q(r)\to \cobound_r(Q) \end{array}\!\!\right)}
\ar[r] & {\rmap^{\Sigma_r}\left(\!\!\begin{array}{l} \bound_r(P)\to P(r), \\ \rule{1mm}{0mm}
Q'(r)\to \cobound_r(Q')\end{array}\!\!\right)}
}
\]
By lemma~\ref{lem-muchmalignedrat}, the upper horizontal arrow one is an $\rhh$-equivalence. By the same lemma and observation (b),
the lower horizontal arrow is also an $\rhh$-equivalence. (Observation (b) is useful here because the two mapping spaces
in the lower row of the square can be written as homotopy limits of diagrams, indexed by finite posets,
of more ordinary mapping spaces as in lemma~\ref{lem-muchmalignedrat}.)
It follows that, for every choice of base point in the lower left-hand term of the square,
the map of vertical homotopy fibers selected by that choice is an $\rhh$-equivalence. Therefore, by theorem~\ref{thm-Goeppl}
and induction on $r$,
\[ \rmap(P_{\le r},Q_{\le r}) \lra \rmap(P_{\le r},Q'_{\le r}) \]
is an $\rhh$-equivalence. \qed

\smallskip
By combining proposition~\ref{prop-ratoper} and theorem~\ref{thm-operpest}
we obtain lemma~\ref{lem-unravelBoaWe} below, which combines well with proposition~\ref{prop-unravelBoaWe}.
(Thanks to M Krannich and A Kupers for help with the formulation.)

\emph{Notation.} For a space $Y$ which has the homotopy type of a CW-space and is path-componentwise simple \cite[\S1.1.2]{MoritaGCC},
let $Y_\QQ$ be the path-componentwise rationalization. Let
\begin{eqnarray*}  J_0 &:= & \pi_0(\map_*(W_\parpun,\haut_\fin(\config(\RR^{2n};r)))),  \\
J &:= & \pi_0(\map_*(W_\parpun,\haut_\fin(\config(\RR^{2n};r)_\QQ))), \\
I_0 &:= & \pi_0(\map_*(W_\parpun,\SOr(2n))), \\
I & := & \pi_0(\map_*(W_\parpun,\SOr(2n)_\QQ)).
\end{eqnarray*}

\begin{lem} \label{lem-unravelBoaWe}
For any $r\ge 1$ and $n\ge 3$,
the inclusion of the homotopy pullback of \eqref{eqn-mainbdryspec3} in the homotopy pullback of
\begin{equation} \label{eqn-mainbdryspec5}
\begin{aligned}
\xymatrix@C=-100pt@R=40pt{
{B\left(\!\!\!\begin{aligned} \begin{array}{l}  \haut_*(W_\parpun)\,\,\, \ltimes \\
I_0\times_I\map_*(W_\parpun,\SOr(2n)_\QQ)
\end{array}\end{aligned}\!\!\right)} \ar@/^1pc/[dr] &&
\ar@/_1pc/[dl]^-{\textup{specialization}} B\haut^\partial_\finplus(\config(W_\parpun;r)) \\
 & {B\left(\!\!\!\begin{aligned}
\begin{array}{l}  \haut_*(W_\parpun) \,\,\, \ltimes \\
J_0\times_J\map_*(W_\parpun,\haut_\fin(\config(\RR^{2n};r)_\QQ))
\end{array} \end{aligned} \!\!\right)}
}
\end{aligned}
\end{equation}
induces an isomorphism in the rational homology of the base point components.
\end{lem}

\proof[First part of the proof] We show that \emph{the inclusion of the homotopy pullback of \eqref{eqn-mainbdryspec3} in the homotopy pullback of}
\begin{equation} \label{eqn-mainbdryspec4}
\begin{aligned}
\xymatrix@C=-90pt@R=40pt{
{B\left(\!\!\!\begin{aligned} \begin{array}{l}  \haut_*(W_\parpun)\,\,\, \ltimes \\
\big(\map_*(W_\parpun,\SOr(2n))\big)_\QQ
\end{array}\end{aligned}\!\!\right)} \ar@/^1pc/[dr] &&
\ar@/_1pc/[dl]^-{\textup{specialization}} B\haut^\partial_\finplus(\config(W_\parpun;r)) \\
 & {B\left(\!\!\!\begin{aligned}
\begin{array}{l}  \haut_*(W_\parpun) \,\,\, \ltimes \\
\big(\map_*(W_\parpun,\haut_\fin(\config(\RR^{2n};r)))\big)_\QQ
\end{array} \end{aligned} \!\!\right)}
}
\end{aligned}
\end{equation}
\emph{induces an isomorphism in rational homology.} Proof of this claim: write
$X$ for $\map_*(W_\parpun,\SOr(2n))$ and $Y$ for $\map_*(W_\parpun,\haut_\fin(\config(\RR^{2n};r)))$.
Both are group-like $A_\infty$ spaces using pointwise multiplication of maps
and the standard $A_{\infty}$ multiplications of $\SOr(2n)$ and $\haut_\fin(\config(\RR^{2n};r))$.
Both the homotopy pullback of \eqref{eqn-mainbdryspec3} and the homotopy pullback
of \eqref{eqn-mainbdryspec4} come with forgetful maps to
\[ B\haut^\partial_\finplus(\config(W_\parpun;r)). \]
The homotopy fiber (over the base point) of the first of these is
$\hofiber[\,BX\to BY\,]$, abbreviated to $Y/\!\!/X$. The homotopy fiber of the second map
is $Y_\QQ/\!\!/X_\QQ$. By a spectral
sequence argument or a filtration argument, it suffices therefore to show that the inclusion(-induced)
map
\begin{equation} \label{eqn-uglycomp}
 Y/\!\!/X \lra Y_\QQ/\!\!/X_\QQ
 \end{equation}
induces an isomorphism in rational homology. Source and target of~\eqref{eqn-uglycomp} are based spaces and
componentwise simple spaces. [Reason: $X$ and $Y$ admit compatible $(n+1)$-fold deloopings due to the fact that
$W_\parpun$ is homotopy equivalent to a wedge of $n$-spheres.]
Noting that $Y$ acts on both, we can reduce further to the statement that~\eqref{eqn-uglycomp} induces a bijection in $\pi_0$ and
isomorphisms in $\pi_j\otimes\QQ$ for all $j\ge 1$. For $\pi_0$ this is true by construction. For $\pi_j\otimes\QQ$,
where $j\ge 1$, it is true by construction and the five lemma.

\proof[Second part of the proof] Write $P$ for $\sE_{2n}$~, the operad of little $2n$-disks. There is a commutative diagram
\[
\xymatrix@C=15pt@M=5pt@R=12pt{
\haut((P_{\le r})_\QQ) \ar[d] \ar[r] & \ar[d] \big(\haut((P_{\le r})_\QQ)\big)_\QQ &
\ar[d] \ar[l] (\haut(P_{\le r}))_\QQ  \\
\haut_\fin(\config(\RR^{2n};r)_\QQ) \ar[r] & \big(\haut_\fin(\config(\RR^{2n};r)_\QQ)\big)_\QQ
& \ar[l] \big(\haut_\fin(\config(\RR^{2n};r))\big)_\QQ
}
\]
The forgetful vertical arrows are weak equivalences by theorem~\ref{thm-operpest}.
The arrows in the upper row induce isomorphisms in $\pi_j$ for $j\ge 3$ by proposition~\ref{prop-ratoper},
if we choose the standard base points. Therefore the arrows in the lower row
induce isomorphisms in $\pi_j$ for $j\ge 3$,
if we choose the standard base points. Now since $W_\parpun$ has the based homotopy type of a wedge
of $n$-spheres where $n\ge 3$, it follows that the resulting maps
\[
\xymatrix@R=10pt{
\map_*(W_\parpun,\haut_\fin(\config(\RR^{2n};r)_\QQ)) \ar[d] \\
\map_*\big(W_\parpun,\big(\haut_\fin(\config(\RR^{2n};r)_\QQ)\big)_\QQ\big) \\
\map_*\big(W_\parpun,\big(\haut_\fin(\config(\RR^{2n};r))\big)_\QQ\big)  \ar[u]
}
\]
are weak homotopy equivalences. Therefore we can write
\begin{eqnarray*}
&& J_0\times_J \map_*(W_\parpun,\haut_\fin(\config(\RR^{2n};r)_\QQ)) \\
&\simeq& J_0\times_J \map_*(W_\parpun,(\haut_\fin(\config(\RR^{2n};r)))_\QQ)  \\
&\simeq & (\map_*(W_\parpun,\haut_\fin(\config(\RR^{2n};r))))_\QQ\,.
\end{eqnarray*}
Similarly but more easily, we can write
\[
I_0\times_I \map_*(W_\parpun,\SOr(2n)_\QQ) \simeq  (\map_*(W_\parpun,\SOr(2n)))_\QQ\,.
\]
Therefore the statement which we are after is equivalent to what we have established in the first part of the proof. \qed

\subsection{Stabilization} \label{subsec-stabi}
For integers $m_1$ and $m_2$ (where $m_2\ge m_1+3$ and $m_1\ge 3$) there is a diagram
\begin{equation} \label{eqn-disso}
\begin{aligned}
\xymatrix@C=12pt@R=14pt{
B\SOr(m_1)  \ar[r] \ar[d] & \ar[d] B\haut_{\fin}(\config(\RR^{m_1};r)) \ar[r] & B\haut_{\fin}(\config(\RR^{m_1};r)_\QQ) \ar[d]\\
B\SOr(m_2) \ar[r] &  B\haut_{\fin}(\config(\RR^{m_2};r))  \ar[r] &  B\haut_{\fin}(\config(\RR^{m_2};r)_\QQ)
}
\end{aligned}
\end{equation}
commutative up to preferred homotopies, where the vertical arrows are stabilization maps.
(The left-hand square is already defined if $m_2\ge m_1\ge 0$. For the right-hand square
we need $m_1\ge 3$ and $m_2-m_1\ge 3$; this proviso is needed only because we are using the
\emph{rationalized} configuration categories. The horizontal arrows on the left appeared earlier, e.g.~in diagram~\eqref{eqn-mainbdryspec2coord}. A linear
isometry $\RR^m\to \RR^m$ can be viewed as a homeomorphism from the manifold $\RR^m$
to itself and as such determines automorphisms of $\config(\RR^m)$ and of $\config(\RR^m;r)$, for any $r\ge 1$.)
The stabilization maps in the middle and right-hand columns are not a trivial matter; they use something inspired
by the Boardman-Vogt tensor product of operads. More details are given in \cite{BoavidaWeiss_conftensor} and in
remark~\ref{rem-lackof} below.

We will see in a moment that
the homotopy groups of the upper right-hand term in diagram~\eqref{eqn-disso}, respectively lower right-hand term,
are clustered in dimensions close to the integer multiples of $m_1-2$, respectively $m_2-2$.
It follows that the right-hand vertical arrow in diagram~\eqref{eqn-disso} is far from being a good
homotopical approximation in many cases where $m_2$ is large and $m_1$ is close to $m_2$ but not equal to $m_2$\,.
This gave rise
to the section title, \emph{dissonance}. It is a case where homological stability fails to hold.
That makes a surprising contrast with the left-hand vertical arrow.
It is precisely that contrast which lets us deduce that the map across the lower row in diagram~\eqref{eqn-disso}
is far from being a good homotopical approximation (under similar assumptions on $m_2$ and $m_1$)
in the rational setting.

\begin{lem} \label{lem-boundestim} If $m\ge 3$ and $r\ge 5$, then $\bound_r(\sE_m)$ is $1$-connected.
\end{lem}
\proof Since $m\ge 3$, the spaces $\sE_m(s)$ for $s\ge 0$ are all $1$-connected.
Consequently the spaces $(N_d\sE_m)_S$ for $S$ in $\tree^c$ are all $1$-connected.
It follows that the projection map
\[  \bound_r(\sE_m)\,\,\, = \!\!\!\!\!\hocolimsub{\twosub{t_r\to S}{S\textup{ in }\tree^c_{\le r-1}}}\!\!\!\!\!\! (N_d\sE_m)_S
\quad\lra \quad  B((t_r\downarrow F_{r-1})^\op) \]
(notation as in definition~\ref{defn-boundcobound}) is $2$-connected. On the other hand,
\[  B((t_r\downarrow F_{r-1})^\op)\times \Sigma_r  \simeq \bound_r(\sE_1)\,. \]
It follows from proposition \ref{prop-bound} that $\bound_r(\sE_1)$ is weakly equivalent, as a space with action
of $\Sigma_r$\,, to a product $S^{r-3}\times\Sigma_r$. (The
Axelrod-Singer-Fulton-Macpherson compactification of the normalized ordered configuration space of $r$ points in $\RR^1$
can be identified with a Stasheff polytope \cite{StashH}, times $\Sigma_r$.) Therefore $B((t_r\downarrow F_{r-1})^\op)$ is weakly equivalent to $S^{r-3}$. \qed

\begin{lem} \label{lem-appligoeppl} In theorem~\ref{thm-Goeppl}, suppose that $P$ is $\sE_m$ and $Q$ is the rationalization of $\sE_m$\,.
Suppose also that $m\ge 3$ and $r\ge 5$. The rationally nontrivial homotopy groups of the homotopy fiber (over the base point) of the forgetful map
\[  \rmap(P_{\le r},Q_{\le r}) \lra  \rmap(P_{\le r-1},Q_{\le r-1}) \]
are in dimensions which belong to intervals of the form
\[  J_{m,r,s}:=[s(m-2)-3r+6, s(m-2)+3]\cap [0,\infty)  \]
where $s=0,1,2,3,\dots$. Consequently the nonzero rational homotopy groups of
$\rmap(\config(\RR^m;r),\config(\RR^m;r))$ \emph{(standard choice of base point)}
appear only in dimensions which belong to one of the intervals $J_{m,r,s}$.
\end{lem}
[Later we will specialize by taking either $r=8$ or $r=6$. For fixed $m>23$ the
intervals $J_{m,8,s}$ are pairwise disjoint, of length 21 unless $s=0$, and with spacing of length $m-23=m-2-21$.
Similarly for fixed $m>15$, the intervals $J_{m,6,s}$ are pairwise disjoint, of length 15 unless $s=0$, and with spacing
of length $m-17=m-2-15$.]

\proof
Relying on proposition~\ref{prop-bound}, we may abuse notation slightly to write
$\partial P(r)$ for $\bound_r(P)$ and pretend that $(P(r),\partial P(r))$ is an oriented
Poincar\'e pair of formal dimension $(r-1)m-1$. The integral cohomology of $P(r)$ is concentrated in dimensions which
are multiples of $m-1$. See \cite[\S III.6]{CohenLadaMay}.
The relative homology of $(P(r),\partial P(r))$ is therefore concentrated in dimensions
of the form $\ell(m-1)+r-2$ where $\ell\in \{0,1,\dots,r-1\}$. All these homology groups are free
abelian and the spaces involved are simply connected. (Here we need lemma~\ref{lem-boundestim} and consequently we need
$r\ge 5$.)
It follows that, up to homotopy equivalences, there is a cell decomposition of $P(r)$ relative
to $\partial P(r)$ where all relative cells are in dimensions of the form $\ell(m-1)+r-2$.
The rationally nontrivial homotopy groups $\pi_t$ of $Q(r)$ are all in dimensions of the form
$t=j(m-2)+1$ where $j=1,2,3,\dots$. (To see this, use the Fadell-Neuwirth fibrations $P(r)\to P(r-1)$,
again \cite[\S III.6]{CohenLadaMay}, Hilton's theorem on the
homotopy groups of wedges of spheres, and induction on $r$.)
The rationally nontrivial homotopy groups of the homotopy fiber of $Q(r)\to \cobound_r~Q$ are therefore in dimensions
clustered in intervals of the form
\[ [j(m-2)-r+3,\,j(m-2)+1]~; \]
here we use lemma~\ref{lem-cobound}.
Therefore by theorem~\ref{thm-Goeppl},
the rationally nontrivial homotopy groups of the homotopy fiber (over the base point) of the forgetful map
\[  \rmap(P_{\le r},Q_{\le r}) \lra  \rmap(P_{\le r-1},Q_{\le r-1}) \]
are in dimensions which belong to intervals of the form
\begin{eqnarray*}
&&  [j(m-2)-r+3,j(m-2)+1] - \ell(m-1)-r+2  \\
& = & [j(m-2)-2r+5-\ell(m-1), j(m-2)-r+3-\ell(m-1)]  \\
& = & [(j-\ell)(m-2)-2r+5-\ell, (j-\ell)(m-2)-r+3-\ell] \\
& \subset & [(j-\ell)(m-2)-3r+6, (j-\ell)(m-2)+3]
\end{eqnarray*}
where $j=1,2,\dots$ but $\ell=0,1,\dots,r-1$ only. Writing $s$ for $j-\ell$ we can restate this by saying that the
nonzero rational homotopy groups of the said homotopy fiber
are in dimensions which belong to intervals of the form
\[  J_{m,r,s}:=[s(m-2)-3r+6, s(m-2)+3]\cap [0,\infty)  \]
where $s=0,1,2,3,\dots$.

Note that $J_{m,r-1,s}\subset J_{m,r,s}$ if both are defined.
It follows that the nonzero rational homotopy groups of
$\rmap(P_{\le r},Q_{\le r})$ (with the standard choice of base point)
appear only in dimensions which belong to one of the intervals
$J_{m,r,s}$, where $s=0,1,2,\dots$. Equivalently, the nonzero rational homotopy groups of
$\rmap(Q_{\le r},Q_{\le r})$ (standard choice of base point: the identity map)
appear only in dimensions which belong to one of the intervals
$J_{m,r,s}$. Equivalently (where we use theorem~\ref{thm-operpest}), the nonzero rational homotopy groups of
$\rmap(\config(\RR^m;r),\config(\RR^m;r))$ (with the standard choice of base point)
appear only in dimensions which belong to one of the intervals
$J_{m,r,s}$.  \qed

\medskip
This calculation feeds into yet another technical lemma. We take $m_2=2n$ and $r=8$ in diagram~\eqref{eqn-disso}.
The space $B\SOr(2n)_\QQ$ is a product of Eilenberg-MacLane spaces. The
Pontryagin classes $p_1,\dots p_{n-1}$ in degrees $4,8,\dots,4n-4$ together with the Euler class $e$
in degree $2n$ define such a splitting. It follows that for $b\in\{1,2,\dots,n-1\}$ there is a unique homotopy class of
based maps
\[ q_b\co B\SOr(2n)_\QQ\to B\SOr(2n)_\QQ \]
such that $q_b^*(p_j)=p_j$ for $j\ne b$ and $q_b^*(e)=e$, whereas $q_b^*(p_b)=0$. We ask how this interacts with
the standard comparison map
\[  w\co B\SOr(2n)_\QQ \lra B\haut_{\fin}(\config(\RR^{2n};8)_\QQ). \]
\emph{Terminology}. For a space $Y$ and an integer $t\ge 0$, let $\tau_{\le t}Y$ be the Postnikov truncation
(obtained by killing the homotopy groups of $Y$ in dimensions greater than $t$). We may also write $\tau_{<t+1}Y$
instead of $\tau_{\le t}Y$. For a based space $X$ and $t>0$, let $\tau_{\ge t}X$ be the homotopy
fiber (over the base point) of the standard map $X\to \tau_{<t}X$. We may also write $\tau_{>t-1}X$ for that.
A map from $X$ to $Y$ \emph{defined in the homotopy range $[s,t]$} is a map from $\tau_{\ge s}X$
to $\tau_{\le t}Y$. An ordinary map $X\to Y$ can always be viewed as a map defined in the homotopy range $[s,t]$.

\begin{lem} \label{lem-estimate} If $n\ge 83$
and $b=n-11$, the map $q_b$ as a map defined in the homotopy range $[n+1,8n-27]$ satisfies $wq_b\simeq w$.
If $n\ge 59$ and $b=\lfloor 3n/4\rfloor$, the map
$q_b$ as a map defined in the homotopy range $[n+1,6n-35]$ satisfies $wq_b\simeq w$. In both cases it can be
arranged that $q_b$ is the identity on the $(n+1)$-skeleton of $\tau_{>n}B\SOr(2n)_\QQ$ and the homotopy from $wq_b$ to $w$
is stationary on the $(n+1)$-skeleton.
\end{lem}

\proof For the first case, $n\ge 83$ and $b=n-11$:
write $X$ for $\tau_{>n}B\SOr(2n)_\QQ$ and write $Y$ for
$\tau_{\le 8n-27}B\haut_{\fin}(\config(\RR^{2n};8)_\QQ)$. We take the freedom to write
$w\co X\to Y$ and $q_b\co X\to X$. The goal is to show $w\simeq wq_b$ in that setting.
Write $X= X_\lambda\times X_\mu \times X_\rho \times X_\eta$   
where
\begin{itemize}
\item $X_\lambda$ is the product of the factors corresponding to $p_j$ where $j<b$;
\item $X_\mu$ is the factor corresponding to $p_b$;
\item $X_\rho$ is the product of the factors corresponding to $p_j$ where $j>b$;
\item $X_\eta$ is the factor corresponding to the Euler class $e$.
\end{itemize}
The maps $w$ and $wq_b$ agree up to homotopy on $X_\lambda\times X_\mu$ because, as we have just
shown, $w$ is nullhomotopic there and consequently $wq_b$ is also nullhomotopic there.
[Explanation: pretend that $X_\lambda\times X_\mu\subset B\SOr(2b)_\QQ$ and remember diagram~\eqref{eqn-disso},
taking $m_2=2n$ and $m_1=2b$. The idea is to exploit dissonance and to use
obstruction theory to conclude that the right-hand vertical arrow in that diagram is nullhomotopic
in the homotopy range $[n,8n-27]$. So we need to have an idea where, i.e. in which degrees,
the nontrivial rational cohomology
of $B\haut_{\fin}(\config(\RR^{2b};8)_\QQ)$ is located and where the nontrivial
rational homotopy of $B\haut_{\fin}(\config(\RR^{2n};8)_\QQ)$ is located. From the above calculation
we have such an idea and consequently it is enough to verify that a linear combination
$\sum a_ix_i$  with positive integer coefficients $a_i$ and
\[ x_i\in  J_{2b,8,0}\cup J_{2b,8,1} \cup J_{2b,8,2} \cup J_{2b,8,3} \cup J_{2b,8,4} \]
will never belong to one of the intervals $J_{2n,8,s}$
where $s=1,2,3$. The verification is mechanical.
Note that $J_{2n,8,4}$ is not of interest because its minimal element is $8n-26$,
and $J_{2b,8,5}$ is not of interest because its minimal element is $10n-138$, which is
even greater. The condition $n\ge 83$ ensures not only that $10n-138>8n-27$, but also
that linear combinations $a_1x_1+a_2x_2+a_3x_3+a_4x_4$ where $x_i\in J_{2b,8,1}$
and $a_i=1$ for $i=1,2,3,4$ do not belong to $J_{2n,8,3}$.] The maps $w$ and $wq_b$ also
agree on $X_\lambda\times X_\rho\times X_\eta$ because $q_b$ is the identity there.
So they agree (up to a homotopy which we choose) on
\[ X_\lambda\times\big(X_\mu\vee X_\rho\times X_\eta\big). \]
The remaining obstructions to agreement of $w$ and $wq_b$ on all of $X$
(as in obstruction theory) are in the reduced cohomology $\tilde H^j$ of
\[  \frac{X_\lambda\times\big(X_\mu\times X_\rho\times X_\eta\big)}{X_\lambda\times\big(X_\mu\vee X_\rho\times X_\eta\big)}
= (X_\lambda)_+\wedge X_\mu\wedge (X_\rho\times X_\eta) \]
with coefficients in the homotopy groups $\pi_j(Y)$. By a straightforward computation,
there is no $j$ such that both the reduced cohomology of $(X_\lambda)_+\wedge X_\mu\wedge(X_\rho\times X_\eta)$ in degree $j$ and the
homotopy group $\pi_j(Y)$ are nontrivial.
(The first nontrivial reduced cohomology group of $(X_\lambda)_+\wedge X_\mu\wedge(X_\rho\times X_\eta)$
is in dimension $j=4b+2n=6n-44$, but $\pi_j(Y)$ is zero for this $j$. The next run of nontrivial cohomology
groups appears in dimensions $j$ greater than $7n$ approximately, but this is not matched by nontrivial $\pi_j(Y)$
until we reach $j=4(2n-2)-18=8n-26$ and remember that $\pi_j(Y)$ for $j\ge 8n-26$ is zero by the definition of $Y$.)

For the second case:
$Y$ is now $\tau_{\le 6n-35}B\haut_{\fin}(\config(\RR^{2n};8)_\QQ)$ and $b$ is now $\lfloor 4n/3\rfloor$.
We assume $n\ge 59$.
In the first step, to show that $w$ and $wq_b$ agree up to homotopy on
$X_\lambda\times X_\mu$ we need to verify that a linear combination
$\sum a_ix_i$  with positive integer coefficients $a_i$ and
\[ x_i\in  J_{2b,8,0}\cup J_{2b,8,1} \cup J_{2b,8,2} \cup J_{2b,8,3} \]
will never belong to one of the intervals $J_{2n,8,s}$
where $s=1,2$. The verification is mechanical. Here $J_{2b,8,4}$ is not of interest because its minimal element is greater
than $6n-35$, and $J_{2n,8,3}$ is not of interest because its minimal element is even greater.
The condition $n\ge 59$ ensures
that linear combinations $a_1x_1+a_2x_2+a_3x_3$ where $x_i\in J_{2b,8,1}$
and $a_i=1$ for $i=1,2,3$ do not belong to $J_{2n,8,2}$. For the second step we observe that
$(X_\lambda)_+\wedge X_\mu\wedge(X_\rho\times X_\eta)$ is $(4b+2n-1)$-connected. Therefore no further obstructions
arise in the obstruction theory scheme since the maximal element of $J_{2n,8,2}$ is $4n-1$, less than $4b+2n$.

For the third part of the lemma: by construction these conditions will be satisfied if the inclusion
of $X_\lambda\times X_\mu$ in $X$ is $(n+1)$-connected. We need $4(n-11)>n+1$, respectively  $4\lfloor 3n/4\rfloor>n+1$.
These lower bounds are of course much lower than the ones above, $n\ge 83$ and $n\ge 59$, respectively.
\qed

\begin{rem} \label{rem-randomconst} The formulation of lemma~\ref{lem-estimate} leaves a few questions.  \newline
\emph{What were the reasons for taking $r=8$~?}
The number 8 makes another appearance in the lemma, in the choice of Postnikov truncation
$\tau_{\le 8n-27}$. This is more fundamental than the choice $r=8$.
The proof did not come into being with a random choice $r=8$. It is closer to the truth
to say that the proof began with $X=\tau_{>n}B\SOr(2n)_\QQ$ and $Y=\tau_{\le f(n)}B\haut_{\fin}(\config(\RR^{2n};r)_\QQ)$
where $r$ and the function $f$ were initially unspecified. Then came the product decomposition
\[ X=X_\lambda\times X_\mu \times X_\rho \times X_\eta \]
and obvious connectivity statements for the factors. This (together with the calculation preceding the lemma) suggested
that $f(n)=8n-s$ for some constant $s$ was a good choice. Then, in the light of section~\ref{sec-conhoto}, it became clear
that $r=8$ was convenient. We could also have taken a larger $r$,
but that would only have made the estimates slightly worse. (The intervals $J_{m,r,s}$ get bigger
as $r$ grows.) --- In the second part of the lemma, where we used the Postnikov truncation $\tau_{\le 6n-35}$,
we could have taken $r=6$ and this could have improved the estimates slightly. \newline
\emph{Why was it necessary to distinguish two cases, $b=n-11$ and $b=\lfloor 3n/4\rfloor$~?}
We will see that the first case, $b=n-11$, helps us to establish cases of theorem~\ref{thm-pmain} where $k$ is
fairly large, and the second case, $b=\lfloor 3n/4\rfloor$, allows us to establish cases
of theorem~\ref{thm-pmain} where $k$ is fairly small in relation to $n$.
\end{rem}

\begin{rem} \label{rem-lackof}
\emph{Rationalized} configuration categories are not explicitly
mentioned in \cite{BoavidaWeiss_conftensor}. Therefore it may be appropriate to review the
construction of the middle column in diagram~\eqref{eqn-disso} and point out that
the right-hand column can be set up by the same mechanism. There is an elementary construction
which to $X$ and $Y$, fiberwise complete Segal spaces over $N\fin$, associates a fiberwise
complete Segal space $X\boxtimes^{\textup{pre}} Y$ over $N\fin$. In the case where
$X=\config(M)$ and $Y=\config(N)$, the Segal space $X\boxtimes^{\textup{pre}} Y$
comes with a map to $\config(M\times N)$ which can be characterized by a universal property;
it is a \emph{conservatization} map. In other words, we may write
\[ \config(M\times N)\simeq \Lambda(\config(M)\boxtimes^{\textup{pre}}\config(N)) \]
where $\Lambda$ is an abstract conservatization procedure.
Consequently we can make a canonical map
\[ \rmap_\fin(\config(M),\config(M)) \lra \rmap_\fin(\config(M\times N),\config(M\times N)) \]
by tensoring with the identity of $\config(N)$:
\[
\xymatrix@R=12pt{
\rmap_\fin(\config(M),\config(M)) \ar[d] \\
\rmap_\fin(\config(M)\boxtimes^{\textup{pre}}\config(N),\config(M)\boxtimes^{\textup{pre}}\config(N)) \ar[d] \\
\rmap_\fin(\config(M)\boxtimes^{\textup{pre}}\config(N),\Lambda(\config(M)\boxtimes^{\textup{pre}}\config(N))) \\
\ar[u]_-\simeq
\rmap_\fin(\Lambda(\config(M)\boxtimes^{\textup{pre}}\config(N)),\Lambda(\config(M)\boxtimes^{\textup{pre}}\config(N)))
}
\]
This takes homotopy invertible path components to homotopy invertible path components and so restricts to a map
\[
\xymatrix{
\haut_\fin(\config(M)) \ar[r] &
{\haut_\fin(\Lambda(\config(M)\boxtimes^{\textup{pre}}\config(N))) \,\,\simeq\,\,
\haut_\fin(\config(M\times N))}.
}
\]
Now there are some easy variations. Instead of $X=\config(M)$ and $Y=\config(N)$ we can take the
rationalizations of these. (The rationalizations are meaningful if $X$ and $Y$ are of dimension
$\ge 3$ and simply connected.) Then we obtain in the same way
\[
\xymatrix@C=9pt{
\haut_\fin(\config(M)_\QQ) \ar[r] &
{\haut_\fin(\Lambda(\config(M)_\QQ\boxtimes^{\textup{pre}}\config(N)_\QQ)) \simeq
\haut_\fin(\config(M\times N)_\QQ)}.
}
\]
We can also apply truncation, replacing $\config(M)$ by $\config(M;r)$ etc; in this case we view
$\config(M;r)$ as a fiberwise complete Segal space over the nerve of $\fin_{\le r}$, the full subcategory of $\fin$
spanned by the objects $\uli\ell$ for $\ell \le r$. Then we obtain
\[
\begin{array}{rcl}
\haut_\fin(\config(M;r)) & \lra &
\haut_\fin(\Lambda(\config(M;r)\boxtimes^{\textup{pre}}\config(N;r))) \\
&& \simeq \haut_\fin(\config(M\times N;r)).
\end{array}
\]
Last not least we can also apply truncation and rationalization, taking $\config(M;r)_\QQ$ for $X$
and $\config(N;r)_\QQ$ for $Y$, etc.,
provided $X$ and $Y$ are of dimension $\ge 3$ and simply connected. More specifically we can choose
$M=\RR^{m_1}$ and $N=\RR^{m_2-m_1}$.
\end{rem}

\subsection{The end}
Under the conditions of lemma~\ref{lem-estimate}, the map
\[ q_b\co B\SOr(2n)_\QQ\to B\SOr(2n)_\QQ \]
for $b=2n-11$ admits a homotopy $wq_b\simeq w$ in the homotopy range $[n,8n-27]$. Also, $q_b$ for $b=\lfloor 3n/4\rfloor$ admits a homotopy $wq_b\simeq w$
in the homotopy range $[n,6n-35]$. In both cases, $q_b$ and the homotopy together
determine a self-map $f_b$ of the homotopy pullback of~\eqref{eqn-mainbdryspec5}, in the homotopy ranges
$[0,7n-27]$ and $[0,5n-35]$ respectively. This is induced by $q_b$ on the upper left-hand term
in diagram~\eqref{eqn-mainbdryspec5} and by the identity on the other two terms. The map $f_b$ takes the base
point component to itself (by the last sentence in lemma~\ref{lem-estimate}).

\proof[Proof of theorem~\ref{thm-gmain}.] The constants are or can be $c_1=83$ and $c_2=11$. In particular, we will assume
$n\ge 83$; this comes from (the proof of) lemma~\ref{lem-estimate}. We distinguish three cases.

\emph{Case 1.} Suppose that $n/4<k+11\le n$. Set $b:=n-11$ as in the first part of lemma~\ref{lem-estimate},
and $a:=k+11$, so that $n\ge a>n/4$ and $b>n/4$. By the Galatius-RandalWilliams theorems,
specifically \cite[Thm.1.1]{GalRW2014}, there exists
$x\in H_{2n+4k}(B\diff_\partial(W_{g,1}))$ such that $\langle\kappa(p_ap_b), x\rangle \ne 0$
while $\langle\kappa(\sP),x\rangle=0$ for all other \emph{monomials} $\sP$ of cohomological degree $4n+4k$ in the Pontryagin classes.
By \cite{BeBe}, the coefficient of $p_ap_b$ in the Hirzebruch polynomial $\sL_{n+k}$ is nonzero.
Therefore $\langle\kappa(\sL_{n+k}),x\rangle\ne 0$.
Let $y$ be the image of $x$ in $H_{2n+4k}(B\diff_\partial(W_\parpun))$. By proposition~\ref{prop-unravelBoaWe} and lemma~\ref{lem-unravelBoaWe}
it is allowed to write
\[ (f_b)_*(y) \in H_{2n+4k}(B\diff_\partial(W_\parpun)). \]
(Here we need and we have $2n+4k<7n-15$ because of a bound in example~\ref{prop-unravelBoaWe}.
We need and we have $2n+4k<7n-27$ because of a qualifier in lemma~\ref{lem-estimate}.)
By part (i) of proposition~\ref{prop-uv} and diagrams~\eqref{eqn-sentriesnu1},~\eqref{eqn-sentriesnu2}
(specialized as indicated in diagrams~\eqref{eqn-sentriesappl3},\eqref{eqn-sentriesappl4}),
the element $(f_b)_*(y)$ has a nonzero scalar product with $\kappa_t(\sL_{n+k})$.
(Here we need and we also have $2n+4k<7n-17$.)
By proposition~\ref{prop-uv} part (ii) and proposition~\ref{prop-muchridiculed}, the element $(f_b)_*(y)$ has zero scalar
product with $\kappa_t(\sP)$ for any decomposable polynomial of cohomological degree $4n+4k$
in the Pontryagin classes. After multiplication by a nonzero integer, if necessary, the class $(f_b)_*(y)$ can be represented by a map
\[ \beta\co M  \lra B\diff_\partial(W_\parpun) \]
where $M$ is a smooth closed stably framed manifold (and nullbordant as such). Let $E_M\to M$ be the bundle of
$2n$-dimensional closed oriented topological manifolds with fiber $W_g$ obtained by composing
\[  M\xrightarrow{\beta} B\diff_\partial(W_\parpun) \hookrightarrow B\homeo_\partial(W_{g,1})\to B\homeo^\omega(W_g) \]
where $\homeo^\omega$ is for orientation preserving homeomorphisms.
Now $E_M\to M$ is a bundle with the properties that we require.

\emph{Case 2.} Suppose that $k+11\le n/4$. Then we have $4n+4k<5n$. Therefore we go for
the second part of lemma~\ref{lem-estimate}, taking $b=\lfloor 3n/4\rfloor$. Let
$a:=n+k-b$, so that $a>n/4$. 
Then proceed as in Case 1. (Various upper bounds on $k$ appear;
the most stringent of these is $2n+4k< 5n-35$, which is satisfied.)

\emph{Case 3.} Suppose that $5n/4>k+11>n$. Set $b=n-11$ as in the first part of lemma~\ref{lem-estimate}.
Let $a:=n+k-b=k+11>n$. Now $p_a=0\in H^{4a}(B\SOr(2n))$, so that we cannot proceed exactly as in Case 1.
Choose integers $a_1,a_2\in [n/2,n-1]$ such that $a=a_1+a_2$.
By the Galatius-RandalWilliams theorem there exists
$x\in H_{2n+4k}(B\diff_\partial(W_{g,1}))$ such that $\langle\kappa(p_{a_1}p_{a_2}p_b),x\rangle \ne 0$
while $\langle\kappa(\sP),x\rangle=0$ for all other monomials $\sP$ of cohomological degree $4n+4k$ in the Pontryagin classes.
By \cite{BeBe}, the coefficient of $p_{a_1}p_{a_2}p_b$ in the Hirzebruch polynomial $\sL_{n+k}$ is nonzero.
Therefore $\langle\kappa(\sL_{n+k}),x\rangle\ne 0$. Continue as in Case 1. (Various upper bounds on $k$ appear;
the most stringent of these is $2n+4k< 7n-27$, which is satisfied.)

\subsection{Miscellaneous remarks}

\begin{rem} \label{rem-impest} Two homotopical properties of $W_\parpun$ were used in two different places.
Firstly, $W_\parpun$ is homotopically $n$-dimensional (although it is geometrically $2n$-dimensional); this was important for the
estimates in proposition~\ref{prop-unravelBoaWe}.
But secondly, $W_\parpun$ is $(n-1)$-connected. This allowed us to replace
$\SOr(2n)$ by its $(n-1)$-connected cover in lemma~\ref{lem-estimate}.
\end{rem}

\begin{rem} \label{rem-sigconf} Let $p\co E\to M$ be a fibration where base $M$ and fiber $F$ are oriented
Poincar\'e duality spaces, and $M$ is connected. It is known that the signature of $E$ depends only on $M$ and
$H^n(F;\QQ)$ with the intersection form and the action of $\pi_1(M)$, where $n$ is half the formal dimension of $F$.
(If the formal dimension of $F$ is not even, the statement remains correct in the sense that the signature of $E$
is zero.) I was told by Andrew Ranicki that this is \cite[Thm. 2.7]{LueRa1992}, to be used with rational coefficients; and that
in the case of a fiber bundle where base and fiber are closed manifolds, the work of W.~Meyer \cite{Meyer72} can be used.
A very transparent proof of this fact can be found in \cite{Rovi}.
The conclusion from this is that part (i) of proposition~\ref{prop-uv} has more radical variants: for example $\kappa_t(\sL_{n+k})$
comes from the cohomology of $B(\pi_0\haut_\partial(W_{g,1}))$, a.k.a.~cohomology of the
discrete group $\pi_0\haut_\partial(W_{g,1})$. From this point of view, it might seem that
section~\ref{sec-conhoto} with its emphasis on $\partial W_{g,1}$ is superfluous.
But there is part (ii) of proposition~\ref{prop-uv}, too. It is not easy to take $\partial W_{g,1}$ out of that.
Therefore section~\ref{sec-conhoto} could be essential after all.
\end{rem}

\begin{rem} In an earlier version of this article, one of the main results of \cite{BeBe} (their Cor.3)
was stated as a conjecture, and a few special cases were proved. These special cases were
used and are used in the proof of theorem~\ref{thm-gmain}.
Evidence for the conjecture came from tables \cite{McTague}.
\end{rem}

\begin{rem} \label{rem-KupKran} The proof of theorem~\ref{thm-gmain} used some upper bounds on the size
of the homotopy groups of $\haut_\fin(\config(\RR^{2n};r)_\QQ)$ for rather small $r$, and the close relationship of
$\haut_\fin(\config(\RR^{2n};r)_\QQ)$ to
the path-componentwise rationalization of $\haut_\fin(\config(\RR^{2n};r))$. The point was to show that the comparison map
from $\SOr(2n)$ to $\haut_\fin(\config(\RR^{2n};r)_\QQ)$ is rather far from inducing a surjection
in rational cohomology. These computations were based on the failure of homological
stability observed in section~\ref{subsec-stabi}. No serious attempt was made to understand
the homotopy type of $\haut_\fin(\config(\RR^{2n};r)_\QQ)$ as a whole. But Manuel Krannich and Alexander Kupers
pointed out (at a workshop held in Bonn, first week of September 2019) that a lot of information on the homotopy type of
$\haut_\fin(\config(\RR^{2n})_\QQ)$ is available thanks to the efforts of Fresse, Turchin and Willwacher \cite{FresseTurchinWillwacher}.
(This is written in the language of operads, so one has to use theorem~\ref{thm-operpest} to translate.
It relies on the formality of the little disk operads, sketched in \cite{KontsevichFormal} and
fully spelled out in \cite{LambrechtsVolic}, which has not been used here.)
Using \cite[Cor.5]{FresseTurchinWillwacher}, they said, one could hope to understand
the comparison map from $\SOr(2n)_\QQ$ to $\haut_\fin(\config(\RR^{2n};r)_\QQ)$ by writing it as a composition
\[ \SOr(2n) \lra \haut_\fin(\config(\RR^{2n})_\QQ) \lra \haut_\fin(\config(\RR^{2n};r)_\QQ)
\]
and examining the first arrow in the composition.
This might lead to a statement much sharper
than theorem~\ref{thm-gmain}. (The reviewer points out that instead of using \cite{FresseTurchinWillwacher} one could
use \cite{FresseWillwacher}.)
\end{rem}

\section{Detecting elements in rational homotopy groups} \label{sec-decomp}
\subsection{Euclidean bundles on spheres}
Reacting to an earlier draft of this article, Diarmuid Crowley asked me (in March 2016)
whether the nonzero Pontryagin classes $p_{n+k}\in H^{4n+4k}(B\TOP(2n);\QQ)$
found in theorem~\ref{thm-pmain} evaluate nontrivially on $\pi_{4n+4k}(B\TOP(2n))\otimes\QQ$.

\begin{prop} \label{prop-Crowleyquest} They do evaluate nontrivially on $\pi_{4n+4k}(B\TOP(2n))\otimes\QQ$.
\end{prop}

\proof We begin with some general and less general observations. Cohomolo\-gy is taken with rational
coefficients throughout.
\begin{itemize}
\item[(i)] Let $X$ be a simply connected based CW-space and $b\in H^s(X)$ where $s\ge 3$. Suppose that $\Omega^2b\in H^{s-2}(\Omega^2X)$
is nonzero. (Think of $b$ as a homotopy class of based maps from $X$ to an Eilenberg-MacLane space $K(\QQ,s)$. This justifies the
notation $\Omega^2b$ for a homotopy class of based maps from $\Omega^2X$ to $\Omega^2 K(\QQ,s)=K(\QQ,s-2)$.) \emph{Then $b$ evaluates
nontrivially on $\pi_s(X)$.} Proof: if $\Omega^2b$ is nonzero, then $\Omega b\in H^{s-1}(\Omega X)$ is indecomposable in $H^*(\Omega X)$.
Since $\Omega X$ is a connected $H$-space, it is rationally a product of Eilenberg-MacLane spaces (well-known)
and it follows that $\Omega b$ evaluates nontrivially on $\pi_{s-1}(\Omega X)$.
\item[(ii)] Let $(X,Y)$ be a based CW-pair and $b\in H^s(X,Y)$.
Let $Z$ be the homotopy fiber of the inclusion $Y\to X$. The functor taking a based pair $(X,Y)$ to
the homotopy fiber of $Y\hookrightarrow X$ has a left adjoint, $Z\mapsto(\cone(Z),Z)$; consequently there is a \emph{counit} map
from $(\cone(Z),Z)$ to $(X,Y)$, for $Z=\hofiber(Y\to X)$.  Let
\[  \varphi\co H^s(X,Y) \lra \tilde H^{s-1}(Z)\cong H^s(\cone(Z),Z) \]
be the homomorphism induced by that map. Suppose that $Z$ is simply connected, $s\ge 4$, and $\Omega^2\varphi(b)\in H^{s-3}(\Omega^2Z)$ is nonzero. \emph{Then $b$
evaluates nontrivially on $\pi_s(X,Y)$.} Proof: use (i) and use $\pi_s(X,Y)\cong \pi_{s-1}(Z)$.
\item[(iii)] The class $p_{n+k}\in H^{4n+4k}(B\TOP(2n))\cong H^{4n+4k}(B\STOP(2n))$ lifts to
\[ H^{4n+4k}(B\STOP(2n),BS\Or(2n)). \]
(We are assuming $k>0$, as in theorem~\ref{thm-gmain}.)
If a chosen lift evaluates nontrivially on
$\pi_{4n+4k}(B\STOP(2n),B\SOr(2n))$, then $p_{n+k}$ itself evaluates
nontrivially on $\pi_{4n+4k}(B\STOP(2n))$. Proof:
the inclusion-induced map
\[ \pi_{4n+4k}(B\STOP(2n))\to \pi_{4n+4k}(B\STOP(2n),B\SOr(2n)) \]
is a rational isomorphism.
\item[(iv)] Let $V$ be a topological manifold with boundary. Then the space of collars on $\partial V$
is contractible. (This space can be defined as the geometric realization of a simplicial set
where a $k$-simplex is an embedding $\Delta^k\times \partial V\times[0,1)\to \Delta^k\times V$
over $\Delta^k$, etc.)
\end{itemize}
Sketch proof of (iv): By Brown \cite{Brown62} and Connelly \cite{Connelly71} the space of collars on $\partial V$ is
nonempty. This makes it easy to reduce to the situation where $V\cong \partial V\times [0,1)$. In that case the space of collars
has a preferred base point and it is easy to produce a contraction of the space of collars
to that base point.

Next, we use the manifold formulation (theorem~\ref{thm-gmain}) with fiber bundle projection
$E \to M$. 
\begin{itemize}
\item[(v)] By construction, the map $E\to M$ has a section $s\co M\to E$ such that the restricted projection
$E\smin s(M) \lra M$ has the structure of a smooth fiber bundle (with noncompact fibers).
\item[(vi)] There exists a (topological) embedding $u\co \RR^{2n}\times M \to E$ over $M$ such that $s(x)=u(0,x)$
for all $x\in M$.
\end{itemize}
Note that (v) and (vi) together are equivalent to proposition~\ref{prop-concess}. Sketch proof of (vi):
Begin with a (topological) embedding
$v\co D^{2n} \times M\to E$ over $M$ such that $s(x)=v(z,x)$ for all $x\in M$,
and for a fixed $z\in \partial D^{2n}$. Use (iv) to choose a fiberwise collar on
$v(\partial D^{2n}\times M)$ in $E\smin v((D^{2n}\smin\partial D^{2n})\times M)$, total space
of a fiber bundle on $M$. Use that
fiberwise collar to extend $v$ to an embedding $u_1\co \RR^{2n}\times M \to E$ over $M$.
Define $u(y,x)=u_1(y+z,x)$.
\begin{itemize}
\item[(vii)] There is a commutative square
\[
\xymatrix@R=17pt{
H^*(B\STOP(2n)) \ar[r]^-{t^*} & H^*(E) \\
H^*(B\STOP(2n),B\SOr(2n)) \ar[u] \ar[d]^-{\Omega^{2n-1}\circ\varphi} \ar[r] & \ar@{<->}[d]_-\cong \ar[u] H^*(E,E\smin s(M))  \\
\tilde H^{*-2n}(\Omega^{2n-1}(\STOP(2n)/\SOr(2n))) \ar[r] & H^{*-2n}(M)
}
\]
\end{itemize}
The upper and middle horizontal arrows in (vii) are induced by the classifying map $t$
for the vertical tangent bundle of $E$. The isomorphism in the right-hand column
uses excision to replace $E$ by the image of the embedding $u\co \RR^{2n}\times M\to E$.
The lower half of the diagram is easier to understand if we replace cohomology
of pairs by reduced cohomology of the corresponding mapping cones (homotopy cofibers),
\[ \begin{array}{l}
 C_1= \cone [B\SOr(2n)\hookrightarrow B\STOP(2n)], \\
 C_2= \cone[E\smin s(M)\hookrightarrow E].
\end{array}
\]
There is a homotopy commutative square of based spaces
\[
\xymatrix@R=12pt{
\Omega^{2n}C_1 & \ar[l]_-{\Omega^{2n}\bar t}  \Omega^{2n} C_2 \\
\Omega^{2n-1}(\STOP(2n)/\SOr(2n)) \ar[u] & \ar[l] \ar[u] M_+
}
\]
where $\bar t\co C_2\to C_1$ is induced by $t$, the right-hand vertical arrow is adjoint to the
weak equivalence $S^{2n}\wedge M_+ \to C_2$\,, and the left-hand vertical arrow is $\Omega^{2n-1}$ of a
more obvious map $\STOP(2n)/\SOr(2n)\to \Omega C_1$. Now we can trade the $\Omega^{2n}$ prefix in the upper row
for an $S^{2n}\wedge$ prefix in the lower row.
This explains the commutative square of (reduced) cohomology groups.

So much for preparations. Choose $b\in H^{4n+4k}(B\STOP(2n),B\SOr(2n))$ which lifts the Pontryagin class
$p_{n+k}\in H^{4n+4k}(B\STOP(2n))$.
Then $b$ has a nonzero image in $H^{2n+4k}(\Omega^{2n-1}(\STOP(2n)/\SOr(2n)))$
since it has a nonzero image in $H^{4n+4k}(E)$. Using (ii), we deduce that $b$ evaluates nontrivially on the relative homotopy
group $\pi_{4n+4k}(B\STOP(2n),B\SOr(2n))$. Using (iii), we deduce that $p_{n+k}$ evaluates nontrivially on $\pi_{4n+4k}(B\STOP(2n))$.
\qed

\begin{expl} This is an example of a simply connected based space $X$ such that the
canonical graded homomorphism from the rational cohomology $H^*(X)$ modulo decomposables to $\hom(\pi_*(X),\QQ)$ is not injective.
(I learned this from Diarmuid Crowley.) Choose a non-torsion element
in $\pi_4(S^2\vee S^2)$ and represent it by a based map $f\co S^4\to S^2\vee S^2$. Such an element exists
by \cite{Hilton55} or \cite{MilnorFK}. Let $X$ be the mapping cone of $f$.
Then $H^5(X)\cong \QQ$ and any nonzero element of that group is indecomposable. The inclusion
$S^2\vee S^2\to X$ induces a homomorphism in $H^5$ which is zero, and a homomorphism in $\pi_5$
which is surjective by the long exact sequence of homotopy groups of the pair $(X,S^2\vee S^2)$.
Therefore every element of $H^5(X)$ evaluates trivially on $\pi_5(X)$.
\end{expl}

\subsection{Spaces of smooth automorphisms of disks}
The inclusion of $B\Or(2n)$ in $B\TOP(2n)$
clearly induces a surjection of the rational cohomology
rings. Therefore we have an isomorphism in rational cohomology
\[  H^*(B\TOP(2n),B\Or(2n))  \xrightarrow{\,\,\cong\,\,} \ker[\,H^*(B\TOP(2n))\to H^*(B\Or(n))\,]. \]
In particular $p_{n+k}\in H^{4n+4k}(B\TOP(2n))$ for $k>0$ can be viewed unambiguously as an element in
$H^{4n+4k}(B\TOP(2n),B\Or(2n))$.

\medskip
Let $(X,Y)$ be a pair of spaces, with base point $*\in Y$. Let $Z$ be the homotopy fiber (over $*$)
of the inclusion $Y\to X$. There is a counit map of based pairs
\[  (\cone(Z),Z)\lra (X,Y). \]
It induces a map $H^*(X,Y)\to H^*(\cone(Z),Z)\cong \tilde H^{*-1}(Z)$ in (rational) cohomology which we
denote $b\mapsto \varphi(b)$, as in the proof of proposition~\ref{prop-Crowleyquest}.
In particular we may take $X=B\TOP(2n)$, $Y=B\Or(2n)$ and $r=2n+1$. Then $Z\simeq \TOP(2n)/\Or(2n)$
and $\Omega^rZ=\Omega^{2n+1}Z\simeq \Omega^{2n+1}(\TOP(2n)/\Or(2n))\simeq \diff_\partial(D^{2n})$
by smoothing theory (remark~\ref{rem-Morlet} and remark~\ref{rem-moreMorlet} below).
For $k>0$ we can view $p_{n+k}$ as an element of $H^{4n+4k}(B\TOP(2n),B\Or(2n))$
and we obtain rational cohomology classes
\[ \Omega^{2n+1}\varphi(p_{n+k})\in H^{2n-2+4k}(\diff_\partial(D^{2n})). \]
The following is an obvious consequence of theorem~\ref{thm-pmain} and proposition~\ref{prop-Crowleyquest}.

\begin{prop} \label{prop-diskdiff} Under the conditions of theorem~\ref{thm-pmain} on $n$ and $k$, the class
$\Omega^{2n+1}\varphi(p_{n+k})\in H^{2n-2+4k}(\diff_\partial(D^{2n}))$ is nonzero. Indeed it
evaluates nontrivially on the rational homotopy group $\pi_{2n-2+4k}(\diff_\partial(D^{2n}))\otimes\QQ$. \qed
\end{prop}

\begin{rem} \label{rem-moreMorlet} Morlet's weak equivalence
\[ \diff_\partial(D^m) \xrightarrow{\,\simeq\,} \Omega^{m+1}(\TOP(m)/\Or(m)) \]
(for $m\ne 4$), already mentioned in remark~\ref{rem-Morlet},
is the map of vertical homotopy fibers in the commutative square
\[
\xymatrix@R=15pt{
\diff_\partial(D^m) \ar[d] \ar[r] & \Omega^m\Or(m) \ar[d] \\
\homeo_\partial(D^{m}) \ar[r] & \Omega^m\TOP(m)
}
\]
The horizontal arrows are obtained by associating to a diffeomorphism or
homeomorphism its derivative. (The notation is deliberately careless about the distinction between $\Or(m)$ and the general
linear group $\GL(m)$.) Note that $\homeo_\partial(D^{m})$ is contractible, by
the Alexander trick. For the present purposes, the horizontal arrows can be defined
as homomorphisms of simplicial groups. In particular $\TOP(m)$ could be interpreted \cite{Kister}
as the simplicial group of homeomorphism germs from $(\RR^m,0)$ to $(\RR^m,0)$. ---
Burghelea and Lashof \cite[Thm.4.4.(a)]{BurghLash} write
\[ \diff_\partial(D^m) \xrightarrow{\,\simeq\,} \Omega^{m+1}(\PL(m)/\Or(m)), \]
presumably following Morlet. This formulation has the advantage that $m=4$ is \emph{not} an exception. As
explained in \cite[Essay V,\S5.0]{KirbySiebenmann} or \cite[Thm.4.6]{BurghLash},
the inclusion $\Omega^{m+1}(\PL(m)/\Or(m))\to \Omega^{m+1}(\TOP(m)/\Or(m))$ is a weak equivalence for $m\ne 4$.
\end{rem}

\begin{appendices}
\section{Composition and descent} \label{sec-compodesc}
This appendix is an abstract justification for a step taken earlier in passing from
diagram~\eqref{eqn-mainbdryspec2} to diagram~\eqref{eqn-mainbdryspec2coord}. Readers who found that step self-explanatory may prefer to skip the appendix.

The word \emph{descent} (from sheaf theory) refers to situations where global data can be recovered from
localized data together with gluing instructions. Our aim here is to describe in homotopical terms something like composition of
maps-subject-to-conditions, provided the conditions allow descent.

\subsection{Simplicial sets with monoid action} \label{subsec-simpmon} Let $\sK$ be a simplicial monoid.
The category of simplicial sets with a right action of $\sK$ comes with a forgetful functor
to the category $\sSet$ of simplicial sets. That functor has a left adjoint given by $X\mapsto X\times \sK$.
The adjunction can be used to export the standard model structure on the category of simplicial
sets (where fibrations are Kan fibrations) to the category of simplicial sets with right action of $\sK$.

\begin{prop} \label{prop-prostration}
The category of simplicial sets with a right action of $\sK$ admits a (unique) cofibrantly generated
model structure where a morphism $X\to Y$ is a weak equivalence, respectively fibration,
if and only if the underlying map of simplicial sets is a weak equivalence, respectively fibration.
The generating cofibrations and acyclic cofibrations can be obtained by taking the standard
selection of generating cofibrations and acyclic cofibrations in the model category of simplicial sets
and applying the functor $-\times \sK$ to them.

This model category is a simplicial model category:
the tensoring over $\sSet$ is given by $X\otimes Y:= X\times Y$, where $X$ is an ordinary
simplicial set and $Y$ is a simplicial set with right $\sK$-action.
\end{prop}

\proof This is an application of \cite[Thm.3.6]{GoerssSchem}. (There is something nontrivial to verify, but it is
left to the reader.) Alternative: \cite[Thm.4.1]{SchwedeShipley}. \qed

\begin{expl} \label{expl-resol} 
Let $X$ be a simplicial set equipped with a right action of $\sK$.
Let $\bar X$ be the bar resolution of $X$, better known as a two-sided bar construction
$B(X,\sK,\sK)$; see \cite[\S4.2]{Riehl} or \cite{HollenderVogt}. This is the geometric
realization of the simplicial object $[r] \mapsto X\times \sK^{r+1}$,
with face operators $d_i$ given by
\[  (x,y_0,\dots,y_r) \mapsto \left\{ \begin{array}{rl} (xy_0,y_1,\dots,y_r) & \textup{if }i=0 \\
(x,y_0,...,y_{i-1}y_i,y_{i+1},\dots,y_r) & \textup{if }0<i\le r
\end{array} \right.
\]
and degeneracy operators $s_i$ given by
$(x,y_0,\dots,y_r) \mapsto (x,y_0,\dots,y_{i-1},1,y_i,\dots,y_r)$.
The geometric realization is to be carried out within the category of simplicial sets. Consequently
$\bar X$ comes with a natural augmentation $\bar X\to X$ which is a weak equivalence,
and with a natural right action of $\sK$ making the augmentation $\bar X\to X$ into a $\sK$-map.
It is easy to see that $\bar X$ is cofibrant in the model category structure of proposition~\ref{prop-prostration};
therefore $\bar X\to X$ qualifies as a functorial cofibrant replacement.

Note in passing that the geometric realization is an instance of a homotopy colimit (for a functor from $\Delta^\op$
to spaces).
\end{expl}

\begin{rem} \label{rem-derivedmaps} In a model category $\sC$ according to \cite{Hovey}, \cite[Def.7.1.3]{Hirschhorn} there exist \emph{functorial} fibrant and cofibrant replacements for objects $c$. It follows that we have
functorial fibrant-cofibrant replacements $c''$ for objects $c$:
\[
\xymatrix@M=10pt@R=12pt@C=15pt{
& \textup{terminal object} \\
c & \ar[l]_-\simeq c'  \ar[u] \ar[r]_-{\textup{cofib.}}^-\simeq & c'' \ar[ul]_-{\textup{fib.}} \\
&    \ar[u]^-{\textup{cofib.}} \textup{initial object}
}
\]
(\emph{Functorial} means that the vertices in the diagram are functors of the variable $c$
and the arrows are natural transformations.)
Suppose that $\sC$ is a simplicial model category. For objects $c$ and $d$ of $\sC$,
the ordinary mapping object $\map_{\sC}(c,d)$ is the simplicial set $\mor_{\sC,\bullet}(c,d)$,
where a $j$-simplex is a morphism $\Delta^j\otimes c\to d$. The \emph{right derived} mapping object could be
defined as $\rmap_\sC(c,d)= \map_\sC(c'',d'')$. 

An explicit definition of $\rmap_\sC(c,d)$ for all $c,d$ can be a straitjacket, though.
Instead we can take the view that $\map_\sC(c,d)$ is a valid interpretation of \emph{right derived mapping space}
from $c$ to $d$ whenever $c$ is cofibrant and $d$ is fibrant.

In the case where $\sC$ is the category of simplicial sets with right action of the simplicial monoid $\sK$, we may write
$\map^\sK(X,Y)$ instead of $\map_\sC(X,Y)$.
\end{rem}

\begin{expl} \label{expl-framed} Suppose that the simplicial monoid $\sK$ is fibrant as a simplicial set.
Let $U$ and $V$ be simplicial sets, and suppose
that $V$ is fibrant as such. Then $U\times\sK$ is cofibrant and $V\times\sK$ is fibrant-cofibrant (in the model
category of simplicial sets with right $\sK$-action). Therefore $\map^\sK(U\times\sK,V\times\sK)$ can be
regarded as a right derived mapping space. Clearly we have
\[  \map^\sK(U\times\sK,V\times\sK)\cong \map(U,V\times\sK)\cong \map(U,V)\times \map(U,\sK). \]
In the case where $U=V$, we also obtain the structure of the monoid of endomorphisms of $V\times\sK$:
\[  \map^\sK(V\times\sK,V\times\sK)\cong \map(V,V)\ltimes \map(V,\sK). \]
Here $\map(V,V)$ is the monoid of endomorphisms of $V$ as a (fibrant) simplicial set,
$\map(V,\sK)$ is a monoid by pointwise multiplication (in $\sK$), and $\map(V,V)$ acts on
the right of $\map(V,\sK)$ by composition. The semidirect product is as in definition~\ref{defn-semidirect}.
\end{expl}

Let $\sJ$ be the union of the path components of $\sK$ which, as elements of the discrete monoid $\pi_0\sK$, are invertible.
Then $\sJ$ is a submonoid of $\sK$ and, as a simplicial monoid in its own right, it is group-like. (We are not assuming that $\sK$ is fibrant.) There are Quillen adjunctions
\[
\xymatrix@C=12pt@M=7pt@R=13pt{
\textup{simplicial sets with right action of $\sK$} \ar@<1ex>[d]^-\rho \\
\textup{simplicial sets with right action of $\sJ$} \ar@<1ex>[u]^-\eta \ar@<1ex>[d]^-\beta \\
\textup{simplicial sets over $B\sJ$} \ar@<1ex>[u]^-\pi \ar@<0.7ex>[r] & \ar@<0.7ex>[l]
{\begin{array}{l} \textup{simplicial sets over a fibrant} \\ \textup{replacement $B^\sharp\sJ$ of $B\sJ$} \end{array}}
}
\]
($\rho$ for restriction, $\eta$ for extension, $\beta$ for Borel construction, $\pi$ for
pullback). In more detail, $\rho$ is given by restricting $\sK$-actions to $\sJ$,
and $\eta$ is given by $X\mapsto X\times_\sJ\sK$.
The functor $\pi$ takes $X$ with reference map to $B\sJ=B(*,\sJ,*)$ to
the pullback of
\[ X \to B(*,\sJ,*)\leftarrow B(*,\sJ,\sJ) \]
where we use the two-sided bar construction as in \cite{HollenderVogt}, \cite[\S4.2]{Riehl}.
The best way to describe $\beta$ is probably to say nothing about it except that it is right adjoint to $\pi$.
(If this is not enough: for a simplicial set $X$ with right action of $\sJ$, a $k$-simplex in $\beta(X)$
is a pair $(u,v)$ where $u\co\Delta[k]\to B\sJ$ is a map of simplicial sets and $v$ is a morphism from $\pi(u)$ to $X$.)
In the lower row, we are assuming a weak equivalence
and cofibration $B\sJ\to B^\sharp\sJ$. The horizontal arrow from left to right is given by composition with that
map; the other is given by pullback along that map. (The category of simplicial sets
over a fixed simplicial sets inherits a preferred model structure from the category of simplicial
sets; see \cite[Expl.1.7]{GoerssSchem}. In this model structure, a morphism is a cofibration/fibration/weak equivalence
if and only if the underlying map of simplicial sets is an injection/Kan fibration/weak equivalence.)

The functor $\pi$ preserves cofibrations and trivial
cofibrations (a.k.a.~acyclic cofibrations). The functor $\rho$ preserves fibrations and trivial fibrations
(a.k.a.~acyclic fibrations). It follows \cite[Prop.8.5.3]{Hirschhorn} that $\pi$ and $\beta$ constitute
a Quillen pair, and also that $\rho$ and $\eta$ constitute a Quillen pair. Similarly the two horizontal
arrows constitute a Quillen pair. Moreover, the functors $\pi$ and $\beta$ constitute a Quillen equivalence, and
the two horizontal arrows also constitute a Quillen equivalence.

\medskip
We come to a description of the
derived morphism spaces in the category of simplicial sets over a fixed base $B$ which is itself fibrant.
(The case that one should have in mind is $B=B^\sharp\sJ$.)
Let $p\co X\to B$ and $q\co Y\to B$ be two
objects in the model category of simplicial sets over $B$. Suppose that $q$ is a fibrant
object, i.e., the map $q$ is a fibration.
(The two objects are automatically cofibrant.)

\begin{prop} Then the mapping object
$\map(X\to B,\,Y\to B)$ is identified with the fiber over $p$ of the map
\[
\map(X,Y) \lra \map(X,B)~;~ f\mapsto qf
\]
which is a fibration in its own right. Moreover $Y$ is fibrant as a simplicial set in its own right. \qed
\end{prop}
Therefore we can also say informally that $\rmap(X\to B,\,Y\to B)$ is identified with the
homotopy fiber over the $0$-simplex $p$ of the map from $\rmap(X,Y)$ to $\rmap(X,B)$ given by composition with $q$.
This description of right derived mapping objects (in the model category of simplicial sets over $B$)
extends to a description of composition. The details are left to the reader.

\smallskip
Let $X$ and $Y$ be simplicial sets with right action of $\sJ$. Suppose that both are
cofibrant as such. Abbreviate $X^e:=X\times_\sJ\sK$
and $Y^e:=Y\times_\sJ\sK$. We ask how $\rmap^\sJ(X,Y)$ is related to
$\rmap^\sK(X^e,Y^e)$.

As a simplicial set, $X^e$ is the disjoint union (coproduct) of $X$ alias $X\times_\sJ\sJ$ and something else
which can be described as $X\times_\sJ(\sK\smin\sJ)$. Consequently we may write $\pi_0X\subset \pi_0X^e$.
An element of $\pi_0\rmap^\sK(X^e,Y^e)$ induces a map $\pi_0X^e\to \pi_0 Y^e$. This will always take
the complement of $\pi_0X$ to the complement of $\pi_0Y$. If it takes
$\pi_0X$ to $\pi_0Y$, we call it \emph{tame}.
(Example: all invertible elements of $\pi_0\rmap^\sK(X^e,Y^e)$ are tame.)

\begin{prop} The comparison map $\rmap^\sJ(X,Y)\to \rmap^\sK(X^e,Y^e)$ restricts to a weak equivalence
of $\rmap^\sJ(X,Y)$ with a union of path components of $\rmap^\sK(X^e,Y^e)$, namely, the tame ones.
\end{prop}
The proof is left as an exercise. It is probably a good idea to begin with the case where
$X$ is free on one generator, in other words $X=\sJ$.

\subsection{Morphism transport}
As in \cite{BoavidaWeissAbdn}, let $\man$ be the following category, enriched over simplicial sets.
The objects are the smooth $m$-manifolds (without boundary), for a fixed $m$; the simplicial set $\mor_{\man,\bullet}(K,L)$ of morphisms from $K$ to
$L$ is the simplicial set of smooth embeddings from $K$ to $L$. More precisely, a $j$-simplex in $\mor_{\man,\bullet}(K,L)$
is a smooth embedding of $\Delta^j\times K$ in $\Delta^j\times L$ which respects the projection to $\Delta^j$.  

\begin{defn} Let $E$ be a functor from $\man$ to
a simplicial model category $\sD$ which respects the enrichments over $\sSet$.
The functor is \emph{weakly additive} if
\begin{itemize}
\item[-] it takes $\emptyset$ to the initial object of $\sD$ (up to weak equivalence);
\item[-] it takes union-intersection squares
\[
\xymatrix@R=12pt@C=20pt{
U\cap V \ar[r] \ar[d] & U \ar[d] \\
V \ar[r] & U\cup V
}
\]
in $\man$ to homotopy pushout squares in $\sD$;
\item[-] it takes monotone unions $M=\bigcup_i M_i$
(where $M_i\subset M_{i+1}$ for $i=0,1,2,\dots$ and each $M_i$ is open in $M$) in $\man$
to the corresponding sequential homotopy colimits. (In other words, the comparison map from $\hocolim_i~E(M_i)$ to $E(M)$
is a weak equivalence.)
\end{itemize}
\end{defn}

\smallskip
We fix a simplicial model category $\sC$ and a functor $F\co\man\to \sC$
respecting the simplicial enrichments. (The functor $F$ in example~\ref{expl-endoconfigloc} below is the motivating instance.)
Here is a first condition which we impose on $F$.
\begin{itemize}
\item[(a)] For every $M$ in $\man$, the object $F(M)$ in $\sC$ is fibrant and cofibrant.
\end{itemize}

\begin{defn} Let $\man[F]$ be the following category enriched over simplicial sets. The objects are the
objects of $\man$ and the simplicial set of morphisms from $K$ to $L$ is
$\mor_{\sC,\bullet}(F(K),F(L))$.
There is an enrichment-preserving functor $\man\to \man[F]$ which is the identity on objects. On morphisms
(of any degree) it is defined by $g\mapsto F(g)$.
\end{defn}

The goal is to understand $\man[F]$ in homotopy theoretic terms.
But we impose more conditions on $F$. The most convenient
all-in-one condition for our purposes is the following.
\begin{itemize}
\item[(b)] $F$ is weakly additive and the functor $M\mapsto \mor_{\sC,\bullet}(F(\RR^m),F(M))$
from $\man$ to $\sSet$ is also weakly additive.
\end{itemize}
We assume (a) and (b) for now.

Let $\sH$ be the simplicial monoid of endomorphisms
of the object $\RR^m$ in $\man$. Let $X_M= \mor_{\man,\bullet}(\RR^m,M)$,
so that $\sH$ acts on the right of $X_M$. Note that
$\sH$ is group-like. Indeed it is weakly equivalent to the singular simplicial set of the orthogonal
group $\Or(m)$, by an obvious inclusion. (In the same spirit: $X_M$ is weakly equivalent to the
singular simplicial set of the total space of the frame bundle of $M$.)
Similarly let $\sK$ be the simplicial monoid of endomorphisms
of $F(\RR^m)$ in $\sC$ and let
\[ Y_M:= \mor_{\sC,\bullet}(F(\RR^m),F(M)) \]
so that $\sK$ acts on the right of $Y_M$. There is a monoid homomorphism $\sH\to \sK$ induced by $F$.

\begin{prop} \label{prop-breakup} For $M$ in $\man$ there is a natural weak equivalence
\[ B(X_M,\sH,\sK)\lra Y_M \]
of simplicial sets with right
action of $\sK$.
\end{prop}
\proof The natural map as such is obvious, or if not, can be written as a composition
\[  B(X_M,\sH,\sK) \to B(Y_M,\sK,\sK) \xrightarrow{\simeq} Y_M\,. \]
It is a weak equivalence for $M=\RR^m$; in that case it takes the form of
a standard resolution map $B(\sH,\sH,\sK)\to \sK$. Therefore it suffices to show that the two functors
\[  M\mapsto B(X_M,\sH,\sK), \qquad M\mapsto Y_M  \]
are weakly additive. For the first one this is clear and for the second one it is part of condition (b). \qed

\smallskip
For $M$ and $N$ in $\man$, there is a map of
type \emph{composition},
\begin{equation} \label{eqn-descent} \mor_{\sC,\bullet}(F(M),F(N))  \lra \map^\sK(Y_M,Y_N).
\end{equation}
Now that we are dealing with simplicial sets equipped with a right action of $\sK$, we adopt the
(model category) conventions of section~\ref{subsec-simpmon}.

\begin{prop} \label{prop-compodescent} For all $M$ and $N$ in $\man$, the composition of \eqref{eqn-descent} with the
resolving map $\map^\sK(Y_M,Y_N)\to \rmap^\sK(Y_M,Y_N)$
is a weak equivalence.
\end{prop}
\proof Fix $N$ throughout. For $M=\RR^m$, we get $Y_M=\sK$. The map~\eqref{eqn-descent} is an isomorphism in this case. The
resolving map from $\map^\sK(\sK,Y_N)\cong Y_N$ to $\rmap^\sK(\sK,Y_N)$ is a weak equivalence.
Therefore it is enough to know that the functors $F$
and $M\mapsto Y_M$ on $\man$ are both weakly additive. This is exactly the content of condition (b).  \qed

\begin{cor} \label{cor-endoandframed} Suppose that $M$ in $\man$ admits a framing, i.e., a trivialization of the tangent vector
bundle. Let $sM$ denote the singular simplicial set of $M$.
Then the simplicial monoid $\map(F(M),F(M))$ is weakly equivalent to a semidirect product $\rmap(sM,sM) \ltimes \rmap(sM,\sK)$.
\end{cor}

\proof The framing implies that $X_M$ in proposition~\ref{prop-breakup} is weakly equivalent, as a simplicial set
with right action of $\sH$, to $sM\times\sH$. Therefore by the same proposition, $Y_M$ is weakly equivalent, as a simplicial
set with right action of $\sK$, to $sM\times \sK$. Then by proposition~\ref{prop-compodescent},
the simplicial monoid $\map(F(M),F(M))$ is weakly equivalent to
\[ \rmap^{\sK}(sM\times\sK,sM\times\sK) \]
and that in turn (by example~\ref{expl-framed}) is weakly equivalent to a semidirect product,
$\rmap(sM,sM) \ltimes \rmap(sM,\sK)$. \qed

\begin{cor} Suppose that $M$ in $\man$ admits a framing and let $U\subset M$ be an open subset.
Let $\rmap_{sU}(sM,sM)$ denote the simplicial monoid of maps $sM\to sM$ which extend the identity on $sU$.
Then
\[  \hofiber[\,\map(F(M),F(M)) \lra \map(F(U),F(M))\,] \]
is weakly equivalent to a semidirect product $\rmap_{sU}(sM,sM) \ltimes \rmap_*(sM/sU,\sK)$. \emph{(The homotopy
fiber is taken over the $0$-simplex in $\map(F(U),F(M))$ determined by the inclusion $U\to M$.)} \qed
\end{cor}

\begin{expl} \label{expl-endoconfigloc} Let $\sD$ be the category of simplicial spaces (where \emph{space} should be read as
\emph{simplicial set}). Let $\sC$ be the category of simplicial spaces over $N\fin$, the nerve of $\fin$. As before, $N\fin$ is viewed as a simplicial space which happens to be discrete in every degree. Let $F$ be the functor from $\man$ to $\sC$
taking $M$ in $\man$ to $\config^\loc(M;k)$ for some fixed $k$; here $k$ can be a positive integer or $\infty$.
For the moment we do not ask whether $F$ satisfies condition (a). To make up for this negligence we verify condition (b)
for $F$ in a formulation where $\mor_\sC(-,-)$ is replaced by $\rmap_\fin(-,-)$.

Let $F^\flat(M)\subset F(M)$ be the simplicial subspace of $F(M)$ projecting to the $0$-simplex $\uli 1$ of $N\fin$
and its degeneracies. Up to weak equivalence, $F^\flat(M)$ is a
constant simplicial space, weakly equivalent in degree $0$ to the underlying space of $M$.
The inclusion $F^\flat \to F$ can be viewed as a natural transformaton of functors from $\man$ to
$\sC$. There is also a
``forgetful'' natural transformation $\eta\co F\to F^\flat$ (but now of functors from $\man$ to $\sD$)
such that the composition $F^\flat\to F\to F^\flat$ is a
natural weak equivalence.

For any morphism $e\co L\to M$ in $\man$, the square
\[
\xymatrix@R=20pt{
F(L) \ar[r]^-{e_*} \ar[d]_\eta &  F(M) \ar[d]^-\eta \\
F^\flat(L)  \ar[r]^-{e_*} &  F^\flat(M)
}
\]
is homotopy cartesian in $\sD$ by \cite[Prop.4.2]{BoavidaWeiss}. Since $F^\flat$ is weakly additive
it follows, with the second Mather cube theorem \cite[Thm.25]{Mather}, that $F$ is weakly additive (as a functor from
$\man$ to $\sD$ or as a functor from $\man$ to $\sC$). The same observation also helps us to show that the functor
$M \mapsto \rmap_\fin(F(\RR^m),F(M))$ is weakly additive.
Namely, for any morphism $e\co L\to M$ in $\man$, the square of spaces
\[
\xymatrix@R=20pt{
\rmap_\fin(F(\RR^m),F(L)) \ar[r]^-{e_*} \ar[d]_{\eta_*} &  \rmap_\fin(F(\RR^m),F(M)) \ar[d]^-{\eta_*} \\
\rmap(F(\RR^m),F^\flat(L))  \ar[r]^-{e_*} & \rmap(F(\RR^m),F^\flat(M))
}
\]
is homotopy cartesian. Therefore it suffices to show that the functor from $\man$ to spaces taking $M$ in $\man$ to the
derived mapping space $\rmap(F(\RR^m),F^\flat(M))$ is weakly additive. (This a space of derived maps
in $\sD$.) Since $F^\flat(M)$ is weakly equivalent to a constant simplicial space, the map
\[ \rmap(F(\RR^m),F^\flat(M)) \to \rmap(\,\hocolim_{\Delta^\op}\,F(\RR^m), \hocolim_{\Delta^\op}\,F^\flat(M))  \]
is a weak equivalence and $\hocolim_{\Delta^\op}~F^\flat(M)~\simeq M$.
Moreover
\[ \hocolim_{\Delta^\op}~F(\RR^m) ~~\simeq~~ *  \]
because $F(\RR^m)$ is a Segal space admitting a weakly terminal object. It follows that there is a zig-zag of natural
weak equivalences relating $M\to \rmap(F(\RR^m),F^\flat(M))$ to the functor taking $M$ in $\man$
to its underlying space. That functor is clearly weakly additive, and the conclusion is that
$M\mapsto \rmap_\fin(F(\RR^m),F(M))$ is also weakly additive.

Finally we can enforce condition (a) by applying some functorial fibrant-cofibrant replacement to (the values of) $F$.
This modification of $F$ will then satisfy conditions (a) and (b) as written.
\end{expl}

\end{appendices}

\end{document}